\documentclass[12pt, a4paper, parskip=half, abstracton]{scrartcl}

\usepackage{array}
\usepackage{marginnote}
\usepackage{xcolor}
\usepackage{amscd,amssymb,amsfonts,amsmath,latexsym,amsthm}
\usepackage{hyperref}
\usepackage[all,cmtip]{xy}
\usepackage{mathrsfs}
\usepackage{graphicx}
\usepackage{bm} 
\usepackage[nameinlink, capitalise]{cleveref}

\usepackage{etoolbox}

\def\hB{\hspace*{\fill}$\qed$}

\usepackage[nottoc]{tocbibind} 

\newcommand{\Fun}{{\mathbf{Fun}}}

\newcommand{\cC}{{\mathcal{C}}}

\newcommand{\rk}{{\mathtt{rk}}}

\newcommand{\Ree}{{\mathtt{Re}}}

\newcommand{\diag}{{\mathtt{diag}}}

\newcommand{\cT}{{\mathcal{T}}}

\newcommand{\Proj}{{\mathtt Proj}}

\newcommand{\End}{{\mathtt{End}}}

\newcommand{\cB}{{\mathcal{B}}}

\newcommand{\Ind}{{\mathtt{Ind}}}

\newcommand{\ev}{{\mathtt{ev}}}

\newcommand{\id}{{\mathtt{id}}}
\newcommand{\Ab}{{\mathbf{Ab}}}

\newcommand{\cA}{{\mathcal{A}}}

\newcommand{\cM}{{\mathcal{M}}}

\renewcommand{\ker}{{\mathtt{ker}}}
\newcommand{\nat}{{\mathbb{N}}}
\newcommand{\Z}{{\mathbb{Z}}}

\newcommand{\R}{{\mathbb{R}}}
\newcommand{\C}{{\mathbb{C}}}

\newcommand{\Hom}{{\mathtt{Hom}}}

\renewcommand{\lim}{\operatorname*{\mathtt{lim}}}
\newcommand{\colim}{\operatorname*{\mathtt{colim}}}

\newcommand{\Sp}{\mathbf{Sp}}

\usepackage{slashed}
\usepackage[utf8]{inputenc}
\usepackage{microtype}
\usepackage[english]{babel}
\usepackage{mathtools}
\usepackage{bm}
\usepackage{esvect}

\title{Group completion of representation spaces and E-theory}
\author{
Ulrich Bunke\thanks{Fakult{\"a}t f{\"u}r Mathematik,
Universit{\"a}t Regensburg,
93040 Regensburg,
ulrich.bunke@mathematik.uni-regensburg.de} 
}

\numberwithin{equation}{section}
\newtheorem{theorem}{Theorem}[section] 
\newtheorem{prop}[theorem]{Proposition}
\newtheorem{lem}[theorem]{Lemma}

\newtheorem{ddd}[theorem]{Definition}
\newtheorem{kor}[theorem]{Corollary}

\theoremstyle{remark}
\theoremstyle{definition}
\newtheorem{prob}[theorem]{Problem}

\newtheorem{ex}[theorem]{Example}
\newtheorem{rem}[theorem]{Remark}

\newcommand{\Ani}{\mathbf{Ani}}

 \newcommand{\group}{\mathrm{group}}

\newcommand{\EE}{\mathrm{E}}
\newcommand{\sepa}{\mathrm{sep}}

\newcommand{\topp}{\mathrm{top}}

\newcommand{\nCalg}{C^{*}\mathbf{Alg}^{\mathrm{nu}}}

\newcommand{\ho}{\mathrm{ho}}

\newcommand{\cR}{\mathcal{R}}

\newcommand{\Pro}{\mathbf{Pro}}

\newcommand{\CMon}{\mathbf{CMon}}
\newcommand{\CGrp}{\mathbf{CGrp}}

\newcommand{\Fib}{{\mathrm{Fib}}}

\newcommand{\incl}{\mathrm{incl}}

\newcommand{\cP}{\mathcal{P}}

\renewcommand{\Proj}{\mathbf{Proj}}

\newcommand{\cZ}{{\mathcal{Z}}}

\newcommand{\cO}{{\mathcal{O}}}
\newcommand{\cU}{{\mathcal{U}}}

\newcommand{\Calg}{{\mathbf{C}^{\ast}\mathbf{Alg}}}

\newcommand{\op}{\mathrm{op}}

\newcommand{\add}{\mathrm{add}}

\newcommand{\KK}{\mathrm{KK}}

\newcommand{\lex}{\mathrm{lex}}

\newcommand{\ac}{\mathrm{ac}}

\newcommand{\bbK}{\mathbb{K}}
\newcommand{\red}{\mathrm{red}}
\newcommand{\fred}{\mathrm{r}}
\newcommand{\sfred}{\mathrm{sr}}

\renewcommand{\Proj}{\mathrm{Proj}}
  
\begin{document}
\maketitle

\begin{abstract}
{We study group completions of mapping anima between separable
$C^*$-algebras. For each separable source algebra $D$, we
construct a functor $\cR_D$ from separable $C^*$-algebras to
commutative groups in anima by group-completing stabilized mapping
anima, applying stabilized reduction, and resolving the source along
a shape system. For nuclear quasidiagonal $D$, and more generally
for sources in a class closed under countable filtered colimits,
finite products and retracts, we prove that $\cR_D$ is homological
and that the natural comparison}
\[
{\cR_D(B)\longrightarrow\Omega^\infty E(D,B)}
\]
{admits a natural section and induces an isomorphism}
\[
{(\pi_0\cR_D(B))_{\mathrm{Hd}}
\xrightarrow{\cong}E_0(D,B)_{\mathrm{Hd}}}
\]
{on Hausdorff quotients of component groups equipped with their
shape topologies.}

{For a countable discrete group $G$, we identify Willett's group
$QR(G)$ of quasi-representations with
$\bigl(\pi_0\cR_{C^*_{\max}(G)}(\C)\bigr)_{\mathrm{Hd}}$.
This gives a homotopy-theoretic interpretation of
quasi-representations and their comparison with $E$-theory.}
\end{abstract}

\tableofcontents

\section{Overview}
\label{sec:overview}
\begingroup
\theoremstyle{plain}
\newtheorem{overviewtheorem}{Theorem}
\renewcommand{\theoverviewtheorem}{\Alph{overviewtheorem}}

In this paper we  work with the category $\nCalg_{\sepa}$ of separable
$C^*$-algebras and possibly nonunital homomorphisms. 
We study the relation between group completions of representation
spaces and mapping spaces of bivariant $E$-theory. One motivation is to  interpret Willett's work \cite{willett} on quasi-representations
and controlled $KK$-theory  in the framework of non-commutative homotopy theory. 

Let $D$ and $B$ be in $\nCalg_{\sepa}$.
A compact representation of  $D$ on the standard Hilbert $B$-module $H_{B}:=B\otimes H$
is a homomorphism $\phi:D\to B\otimes \bbK$. A homotopy between compact representations is a path
of homomorphisms  $\phi_{u}:D\to B\otimes \bbK$, $u\in [0,1]$, such that $u\mapsto \phi_{u}(d)$ is norm-continuous for all $d$ in $D$.
Given two compact  representations $\phi,\psi$ we define their sum by
\[
{\phi\oplus \psi:D\xrightarrow{(\phi, \psi)} (B\otimes\bbK)\oplus(B\otimes\bbK)\xrightarrow{\diag}B\otimes {M}_{2}(\bbK)\cong B\otimes \bbK\ .}
\]
This operation  equips  the
set of   homotopy classes of compact representations with  a commutative monoid structure. We are interested in its group completion.

Equipping $\Hom_{\nCalg_{\sepa}}(D,B\otimes \bbK)$ with the point-norm topology we obtain a topological space of representations. The block sum operation refines it to an $E_{\infty}$-space. This $E_{\infty}$-space presents an object  ${\cM}_{D}(B)$ in the $\infty$-category  $\CMon(\Ani)$ of commutative monoids in anima. 
The monoid of homotopy classes of compact representations is then  recovered  by taking the monoid of connected components $\pi_{0}{\cM}_{D}(B)$.   

 Even if one is interested only in the group completion $(\pi_{0}{\cM}_{D}(B ))^{\group}$ it turns out to be useful to consider
the   representation space as a homotopy theoretic object. 
The group completion of 
 commutative monoids in anima is the left adjoint in an adjunction
\begin{equation}\label{gerfwregwerfwewerf}
 (-)^{\group}:\CMon(\Ani)\leftrightarrows\CGrp(\Ani):\incl
\end{equation}\cite[Corollary~4.4 and Remark~4.5]{gepner-groth-nikolaus}.
Since group completion commutes with taking connected components we recover the group completion of the monoid of components by 
\[
\pi_0({\cM}_D(B)^{\group})\cong(\pi_0 {\cM}_D(B))^{\group}.
\]
The central question of the present paper is to understand   the   representation group
${\cM}_{D}(B )^{\group}$ in $\CGrp(\Ani)$. It turns out to be useful  not to consider 
${\cM}_{D}(B)$ as an isolated object, but rather as a functor in the $B$-variable, and partly, also in the $D$-variable.

Considering group completions of commutative monoids in anima fits into the general philosophy of $K$-theory.
Applying group completion to the underlying symmetric monoidal
groupoid of finitely generated projective modules over a ring
produces the infinite loop space of its algebraic $K$-theory
spectrum. Similarly, the topological groupoid of vector bundles over a compact Hausdorff
space with direct sum gives the infinite loop space of the topological $K$-theory spectrum of the space
\cite[Section~8]{gepner-groth-nikolaus};
\cite[Theorem~IV.7.1]{weibel-kbook};
\cite[Sections~2 and 4]{segal-cohomology}.
For $D=\C$ we can interpret $\Hom_{\nCalg}(\C,B\otimes \bbK)$ as the space of projections in $B\otimes \bbK$, interpreting a compact representation $\phi$ as the projection $\phi(1)$. The
commutative monoid
${\cM}_{\C}(B )\simeq \Proj^s(B)$ is the commutative monoid of stable projections over $B$.  For unital  $B$  
 we have the classical group-completion picture
\begin{equation}\label{z345t345z354z35} \Proj^s(B)^{\group}\simeq\Omega^\infty K^{\topp}(B)
\end{equation} of $C^{*}$-algebra $K$-theory,
\cite[Corollary~10.8]{orbita}. 
{For classical spectrum models of $C^{*}$-algebra $K$-theory,
see \cite{bunke-joachim-stolz-k}; for structured and equivariant models,
see \cite{joachim-ko,dellambrogio-emerson-kandelaki-meyer}.}
For the purpose of the present paper, we simply define the $C^{*}$-algebra $K$-theory spectrum   in terms of mapping spectra of the stable $E$-theory category
{\cite[Definition~9.3]{orbita}} by 
 $K^{\topp}(B):= E(\C,B)$.
  It is therefore natural, also for more general $D$,  to study the natural  comparison map
 \begin{equation}\label{gwerrwegw}\eta_{D,B}:{\cM}_{D}(B )^{\group} \longrightarrow \Omega^{\infty}E(D,B)\ .
\end{equation}
  constructed in \eqref{eq:intro-first-comparison}.
 Recall that the stable $E$-theory $\infty$-category is defined in terms of universal properties and constructed as a Dwyer--Kan localization forcing them
 {\cite[Definition~7.9]{orbita}}. This description  makes understanding its mapping anima a difficult task. 
 One motivation of the   present paper is   to 
  find   descriptions of these mapping anima in terms of  the group-completed mapping spaces 
  ${\cM}_{D}(B )$.

 In general, the naive comparison \eqref{gwerrwegw} is not an equivalence. Recall that already for \eqref{z345t345z354z35} we must assume that $B$ is unital. In order to extend this equivalence to non-unital $B$, we must apply the following reduced construction: {For a functor $F:\nCalg_{\sepa}\to\cC$ taking values in a pointed
$\infty$-category admitting fibres of split epimorphisms, its {reduced construction $F^{\fred}$, defined in Section~\ref{sec:reduced-construction}, is computed by}}
\begin{equation}\label{hetrgertget}{F^{\fred}(B)\simeq \Fib_0\bigl(F(B^u)\to F(\C)\bigr),}
\end{equation}
 {where  the arrow is induced by the structure map $B^{u}\to \C$ of  the unitalization of $B$.}
 The equivalence
\begin{equation}\label{jrzjzthtzhrthtz} (\Proj^s(-)^{\group})^{\fred}\simeq \Omega^\infty K^{\topp}(-)\simeq  \Omega^{\infty}E(\C,-):\ \nCalg_{\sepa}\to  \CGrp(\Ani)
\end{equation}
is then the  usual projection-space picture of $C^{*}$-algebra $K$-theory
\cite[Theorem~10.7 and Proposition~10.9]{orbita}.

Although ${\cM}_D^{\group}$ is homotopy invariant and inverts corner
inclusions, the reduced construction need not preserve the latter
property. We therefore introduce the stabilized reduced construction \begin{equation}\label{gewrfrefwrf}F^{\sfred}(-):=F^{\fred}(-\otimes\bbK)\ .
\end{equation} Indeed, when $F$ is homotopy invariant, the stabilized reduced construction $F^{\sfred}$ is again homotopy invariant and in addition inverts corner
inclusions. We are led to introduce the  reduced representation group functor
\[
{\cT_{D}:=({\cM}_{D}^{\group})^{\sfred}:\nCalg_{\sepa}\to \CGrp(\Ani)}
\]
and to consider the comparison
 \begin{equation}\label{gwerrwegw1}\tau_D\colon \cT_{D}(- ) \to \Omega^{\infty}E(D,-):\nCalg_{\sepa}\to \CGrp(\Ani)
\end{equation}
constructed  in \eqref{eq:intro-tau}. We now observe that the target of  the  comparison  \eqref{gwerrwegw1}  sends pullbacks of surjective maps to cartesian squares. 
So given $D$, if \eqref{gwerrwegw1} is an equivalence, also the reduced representation group functor $\cT_{D} $ must have this exactness property. We will show by example that this is not true in general.

We say that $D$ is $\tau$-admissible if $ \eqref{gwerrwegw1} $ is an equivalence.
 In Proposition \ref{prop:intro-semiprojective} we show that nuclear, quasidiagonal and semiprojective 
 $D$ in $\nCalg_{\sepa}$ are $\tau$-admissible. 
  However,  $\tau$-admissibility  is restrictive. 
  By Examples~\ref{ex:intro-nonexact} and 
 \ref{ex:uhf-unresolved-nonexact} the nuclear semiprojective Cuntz algebra $\cO_{\infty}$ and the  nuclear quasidiagonal UHF algebra $M_{2^\infty}$  are not $\tau$-admissible. 
 This illustrates the need for the resolved construction.

To obtain exactness under weaker hypotheses on $D$, we pass to
the resolved functor \[
{{\cR}_D:\nCalg_{\sepa}\to \CGrp(\Ani)\ .}
\]
For a shape system $(D_n)_{n\in\nat}$  for
$D$ (see Definition \ref{def:intro-shape-system}), we compute it by the formula
\begin{equation}\label{gwerferfgwergewrg}
 {\cR}_D\simeq\lim_n \cT_{D_{n}}.
\end{equation}We have a  canonical transformation
\begin{equation}\label{eq:intro-rho}
 \rho_D:\cT_D\longrightarrow {\cR}_D: \nCalg_{\sepa}\to \CGrp(\Ani)
\end{equation}
constructed in \eqref{eq:T-to-R} and call $D$ $\rho$-admissible if $\rho_D$ is an equivalence.  We have a comparison
\begin{equation}\label{gwerrwegw3}\alpha_{D}:{\cR}_{D} \to \Omega^{\infty}E(D,-):\nCalg_{\sepa}\to \CGrp(\Ani)\ .
\end{equation} constructed in \eqref{eq:intro-alpha}.
{Theorem~\ref{thm:overview-one-splitting} below formulates conditions on $D$
under which this map admits a natural section or is an equivalence. {In the latter case we call $D$ $\alpha$-admissible; see Definition~\ref{def:alpha-admissible}.}
{By Example~\ref{UHF}, the UHF algebra $M_{2^{\infty}}$ is $\alpha$-admissible, while $\cO_{\infty}$ is not by {Proposition~\ref{ex:U-nq-oinfinity}}.}

In the case when Theorem~\ref{thm:overview-one-splitting} predicts a splitting,
 Theorem~\ref{thm:overview-one-hausdorff} gives   conditions  on $D$ implying that the induced map on Hausdorff quotients of component groups  is an isomorphism.
}

{The constructions introduced above may appear ad hoc. We now
describe them using categorical language. This approach is
useful for constructing the comparison maps to the $E$-theory mapping
spaces and describing the classes of algebras for which they are
equivalences in terms of homotopy theoretic conditions.}

{We consider the Dwyer--Kan localization
$L_{h,K}:\nCalg_{\sepa}\to L_K\nCalg_{\sepa,h}$ which inverts homotopy
equivalences and corner inclusions. Its target is semiadditive, so
the mapping anima carry natural commutative monoid structures
as explained in Section~\ref{sec:universal-functors}.}

\begin{ddd}\label{def:mapping-monoid}
{For $D$ in $\nCalg_{\sepa}$, we define the representation monoid functor by}
\[
 {{\cM}_D(-):=\Hom_{L_K\nCalg_{\sepa,h}}(D,-)
 :\nCalg_{\sepa}\to\CMon(\Ani).}
\]
\end{ddd}
Remark~\ref{kheoprhtertg} identifies this functor with the
representation monoid functor  obtained from the topological space of compact representations above.
We now consider the universal 
additive completion
$L^{\add}:L_K\nCalg_{\sepa,h}\longrightarrow (L_K\nCalg_{\sepa,h})^{\add}$.
By \eqref{eq:additive-completion-mapping}, its effect on mapping monoids is group completion. 
In particular
we obtain  the following categorical interpretation of the    representation group functor:
\begin{equation}\label{gerfwerfwerfwe}
 {\cM}_D(-)^{\group}\simeq \Hom_{(L_K\nCalg_{\sepa,h})^{\add}}(D,-)\colon\nCalg_{\sepa}\longrightarrow\CGrp(\Ani).
 \end{equation}  
 
 {A functor from $\nCalg_{\sepa}$ to an additive $\infty$-category
is homological if it is reduced,  is homotopy invariant,
inverts corner inclusions, sends short exact sequences to fibre
sequences, and preserves Schochet fibrant cartesian squares
(Definition~\ref{def:homological}). The $E$-theory functor
\[
{e_{\sepa}:\nCalg_{\sepa}\to \EE_{\sepa}}
\]
is the universal homological functor; see 
{Proposition~\ref{prop:universal-homological-existence}.}}
{By the universal properties of $L_{h,K}$ and $L^{\add}$,  the $E$-theory functor 
 factors through $L^{\add}\circ L_{h,K}$
via $e_{\mathrm{add}}:(L_K\nCalg_{\sepa,h})^{\add}\to \EE_{\sepa} $ as in \eqref{eq:additive-E-theory-factorization}.
Under the   identification  \eqref{gerfwerfwerfwe}, the map
induced by $e_{\mathrm{add}}$ on mapping groups is the comparison
\eqref{gwerrwegw}.}

 {Categorically, the reduced construction in \eqref{hetrgertget} is restriction to and  right Kan extension
from unital algebras followed by reduction at zero, as described
in Section~\ref{sec:reduced-construction}. Further recall the stabilized reduced construction \eqref{gewrfrefwrf}.}

\begin{ddd}\label{def:intro-T}
{For $D$ in $\nCalg_{\sepa}$, we define the  reduced representation
group functor by}
\[
 {\cT_D:=({\cM}_D^{\group})^{\sfred}
 :\nCalg_{\sepa}\to\CGrp(\Ani).}
\]
\end{ddd}


%
%
%

 {To define the resolved functor ${\cR}_D$ categorically, we introduce the concept of {additively continuous functors}. A functor to an additive
$\infty$-category is {additively continuous} if it is homotopy
invariant, inverts corner inclusions, and preserves finite
products and countable filtered colimits
(Definition~\ref{def:admissible}).}
Proposition~\ref{prop:universal-unreduced-category} establishes
the existence of the universal {additively continuous functor}, which we denote by
\begin{equation}\label{t435ztreh}
j:\nCalg_{\sepa}\longrightarrow\mathcal U.
\end{equation}
\begin{ddd}\label{kohprtgerthertgertgertge}
{We define the resolved   stable reduced representation group  functor by
\[
 {\cR}_D(-):=\Hom_{\cU}(j(D),j^{\sfred}(-))
 :\nCalg_{\sepa}\to\CGrp(\Ani).
\]}
\end{ddd}
The shape-system formula \eqref{eq:R-shape-formula} in
Lemma~\ref{lem:universal-unreduced-mapping} identifies this functor
with the limit in \eqref{gwerferfgwergewrg}.   

Since the $E$-theory functor is {additively continuous}, the universal
property \eqref{eq:universal-unreduced-property} gives the factorization
$e_{\mathcal U}:\mathcal U\to\EE_{\sepa}$ described in
Section~\ref{sec:universal-functors}. This functor induces the comparison transformation
\[
 \alpha_D:{\cR}_D\longrightarrow\Omega^\infty E(D,-),
\]%
constructed in \eqref{eq:intro-alpha}.


\begin{ddd}\label{def:U-nq}
{We define full subcategories
\[
{\mathcal U_{\mathrm{nqs}}\subseteq\mathcal U_{\mathrm{nq}}\subseteq\mathcal U}
\]
closed 
under equivalences, retracts, finite products, and countable
filtered colimits: The subcategory $\mathcal U_{\mathrm{nq}}$ is generated
by the images under $j$ of separable nuclear quasidiagonal algebras.
For $\mathcal U_{\mathrm{nqs}}$ we in addition require the generators to be semiprojective.}
\end{ddd}

For $D$ in  $\nCalg_{\sepa}$  we will write, by abuse of notation,  $D\in \mathcal U_{\mathrm{?}}$ for $j(D)\in \mathcal U_{\mathrm{?}}$ for $?\in \{\mathrm{nqs},\mathrm{nq}\}$.

 {We consider $D$ in $\nCalg_{\sepa}$.}

\begin{overviewtheorem}\label{thm:overview-one-splitting}
\begin{enumerate}
\item {If $D\in\mathcal U_{\mathrm{nq}}$, then ${\cR}_D$ is homological.}
\item {If ${\cR}_D$ is homological, then $\alpha_D$ admits a natural section. Consequently,}
\[
 {{\cR}_D\simeq\Omega^\infty E(D,-)\oplus{\widehat{\cR}}_D.}
\]
\item {If $D\in\mathcal U_{\mathrm{nqs}}$, or if ${\cR}_D$ is homological and $D$ is $\rho$-admissible, then $D$ is $\alpha$-admissible and ${\widehat{\cR}}_D\simeq0$.}
\end{enumerate}
\end{overviewtheorem}
The functor ${\widehat{\cR}}_D:=\Fib_0(\alpha_D)$  in statement 2. is called the
remainder functor.

{Homologicality is established in
Theorem~\ref{thm:intro-exactness}. The section is obtained
from the universal class by the Yoneda lemma in
Proposition~\ref{prop:intro-splitting}. {Vanishing of the remainder follows from Corollary~\ref{cor:U-nqs-remainder} for $D\in\mathcal U_{\mathrm{nqs}}$, and from Proposition~\ref{prop:rho-homological-alpha} when ${\cR}_D$ is homological and $D$ is $\rho$-admissible.}}

  The class $\mathcal U_{\mathrm{nqs}}$ contains finite-dimensional algebras {by Lemma~\ref{lem:finite-dimensional-U-nqs}}, $C_0(\R)$, and 
$C(S^1)$. It also contains $C^*_{\max}(F_r)$ for finite or countably
infinite rank $r$ (Example~\ref{ex:U-nq-free-products}) and the UHF
algebra $M_{2^\infty}$, since it is a countable filtered colimit of
finite-dimensional algebras. 
   Further examples are discussed in 
Examples~\ref{ex:U-nq-profinite-crossed-product} and  \ref{ex:U-nq-irrational-rotation}. By Proposition~\ref{prop:U-nq-max-tensor}, the class
$\mathcal U_{\mathrm{nq}}$ is closed under taking maximal tensor products. 
{On the other hand, $ \mathcal U_{\mathrm{nq}}$ is a proper
subcategory: $\mathcal O_\infty\not\in  \mathcal U_{\mathrm{nq}}$
since ${\cR}_{\mathcal O_\infty}$ is not homological;  {see {Proposition~\ref{ex:U-nq-oinfinity}}.}}

We consider $D$ and $B$ in $\nCalg_{\sepa}$.
Both component groups $\pi_0{\cR}_D(B)$ and $E_0(D,B)$ carry a shape
topology and a condensed topology, defined in
Definitions~\ref{def:intro-shape-topology} and \ref{def:condensed-topology},
respectively. The two topologies coincide on $\pi_0{\cR}_D(B)$ by
Proposition~\ref{prop:R-condensed-shape} and on $E_0(D,B)$ by
\cite[Propositions~4.25 and 4.42]{budu}.
{The notation $(-)_{\mathrm{Hd}}$ for the Hausdorff quotient is therefore unambiguous.}

We consider  $D$ in $\nCalg_{\sepa}$.
\begin{overviewtheorem}\label{thm:overview-one-hausdorff}
{{If   $D\in\mathcal U_{\mathrm{nq}}$,}
then the comparison $\alpha_{D,B}$
induces an isomorphism of topological groups
\[
 (\pi_0{\cR}_D(B))_{\mathrm{Hd}}
 \xrightarrow{\cong}E_0(D,B)_{\mathrm{Hd}}.
\]
Under the splitting of Theorem~\ref{thm:overview-one-splitting},
the group $\pi_0{\widehat{\cR}}_D(B)$ is contained in the closure
of zero in $\pi_0{\cR}_D(B)$.}
\end{overviewtheorem}
{We prove this theorem in Section~\ref{sec:hausdorff-source-continuity}.}

\begin{prob}\label{prob:intro-hausdorff}
Under the assumptions of Theorem~\ref{thm:overview-one-hausdorff},
is
\[
 \alpha_{D,B}:{\cR}_D(B)\longrightarrow \Omega^\infty E(D,B)
\]
an equivalence of commutative groups in anima?\hB
\end{prob}

\begin{prob}
{We do not know whether
\[
 \mathcal U_{\mathrm{nqs}}=\mathcal U_{\mathrm{nq}}.
\]
An affirmative answer would imply vanishing of the remainder
throughout $\mathcal U_{\mathrm{nq}}$, 
see 
{Remark~\ref{rem:U-nq-nqs}} for a further discussion.

Does homologicality of ${\cR}_D$ imply $\alpha$-admissibility of $D$?
See Remark~\ref{rem:U-nq-nqs}.
}\hB\end{prob}

{For a countable discrete group $G$, Willett's group $QR(G)$ {\cite[Definition~2.11]{willett}} is defined in terms of  finite-dimensional approximate
representations, stable approximate unitary equivalence,
group completion, and by forming an inverse limit over the controls; {see Definition~\ref{def:willett-QR} for details.}}

\begin{overviewtheorem}\label{thm:overview-one-quasirepresentations}\label{thm:overview-qr-comparison}
{For every countable discrete group $G$, there is a
canonical isomorphism of topological groups
\[
 QR(G)\cong
 \bigl(\pi_0{\cR}_{C^*_{\max}(G)}(\C)\bigr)_{\mathrm{Hd}}.
\]}
\end{overviewtheorem}

 The proof is given in Section~\ref{sec:willett-equivalence}.

{
\begin{kor}\label{cor:overview-qr-E}
{{Let $G$ be a countable discrete group such that $C^*_{\max}(G)\in\mathcal U_{\mathrm{nq}}$.}
Then there is a canonical isomorphism of topological groups
\[
 QR(G)\cong E_0(C^*_{\max}(G),\C)_{\mathrm{Hd}}.
\]}
\end{kor}
}

\begin{proof}
{Apply Theorem~\ref{thm:overview-one-hausdorff} with $D=C^*_{\max}(G)$ and $B=\C$, and combine the resulting Hausdorff comparison with Theorem~\ref{thm:overview-qr-comparison}.}
\end{proof}

This corollary is an $E$-theoretic version of  Willett's
\cite[Theorem~10.3]{willett} stating $$QR(G) \cong KL(C^*_{\max}(G),\C)$$ under the assumption that $C^*_{\max}(G)$ is quasidiagonal and has the local lifting property (LLP). On the one hand, the assumption of nuclearity for the generators of $\cU_{\mathrm{nq}}$ is stronger than the LLP so that there might be groups covered by Willett's theorem, but not by 
 Corollary \ref{cor:overview-qr-E}. On the other hand, the class $\cU_{\mathrm{nq}}$ contains many maximal group $C^{*}$-algebras which are not themselves nuclear and quasi-diagonal.


 {For example, let $G_1,\ldots,G_r$ be countable groups whose full group $C^*$-algebras are nuclear and quasidiagonal, and put $\Gamma:=G_1*\cdots*G_r$. The identity $C^*_{\max}(\Gamma)\cong C^*_{\max}(G_1)*_{\C}\cdots*_{\C}C^*_{\max}(G_r)$ and {Lemma~\ref{lem:U-free-product-closure}} give $C^*_{\max}(\Gamma)\in\mathcal U_{\mathrm{nq}}$. Concrete examples are free products of finite and finitely generated free abelian groups, including $\Z/2*\Z/3$ and the free groups $F_r=\Z*\cdots*\Z$. Willett's theorem also applies to these free products: their full group algebras are quasidiagonal and have LLP, since both properties are preserved by unital full free products \cite[Corollaries~12 and~14]{boca-free-products}.}

 {Corollary~\ref{cor:overview-qr-E} also applies to finite direct products of such free products. Indeed, their full group algebras are maximal tensor products of algebras in $\mathcal U_{\mathrm{nq}}$, so Proposition~\ref{prop:U-nq-max-tensor} applies. In particular, $C^*_{\max}(F_r\times F_s)\in\mathcal U_{\mathrm{nq}}$ for $r,s\geq2$; see Example~\ref{ex:U-nq-tensor-products}. Roydor has proved that these algebras do not have LLP \cite[Theorem~A]{roydor-non-llp}. Thus our comparison with Hausdorffized $E$-theory applies to these groups, including $F_2\times F_2$, whereas Willett's Theorem~10.3 does not apply.}

We conclude the paper by giving pointers to related results in the $C^{*}$-algebra literature.
 
\section{\texorpdfstring{{Representation group functors and $E$-theory}}{Representation group functors and E-theory}}
\label{sec:introduction}


\subsection{Universal functors on $C^*$-algebras}
\label{sec:universal-functors}
We let $\nCalg_{\sepa}$  denote   the category
 of separable $C^*$-algebras and possibly nonunital homomorphisms. 
Let $\cA$ be an $\infty$-category and $F:\nCalg_{\sepa}\to\cA$ be some  functor.
We use the following terminology.
\begin{ddd}\label{def:functor-properties}\mbox{}
 \begin{enumerate}
\item $F$ is \emph{reduced} if $F(0)$ is a zero object
of $\cA$.
\item $F$ is \emph{homotopy invariant} if it sends
homotopy equivalences to equivalences.
\item $F$ \emph{inverts corner inclusions} if it sends  all
upper-left corner inclusions $B\to B\otimes\bbK$ to equivalences.
\item \label{def:functor-exactness} $F$ is \emph{exact} if, for every short exact
sequence $0\to I\to B\to Q\to0$, the square
\[
 \xymatrix{
 F(I)\ar[r]\ar[d]&F(B)\ar[d]\\
 F(0)\ar[r]&F(Q)
 }
\]
is cartesian. 
\item \label{hjjtjrtzh} $F$ is \emph{Schochet exact}
if it sends every cartesian square in $\nCalg_{\sepa}$ whose
right vertical map is a Schochet fibration to a cartesian
square in $\cA$.
\item $F$ \emph{preserves countable filtered colimits}
if it sends every countable filtered colimit cocone in
$\nCalg_{\sepa}$ to a colimit cocone in $\cA$.
\end{enumerate}
\end{ddd}

\begin{rem}
{For details on  homotopy invariance, inversion of corner inclusions (called
stability there), and reducedness in an $\infty$-categorical context, see
\cite[Definition~2.4(1)--(3)]{orbita}. For a reduced functor,
the cartesian-square condition for exactness is the fibre-sequence
condition of \cite[Definition~2.4(4)]{orbita}.
Schochet fibrations and Schochet exactness are discussed in
\cite[Section~3, in particular Definition~3.14]{orbita}.
\hB}
\end{rem}

\begin{ddd}\label{def:admissible}
A functor $F:\nCalg_{\sepa}\to\cA$ to an additive
$\infty$-category is {\emph{additively continuous}} if it is homotopy invariant, inverts corner inclusions, and preserves finite products
and countable filtered colimits.
\end{ddd}

\begin{ddd}\label{def:homological}
A functor $F:\nCalg_{\sepa}\to\cA$ to an additive
$\infty$-category is \emph{homological} if it is reduced,
homotopy invariant, exact,  Schochet exact, and inverts corner inclusions.
\end{ddd}

For stable targets, Schochet exactness follows from homotopy
invariance and exactness \cite[Lemma~5.2 and Remark~5.11]{orbita}.

We first construct the universal homotopy-invariant functor that inverts corner inclusions as a composition
\begin{equation}\label{eq:homotopy-stable-localization}
 L_{h,K}:=L_K\circ L_h:
 \nCalg_{\sepa}\xrightarrow{L_h}\nCalg_{\sepa,h}
 \xrightarrow{L_K}L_K\nCalg_{\sepa,h}.
\end{equation}
The Dwyer--Kan localization $L_h$ inverts homotopy
equivalences. Its target is presented by the topologically
enriched category of separable $C^*$-algebras with
point-norm homomorphism spaces
\cite[Proposition~3.5 and Corollary~3.8]{orbita}.
The second localization $L_K$ inverts  the upper-left corner inclusions $B\to B\otimes\bbK$. It is the left Bousfield
localization represented by $B\mapsto B\otimes\bbK$.
The localized categories  are left exact, and
$L_K\nCalg_{\sepa,h}$ is in addition semiadditive
\cite[Corollaries~4.3 and 4.11, Proposition~4.4]{orbita}.

For every $\infty$-category $\cA$, precomposition
induces an equivalence
\begin{equation}\label{eq:homotopy-stable-universal-property}
 L_{h,K}^*:
 \Fun(L_K\nCalg_{\sepa,h},\cA)
 \xrightarrow{\simeq}
 \Fun^{h,K}(\nCalg_{\sepa},\cA).
\end{equation}
The right-hand side is the full subcategory of
homotopy-invariant functors {that invert corner inclusions}.  

\begin{ddd}\label{def:universal-homological-functor}
We define
$e_{\sepa}:\nCalg_{\sepa}\to\EE_{\sepa}$
as the universal homological functor.
\end{ddd}
By definition, for every additive, left-exact
$\infty$-category $\cA$ we have an 
equivalence
\begin{equation}\label{eq:homological-universal-property}
 e_{\sepa}^*:
 \Fun^{\lex}(\EE_{\sepa},\cA)
 \xrightarrow{\simeq}
 \Fun^{\mathrm{hom}}(\nCalg_{\sepa},\cA).
\end{equation}
Here the superscript  $\mathrm{hom}$ indicates the full subcategory of  homological functors.

\begin{prop}\label{prop:universal-homological-existence}
The universal homological functor exists, and
its target $\EE_{\sepa}$ is stable.
\end{prop}
\begin{proof}
Starting from $L_{h,K}$, we form a  further left exact  localization
\[
 \bar e_{\sepa}:L_K\nCalg_{\sepa,h}
       \longrightarrow\EE_{\sepa},
 \qquad e_{\sepa}=\bar e_{\sepa}\circ L_{h,K},
\]
which imposes exactness and Bott periodicity. Its target is stable
\cite[Propositions~5.9, 7.3 and 7.5, Definition~7.9]{orbita}.
The resulting functor is homological and satisfies
\eqref{eq:homological-universal-property} by
\cite[equation~(7.9)]{orbita}.
\end{proof}

The construction in the proof identifies
$e_{\sepa}$ with the usual $\infty$-categorical
$E$-theory functor of \cite[Definition~7.9]{orbita}.
Thus $\EE_{\sepa}$ is the separable $E$-theory category,
and its mapping spectra are the $E$-theory spectra.
{By \cite[Remark~14.8]{orbita}, the composite}
\[
 {\nCalg_{\sepa}\xrightarrow{e_{\sepa}}\EE_{\sepa}
 \xrightarrow{\ho}\ho(\EE_{\sepa})}
\]
{is equivalent to the
classical separable $E$-theory functor of  \cite{higson-excision} and 
\cite{connes-higson-asymptotic}.}
{For its classical universal property, see \cite[Theorem~3.6]{higson-excision} and the textbook treatment in \cite[Theorem~25.6.1]{blackadar-kbook}.}

The universal additive completion  $\mathcal B\to\mathcal B^{\mathrm{add}}$  of a  semiadditive $\infty$-category $\mathcal B$  is obtained by
change of enrichment along symmetric monoidal group
completion $(-)^{\group}$ from \eqref{gerfwregwerfwewerf};
\cite[Theorem~5.1]{gepner-groth-nikolaus};
\cite[Corollaries~3.6.19 and 5.7.6]{gepner-haugseng-enriched}.
 Applying the universal additive completion  to
$L_K\nCalg_{\sepa,h}$ gives the second functor in
\begin{equation}\label{eq:additive-completion-category}
 L_{h,K}^{\mathrm{add}}:
 \nCalg_{\sepa}\xrightarrow{L_{h,K}}L_K\nCalg_{\sepa,h}
 \xrightarrow{L^{\add}}(L_K\nCalg_{\sepa,h})^{\mathrm{add}},
\end{equation}
For  $D$ and $B$  in $L_K\nCalg_{\sepa,h}$ we have an equivalence of mapping groups in anima
\begin{equation}\label{eq:additive-completion-mapping}
 \Hom_{(L_K\nCalg_{\sepa,h})^{\mathrm{add}}}(D,B)
 \simeq\Hom_{L_K\nCalg_{\sepa,h}}(D,B)^{\group}.
\end{equation}
For every additive $\infty$-category $\cA$, precomposition
induces an equivalence
\begin{equation}\label{eq:additive-completion-universal-property}
 (L_{h,K}^{\mathrm{add}})^*:
 \Fun^\oplus((L_K\nCalg_{\sepa,h})^{\mathrm{add}},\cA)
 \xrightarrow{\simeq}
 \Fun^{h,K,\times}(\nCalg_{\sepa},\cA).
\end{equation}
The superscripts $\oplus$ and $\times$ indicate preservation of finite direct sums and finite products, respectively, including the zero and terminal objects. 

The functor $e_{\sepa}$ is homotopy invariant,
inverts corner inclusions, and preserves finite products. Its target
$\EE_{\sepa}$ is stable, hence additive, by
Proposition~\ref{prop:universal-homological-existence}.
Thus \eqref{eq:additive-completion-universal-property} gives the dotted arrow in the   
   commutative diagram
\begin{equation}\label{eq:additive-E-theory-factorization}
 \xymatrix@C=5em@R=3em{
 \nCalg_{\sepa}\ar[r]^{L_{h,K}^{\mathrm{add}}}
       \ar[dr]_{e_{\sepa}}
   &(L_K\nCalg_{\sepa,h})^{\mathrm{add}}
       \ar@{.>}[d]^{e_{\mathrm{add}}}\\
   &\EE_{\sepa}.
 }
\end{equation}

\begin{ddd}\label{def:universal-admissible}
We define $j:\nCalg_{\sepa}\to\mathcal U$
as the universal {additively continuous functor}.
\end{ddd}
By definition, 
 for every additive $\infty$-category $\cA$
admitting countable filtered colimits, precomposition
induces an equivalence
\begin{equation}\label{eq:universal-unreduced-property}
 j^*:\Fun^{\oplus,\mathrm{cfil}}(\mathcal U,\cA)
 \xrightarrow{\simeq}
 \Fun^{{\ac}}(\nCalg_{\sepa},\cA).
\end{equation}
The superscript $\ac$ indicates additively continuous functors (Definition~\ref{def:admissible}), and $\mathrm{cfil}$ indicates preservation of countable filtered colimits. 

\begin{prop}\label{prop:universal-unreduced-category}
The universal {additively continuous functor} exists, and
its target  admits countable filtered colimits.
\end{prop}

\begin{proof}
{We apply Proposition~\ref{prop:adjoin-prescribed-colimits} to
$(L_K\nCalg_{\sepa,h})^{\add}$ and the images under $L_{h,K}^{\add}$
of countable filtered colimit cocones in $\nCalg_{\sepa}$.}
We  therefore consider the additive presheaf category
\[
 \mathcal P:=
 \Fun^\oplus\!\left(
 ((L_K\nCalg_{\sepa,h})^{\mathrm{add}})^{\op},
 \CGrp(\Ani)\right)\]
For $B$ in $(L_K\nCalg_{\sepa,h})^{\add}$, the Yoneda embedding is given by
 \begin{equation}\label{eq:additive-yoneda-presheaf}
 h_B:=\Hom_{(L_K\nCalg_{\sepa,h})^{\mathrm{add}}}(-,B).
\end{equation}
We let 
\begin{equation}\label{gerferfwerfwer}
\mathsf L:\mathcal P\leftrightarrows \mathcal P_{\mathrm{cont}}:\incl
\end{equation}
be the left Bousfield localization at the collection of maps
\begin{equation}\label{eq:unreduced-continuity-relations}
 \colim_i h_{L_{h,K}^{\add}(B_i)}\longrightarrow h_{L_{h,K}^{\add}(\colim_i B_i)}
\end{equation} for all countable filtered diagrams
  $(B_i)$  
in $\nCalg_{\sepa}$. 
{Write $\upsilon$ for the unit of the adjunction
\eqref{gerferfwerfwer}}
\begin{equation}\label{eq:continuous-presheaf-unit}
 {\upsilon:\id_{\mathcal P}\longrightarrow\incl\circ\mathsf L.}
\end{equation} 
 Then  $\mathcal U$  is defined as the full
subcategory of $\mathcal P_{\mathrm{cont}}$ generated
by the objects $\mathsf Lh_B$ for $B$ in $(L_K\nCalg_{\sepa,h})^{\mathrm{add}}$  under countable filtered
colimits. We further define the functor
\begin{equation}\label{gerwfgwergwerg}j:=\mathsf L\circ h_{-}\circ L^{\add}_{h,K}:\nCalg_{\sepa}\to \cU.
\end{equation}
{Proposition~\ref{prop:adjoin-prescribed-colimits} gives additivity
and countable filtered colimits in $\mathcal U$. It also gives
preservation of these colimits by $j$, since the prescribed cocones
are the images of countable filtered colimit cocones of $C^*$-algebras.
We obtain the chain of equivalences}
\begin{equation}\label{eq:universal-additive-continuity-proof-chain}
 {
 \Fun^{\oplus,\mathrm{cfil}}(\mathcal U,\cA)
 \xrightarrow[\simeq]{(\mathsf L\circ h_-)^*}
 \Fun^{\oplus,\mathrm{rel}}\bigl((L_K\nCalg_{\sepa,h})^{\add},\cA\bigr)
 \xrightarrow[\simeq]{(L_{h,K}^{\add})^*}
 \Fun^{{\ac}}(\nCalg_{\sepa},\cA).
 }
\end{equation}
Here $\mathrm{rel}$ indicates additive functors sending the prescribed cocones to colimit cocones.
{The first equivalence is
\eqref{eq:adjoin-prescribed-colimits-universal} for the relations
\eqref{eq:unreduced-continuity-relations}. The second is the restriction
of \eqref{eq:additive-completion-universal-property} to functors
respecting those relations. Their composite is $j^*$ from
\eqref{eq:universal-unreduced-property}.}
\end{proof}

Since $j$ is homotopy invariant, inverts corner inclusions,
and preserves finite products, and $\mathcal U$ is
additive, \eqref{eq:additive-completion-universal-property}
gives a unique additive functor
\begin{equation}\label{gergewrfewrgwerge} j_{\mathrm{add}}:
 (L_K\nCalg_{\sepa,h})^{\mathrm{add}}\longrightarrow\mathcal U,
 \qquad j\simeq j_{\mathrm{add}}\circ L_{h,K}^{\mathrm{add}}.
\end{equation}

The category $\EE_{\sepa}$ admits countable filtered colimits, and $e_{\sepa}$ preserves them, by \cite[Corollary~3.45(2) and Propositions~3.67--3.68]{budu}.
It is homotopy invariant, inverts corner inclusions, and preserves
finite products by left exactness. Since $\EE_{\sepa}$
is stable, hence additive, $e_{\sepa}$ is {additively continuous}.
The universal property
\eqref{eq:universal-unreduced-property} therefore gives
{a unique functor preserving finite direct sums and countable
filtered colimits}
\begin{equation}\label{eq:universal-E-factorization}
 {e_{\mathcal U}:\mathcal U\longrightarrow\EE_{\sepa},
 \qquad e_{\mathcal U}\circ j\simeq e_{\sepa}.}
\end{equation}
Uniqueness in
\eqref{eq:additive-completion-universal-property}
identifies $e_{\mathcal U}\circ j_{\mathrm{add}}$ with
$e_{\mathrm{add}}$, since both composites with
$L_{h,K}^{\mathrm{add}}$ agree with $e_{\sepa}$.
{We summarize these constructions in the following commutative diagram}
\[
 \xymatrix@C=3em@R=3.5em{
 \nCalg_{\sepa}
       \ar[r]^-{L_{h,K}^{\mathrm{add}}}
       \ar@/^2pc/[rr]^j
       \ar[drr]_{e_{\sepa}}
   &(L_K\nCalg_{\sepa,h})^{\mathrm{add}}
       \ar[r]^-{j_{\mathrm{add}}}
       \ar[dr]^{e_{\mathrm{add}}}
   &\mathcal U\ar@{.>}[d]^{e_{\mathcal U}}\\
   &&\EE_{\sepa}.
 }
\]
\subsubsection*{{Appendix: Adjoining colimits subject to prescribed relations}}
\label{sec:universal-functors-appendix}

Let $\mathcal B$ be a small additive $\infty$-category and fix a
small collection  of cocones $(b_i)\to b$ with countable filtered
index categories. Set
$\mathcal P_{\mathcal B}:=\Fun^\oplus(\mathcal B^{\op},\CGrp(\Ani))$
and $h_b:=\Hom_{\mathcal B}(-,b)$. 
Since ${\mathcal P_{\mathcal B}}$ is presentable, there exists a left Bousfield localization \begin{equation}\label{hertgertge}
\mathsf L:\mathcal P_{\cB} \leftrightarrows \mathcal P_{\mathrm{loc}}:\incl
\end{equation}
{We write}
\begin{equation}\label{eq:prescribed-colimits-localization-unit}
{\upsilon:\id_{\mathcal P_{\mathcal B}}\longrightarrow\incl\circ\mathsf L}
\end{equation}
{for its unit. The localization is taken}

at the maps $\colim_i h_{b_i}\to h_b$ associated to the prescribed
cocones. The full subcategory ${\mathcal P_{\mathrm{loc}}}$ consists of those 
 presheafs $P$ that satisfy 
 $P(b)\simeq \lim_i P(b_i)$ for these cocones. Let $\mathcal U_{\mathcal B}$ be the full subcategory
of $\mathcal P_{\mathrm{loc}}$ generated under countable filtered
colimits by the objects $\mathsf L h_b$ for all $b$ in $\cB$, and define
$k\colon\mathcal B\to\mathcal U_{\mathcal B}$ by $k(b):=\mathsf L h_b$.

\begin{prop}\label{prop:adjoin-prescribed-colimits}
The category $\mathcal U_{\mathcal B}$ is additive and admits
countable filtered colimits. The functor $k$ is additive and sends the
prescribed cocones to colimit cocones. For every additive
$\infty$-category $\mathcal D$ admitting countable filtered colimits,
precomposition induces an equivalence
\begin{equation}\label{eq:adjoin-prescribed-colimits-universal}
 k^*\colon
 \Fun^{\oplus,\mathrm{cfil}}(\mathcal U_{\mathcal B},\mathcal D)
 \xrightarrow{\simeq}\Fun^{\oplus,\mathrm{rel}}(\mathcal B,\mathcal D),
\end{equation}
where the right-hand side consists of the additive functors
sending the prescribed cocones to colimit cocones.
\end{prop}

\begin{proof}
Finite sums in $\mathcal P_{\mathcal B}$ are pointwise products
and preserve locality. Thus $\mathcal P_{\mathrm{loc}}$ is additive.
The functors $h_{-}$ and $\mathsf L$ preserve finite sums. Since binary
sums preserve filtered colimits separately in each variable,
$\mathcal U_{\mathcal B}$ is closed under finite sums. Applying $\mathsf L$
to the defining relations shows that $k$ sends the prescribed cocones
to colimit cocones. For $P$ in $\mathcal P_{\mathrm{loc}}$, we have
\begin{equation}\label{eq:localized-additive-evaluation}
 \Hom_{\mathcal P_{\mathrm{loc}}}(k(b),P)
 \simeq\Hom_{\mathcal P_{\mathcal B}}(h_b,P)\simeq P(b).
\end{equation}

Let $F$ be in $\Fun^{\oplus,\mathrm{rel}}(\mathcal B,\mathcal D)$.
We define $\Phi_F\colon\mathcal D\to\mathcal P_{\mathrm{loc}}$ by
$\Phi_F(X):=\Hom_{\mathcal D}(F(-),X)$. The prescribed colimit
conditions on $F$ imply locality of $\Phi_F(X)$.
By \eqref{eq:localized-additive-evaluation},
$\Hom_{\mathcal P_{\mathrm{loc}}}(k(b),\Phi_F(-))$ is corepresented
by $F(b)$. If $V\simeq\colim_i V_i$ and
$\Hom_{\mathcal P_{\mathrm{loc}}}(V_i,\Phi_F(-))$ is corepresented
by $\bar F(V_i)$, then
\begin{equation}\label{eq:adjoin-colimits-corepresentability}
 \begin{aligned}
 \Hom_{\mathcal P_{\mathrm{loc}}}(V,\Phi_F(X))
 &\simeq\lim_i\Hom_{\mathcal P_{\mathrm{loc}}}(V_i,\Phi_F(X))\\
 &\simeq\lim_i\Hom_{\mathcal D}(\bar F(V_i),X)
 \simeq\Hom_{\mathcal D}(\colim_i\bar F(V_i),X)
 \end{aligned}
\end{equation}
shows that ${\Hom_{\mathcal P_{\mathrm{loc}}}(V,\Phi_F(-))}$ is corepresentable.
Hence this functor is corepresentable for every $V$ in
$\mathcal U_{\mathcal B}$. 
Full faithfulness of co-Yoneda assembles
the corepresenting objects into $\bar F\colon\mathcal U_{\mathcal B}\to
\mathcal D$, functorially in $F$. Formula~\eqref{eq:adjoin-colimits-corepresentability}
shows that $\bar F$ preserves countable filtered colimits. The same
calculation for finite sums proves additivity. On generators we
obtain $\bar F\circ k\simeq F$, naturally in $F$.

Conversely, let $G$ be in
$\Fun^{\oplus,\mathrm{cfil}}(\mathcal U_{\mathcal B},\mathcal D)$.
For $V$ in $\mathcal U_{\mathcal B}$, we use
\eqref{eq:localized-additive-evaluation} and apply $G$ to obtain
\begin{equation}\label{eq:adjoin-colimits-comparison-construction}
 V\simeq\Hom_{\mathcal U_{\mathcal B}}(k(-),V)
 \xrightarrow{G}\Hom_{\mathcal D}(Gk(-),G(V))
 =\Phi_{G\circ k}(G(V)).
\end{equation}
Composing with $\Phi_{G\circ k}(f)$ for $f:G(V)\to X$
gives a map of mapping anima
\begin{equation}\label{eq:adjoin-colimits-comparison}
 \Hom_{\mathcal D}(G(V),X)\longrightarrow
 \Hom_{\mathcal P_{\mathrm{loc}}}(V,\Phi_{G\circ k}(X)).
\end{equation}
For $V=k(b)$ this map is equivalent {to} the identity. Thus \eqref{eq:adjoin-colimits-comparison} is an equivalence on
the generators. Both sides send countable filtered colimits in $V$
to limits, so it is an equivalence for every $V$ in $\mathcal U_{\mathcal B}$.
Corepresentability gives $\overline{G\circ k}\simeq G$ naturally in $G$.
This proves \eqref{eq:adjoin-prescribed-colimits-universal}.
\end{proof}

We use the notation from \eqref{hertgertge} and the setting of Proposition~\ref{prop:adjoin-prescribed-colimits}. {Let $b$ be in $\mathcal B$.}
\begin{lem}\label{lem:adjoin-colimits-representable-unit}
We have an equivalence of arrows in $\mathcal P_{\mathcal B}$
\begin{equation}\label{eq:adjoin-colimits-representable-unit}
 \bigl(h_b\xrightarrow{k}\Hom_{\mathcal U_{\mathcal B}}(k(-),k(b))\bigr)
 \simeq
 \bigl(h_b\xrightarrow{\upsilon_{h_b}} \incl\mathsf L h_b\bigr).
\end{equation}
\end{lem}
\begin{proof}
By \eqref{eq:localized-additive-evaluation}, the target of the left arrow of
\eqref{eq:adjoin-colimits-representable-unit} identifies with
$ \incl\mathsf L h_b$. Under this identification, the left arrow of
\eqref{eq:adjoin-colimits-representable-unit}
corresponds to the identity of $k(b)=\mathsf L h_b$: evaluation at
$\id_b$ sends it to $\id_{k(b)}$. Hence it is the unit at $h_b$.
\end{proof}
Let \[
{L\colon\mathcal P\leftrightarrows\mathcal P_{\mathrm{loc}}
\colon\incl}
\]
be a  left Bousfield  localization of additive $\infty$-categories
with unit $\upsilon:\id\to \incl L$. Let $p\colon X\to Y$ be a split epimorphism in
$\mathcal P$ {such that $\Fib_0(p)$ exists}.
\begin{lem}\label{lem:reflection-split-fibres}
 We have an equivalence of arrows in $\mathcal P$
\begin{equation}\label{eq:reflection-split-fibres}
 \bigl(\Fib_0(p)\xrightarrow{\Fib_0(\upsilon)}\Fib_0(\incl Lp)\bigr)
 \simeq \bigl(\Fib_0(p)\xrightarrow{\upsilon_{\Fib_{0}(p)}}\incl L\Fib_0(p)\bigr).
\end{equation}
\end{lem}
\begin{proof}
A section of $p$ gives $X\simeq \Fib_{0}(p)\oplus Y$, with $p$ the projection
onto $Y$. Since $L$ preserves finite sums, $Lp$ identifies with the
projection $L\Fib_{0}(p)\oplus LY\to LY$. Thus $\Fib_0(Lp)\simeq L\Fib_{0}(p)$.
The right adjoint $\incl$ preserves this fibre, and naturality of
$\upsilon$ along $\Fib_{0}(p)\to X$ identifies the left arrow of \eqref{eq:reflection-split-fibres}
with $\upsilon_{\Fib_{0}(p)}$.
\end{proof}

\subsection{Group completion of mapping spaces}\label{kopjrtzjrtzhzhtz}
In this section we study the
  representation monoid functor of Definition~\ref{def:mapping-monoid} and the relation of its group completion with $E$-theory mapping spaces.

\begin{rem}\label{kheoprhtertg} Let $D,B$ be in $\nCalg_{\sepa}$. In this remark we explain why 
the commutative monoid ${\cM}_D(B)$ as introduced in Definition~\ref{def:mapping-monoid} is represented by
the $E_\infty$-space of homomorphisms
$\Hom_{\nCalg_{\sepa}}(D,B\otimes\bbK)$ with the
point-norm topology. Since $L_K\simeq -\otimes \bbK$ is the left
Bousfield localization,
\begin{equation}\label{eq:representation-monoid-underlying-anima}
 {\cM}_D(B)\simeq
 \Hom_{\nCalg_{\sepa,h}}(D,B\otimes\bbK)
\end{equation}
on underlying anima. The right-hand side  is the anima  represented
by the topological homomorphism space $\Hom_{\nCalg_{\sepa}}(D,B\otimes\bbK)$
\cite[Proposition~3.5 and Corollary~3.8]{orbita}.

An explicit $E_\infty$-structure  on $\Hom_{\nCalg_{\sepa}}(D,B\otimes\bbK) $ is supplied by an action of 
the complex linear isometries operad  
\cite[Section~I.1]{may-einfty}. Write $\bbK=\mathbb K(H)$
for a separable infinite-dimensional complex Hilbert
space $H$. The $n$-th space   of this operad is the space  of linear isometric embeddings
$H^{\oplus n}\to H$, with the strong operator topology.
An
embedding $v:H^{\oplus n}\to H$ acts by
\[
 (\varphi_1,\ldots,\varphi_n)\longmapsto
 \left(a\longmapsto
 (1\otimes v)\operatorname{diag}
       (\varphi_1(a),\ldots,\varphi_n(a))(1\otimes v^*)
 \right).
\]
These operations are continuous for the point-norm
topology and compatible with operadic composition.
In particular, choosing a unitary
$H\oplus H\cong H$ in the second space of the operad  gives the block sum, which induces
addition on $\pi_0{\cM}_D(B)$. 

We finally argue that the commutative monoid structure on ${\cM}_{D}(B)$ coming from this $E_{\infty}$-structure coincides with the commutative monoid structure coming from semiadditivity.
We use the functoriality in $B$.
By \eqref{eq:representation-monoid-underlying-anima}, the underlying anima-valued functor of ${\cM}_D$ is represented by the topological space-valued functor $\Hom_{\nCalg_{\sepa}}(D,-\otimes\bbK)$.
Its $\CMon(\Ani)$-valued  refinement induced  by the $E_{\infty}$-structures on the topological  spaces $\Hom_{\nCalg_{\sepa}}(D,-\otimes\bbK) $ 
preserves finite products  and descends to
$L_K\nCalg_{\sepa,h}$, since its underlying functor is
homotopy invariant and inverts corner inclusions. On a semiadditive
$\infty$-category, every finite-product-preserving anima-valued
functor has a unique finite-product-preserving refinement to
commutative monoids \cite[Corollary~5.15]{harpaz-spans}.
Applied to $\Hom_{L_K\nCalg_{\sepa,h}}(D,-)$ and the functor described above, this identifies
the operadic structure with the canonical structure from
semiadditivity, including all higher coherences.
 \hB
\end{rem}
We now  form the objectwise   group completion
\[
{{\cM}_-(-)^{\group}:(\nCalg_{\sepa})^{\op}\times \nCalg_{\sepa}\to \CGrp(\Ani)\ .}
\]
By
\eqref{eq:additive-completion-mapping} we identify these with the mapping groups in the additive completion:
\begin{equation}\label{eq:mapping-group-additive-completion}
 {\cM}_{-}(-)^{\group}\simeq
 \Hom_{(L_K\nCalg_{\sepa,h})^{\mathrm{add}}}\bigl(L_{h,K}^{\add}(-),L_{h,K}^{\add}(-)\bigr).
\end{equation}

\begin{rem}\label{ku-str}We consider the maximal tensor product on $\nCalg_{\sepa}$. The localization $L_{h,K}$ 
and the additivization $L^{\add}$ refine to symmetric monoidal functors.
The tensor unit of $(L_K\nCalg_{\sepa,h})^{\mathrm{add}}$ is represented by $\C$, and we have
\[
{ku\simeq \End_{(L_K\nCalg_{\sepa,h})^{\mathrm{add}}}(\C)\ .}
\]
We get a natural refinement of ${\cM}_{-}(-)^{\group}$ to a $ku$-module valued functor. 
\hB
\end{rem}

Write
\[
 E(D,B):=\operatorname{map}_{\EE_{\sepa}}
       (e_{\sepa}(D),e_{\sepa}(B))\in\Sp
\]
for the mapping spectrum in $E$-theory and set $E_j(D,B):=\pi_jE(D,B)$.
Its underlying mapping group is
\begin{equation}\label{eq:E-mapping-group}
 \Hom_{\EE_{\sepa}}(e_{\sepa}(D),e_{\sepa}(B))
 \simeq\Omega^\infty E(D,B)
 \quad\text{in }\CGrp(\Ani).
\end{equation}

\begin{rem}
There is also a stable $\infty$-category $\EE$ and an
$E$-theory functor $e:\nCalg\to\EE$ defined for all
$C^*$-algebras \cite[Definition~8.6]{orbita}. The canonical
inclusion $\EE_{\sepa}\hookrightarrow\EE$ is fully faithful and
satisfies $e|_{\nCalg_{\sepa}}\simeq e_{\sepa}$, with the inclusion
understood. Consequently, for separable $D$ and $B$,
\[
 \operatorname{map}_{\EE_{\sepa}}
       (e_{\sepa}(D),e_{\sepa}(B))
 \simeq\operatorname{map}_{\EE}(e(D),e(B)),
\]
so the notation $E(D,B)$ is compatible with the notation for the mapping spectrum in the
larger category \cite[Corollary~8.2(1) and equation~(8.6)]{orbita}. \hB
\end{rem}

The additive functor $e_{\mathrm{add}}$ of
\eqref{eq:additive-E-theory-factorization} induces a morphism of
mapping group functors. Using \eqref{eq:mapping-group-additive-completion}
and \eqref{eq:E-mapping-group}, we obtain the natural transformation
\begin{equation}\label{eq:intro-first-comparison}
 \eta:{\cM}_{-}(-)^{\group}\longrightarrow\Omega^\infty E(-,-)
  \colon 
 (\nCalg_{\sepa})^{\op}\times\nCalg_{\sepa}\longrightarrow \CGrp(\Ani).
\end{equation}
Its component $\eta_{D,B}$ is the comparison map \eqref{gwerrwegw}.
Our first question is when this map is an equivalence.

\begin{ex}\label{ex:intro-projections}
For $D=\C$ and $B$ in $\Calg_{\sepa}$ (i.e., $B$ is unital),   we have an equivalence:
\[
 {\cM}_\C(B)^{\group}\xrightarrow[\simeq]{\eta_{\C,B}}\Omega^\infty E(\C,B)
 \simeq\Omega^\infty K^{\topp}(B).
\]
This is the usual 
 group completion description
of operator $K$-theory
\cite[Proposition~10.9]{orbita}.

Unitality of the coefficient algebra matters at this point. Consider for example the non-unital algebra  $B=C_0(\R^2)$. The algebra $C_0(\R^2)\otimes\bbK$ has no
nonzero projections. Hence
\[
 {\cM}_\C(C_0(\R^2))^{\group}\simeq0,
 \qquad  \Omega^\infty K^{\topp}(C_0(\R^2))\simeq \Omega^{\infty+2}KU.
\] \hB
 \end{ex}
Applying $e_{\mathrm{add}}$ to mapping groups extends the comparison \eqref{eq:intro-first-comparison} in the first variable to $(L_K\nCalg_{\sepa,h})^{\add}$. We also denote its components by $\eta_{D,B}$ for $D$ in this additive completion and $B$ in $\nCalg_{\sepa}$.
\begin{ddd}\label{def:eta-admissible}
An object $D$ of the additive completion is
\emph{$\eta$-admissible} if $\eta_{D,B}$ is an equivalence for every
$B$ in $\Calg_{\sepa}$. Write
\[
 (L_K\nCalg_{\sepa,h})^{\add}_{\eta}\subseteq (L_K\nCalg_{\sepa,h})^{\add}
\]
for the full subcategory of these objects. We also call a
separable $C^*$-algebra $\eta$-admissible if its image belongs to
this subcategory.
\end{ddd}

Both mapping group functors take finite products
in the source variable to finite products, so $\eta$-admissible
objects are closed under finite products and retracts. On the level of $C^{*}$-algebras $\eta$-admissibility is also homotopy and Morita invariant.

\begin{ex}\label{ex:additive-free-products}
We know by Example~\ref{ex:intro-projections} that $\C$ is $\eta$-admissible. Another $\eta$-admissible example is $C_{0}(\R)$. Indeed, 
 the monoid ${\cM}_{C_0(\R)}(B)$ is already group-like and is the
 commutative group  in anima of stable unitaries of $B$. For every   $B$ in $\Calg_{\sepa}$, the
comparison is a known equivalence
\[
 {\cM}_{C_0(\R)}(B) \simeq {\cM}_{C_0(\R)}(B)^{\group}
 \stackrel{\eta_{C_{0}(\R),B}}{\simeq }\Omega^\infty E(C_0(\R),B)
 \simeq\Omega^{\infty-1}K^{\topp}(B)
\]
\cite[Example~4.9, Proposition~9.4 and Lemma~9.5]{orbita}.
This also follows from
Proposition~\ref{prop:intro-semiprojective}, since $C_0(\R)$
is nuclear, quasidiagonal, and semiprojective.
  \hB
\end{ex}
 
 \begin{ex}\label{ex:additive-free-products1}
For separable $C^*$-algebras $A$ and $B$, the canonical
map $A*B\to A\oplus B$ is an equivalence after applying
$L_{h,K}$, and hence also in the additive completion.
Consequently, finite free products of $\eta$-admissible algebras are again $\eta$-admissible. 
 In particular,  
$\C*\C$ and $\C*\C*\C$ are $\eta$-admissible.
  \hB
\end{ex}
 \begin{ex}
\label{ex:additive-free-products2}
$\eta$-admissibility is also preserved by  
 unital free products of
augmented algebras. Let $A$ and $B$ be separable unital
$C^*$-algebras with augmentations
$\epsilon_A:A\to\C$ and $\epsilon_B:B\to\C$.
The assertion follows from the claim that $A*_{\C} B$ is a retract of $A*B$ in $(L_{K}\nCalg_{\sepa,h})^{\add}$.
 The retraction is given by   the canonical map
\[
 q:A*B\longrightarrow A*_{\C}B.
\]
 Let
$i_A:A\to A*_{\C}B$ and $i_B:B\to A*_{\C}B$ be the canonical maps, and
let $\eta_A,\eta_B,\eta$ denote the unit maps of $A$, $B$, and $A*_{\C}B$, respectively.
The augmentations define homomorphisms
\[
 f_A:A*_{\C}B\longrightarrow A,\qquad
 f_Ai_A:=\id_A,\quad f_Ai_B:=\eta_A\epsilon_B,
\]
\[
 f_B:A*_{\C}B\longrightarrow B,\qquad
 f_Bi_A:=\eta_B\epsilon_A,\quad f_Bi_B:=\id_B,
\]
and an augmentation $\epsilon:A*_{\C}B\to\C$.
Write $[f]$ for the morphism induced by $f$ in the homotopy category of the additive completion. We show that
\begin{equation}\label{eq:augmented-free-product-unit-correction}
 [i_A][f_A]+[i_B][f_B]=\id_{A*_{\C}B}+[\eta\epsilon].
\end{equation}
For $t$ in $[0,\pi/2]$, set
\[
 R(t):=\begin{pmatrix}\cos t&-\sin t\\ \sin t&\cos t\end{pmatrix}.
\]
The homomorphisms to $M_2(A*_{\C}B)$ given by
\[
 a\longmapsto\operatorname{diag}(i_A(a),\epsilon_A(a)1_{A*_{\C}B}),
 \qquad
 b\longmapsto
 R(t)\operatorname{diag}(i_B(b),\epsilon_B(b)1_{A*_{\C}B})R(t)^*
\]
are unital, so they define a homomorphism
\[
 H:A*_{\C}B\longrightarrow C([0,\pi/2],M_2(A*_{\C}B)).
\]
Its endpoints are
\[
 H_0=\operatorname{diag}(\id_{A*_{\C}B},\eta\epsilon),
 \qquad H_{\pi/2}=\operatorname{diag}(i_Af_A,i_Bf_B).
\]
Homotopy invariance and matrix stability prove
\eqref{eq:augmented-free-product-unit-correction}.
Under $A*B\simeq A\oplus B$ in the additive completion,
$[q]$ has components $[i_A],[i_B]$. Since
$\epsilon_Af_A=\epsilon$, the morphism
\[
 s:=\bigl([f_A]-[\eta_A][\epsilon],[f_B]\bigr):
 A*_{\C}B\longrightarrow A\oplus B\simeq A*B
\]
satisfies
\[
 \begin{aligned}
 [q]s
 &= [i_A]\bigl([f_A]-[\eta_A][\epsilon]\bigr)+[i_B][f_B]\\
 &= [i_A][f_A]+[i_B][f_B]-[\eta\epsilon]\\
 &= \id_{A*_{\C}B},
 \end{aligned}
\]
where we used $i_A\eta_A=\eta$ and
\eqref{eq:augmented-free-product-unit-correction}. Thus $A*_{\C}B$ is a retract of $A\oplus B$
in the additive completion, proving the claim.
\hB\end{ex}
 \begin{ex}
\label{ex:additive-free-products3}
We apply Example \ref{ex:additive-free-products2}  to the full group $C^*$-algebras with the
augmentations given by the trivial representations. The algebra 
 $C^*_{\max}(F_1)\cong C(S^1)$ is $\eta$-admissible by Proposition~\ref{prop:intro-semiprojective},
since $C(S^1)$ is nuclear, quasidiagonal, and semiprojective.
The isomorphisms
\[
 C^*_{\max}(F_{r+1})
 \cong C^*_{\max}(F_r)*_{\C}C^*_{\max}(F_1)
\]
then show by induction that   $C^*_{\max}(F_r)$ is $\eta$-admissible for every $r$ in 
  $\nat_{\geq1}$.
 \hB
\end{ex}

\subsection{Unitalization and the reduced construction}
\label{sec:reduced-construction}
 {The failure of the projection model \eqref{z345t345z354z35} for non-unital algebras (Example~\ref{ex:intro-projections}) and its correction by reduction in \eqref{jrzjzthtzhrthtz} motivate  application of  the reduced construction to ${\cM}_D$ for general $D$.} 

In the following we introduce the reduced construction using categorical language. This will become relevant
 when we introduce $j^{\sfred}$ in \eqref{gerfwergwergwer} below.
 Let $\cC$ be a pointed $\infty$-category admitting fibres of
split epimorphisms. We define the endofunctor and  
natural transformation to the identity
\[
 (-)^{\red}:\Fun(\nCalg_{\sepa},\cC)\longrightarrow
 \Fun(\nCalg_{\sepa},\cC),\qquad (-)^{\red}\longrightarrow\id
\]  which sends $G$ in the functor category to the universal reduced functor
  $G^{\red}\to G$ over $G$.  Thus for each reduced $T$, composition with the structure map induces
an equivalence
\begin{equation}\label{eq:reduction-universal-property}
 {\operatorname{Nat}(T,G^{\red})\xrightarrow{\simeq}
 \operatorname{Nat}(T,G).}
\end{equation}
The reduction exists and is given by
\begin{equation}\label{eq:reduction-fibre-formula}
 {G^{\red}(B)\simeq\Fib_0\bigl(G(B)\to G(0)\bigr).}
\end{equation}
The map is induced by $B\to0$, and the fibre exists since this zero map has  a split   $0\to B$.

If $\cC$ is additive, the section gives the natural splitting
\begin{equation}\label{eq:reduction-natural-splitting}
 {G\simeq G^{\red}\oplus\operatorname{const}G(0).}
\end{equation}
The projection  onto the first summand also induces
\[
{\operatorname{Nat}(G^{\red},T)\stackrel{\simeq}{\to}\operatorname{Nat}(G,T)}
\]
for all  reduced $T$.

We have the adjunction \[
{(-)^{u}:  \nCalg_{\sepa}  \leftrightarrows  \Calg_{\sepa} :i\ ,}
\]
where $i$ is the inclusion and $(-)^{u}$ is unitalization.
The unit of this adjunction at $A$ in $\nCalg_{\sepa}$  is  the  unitalization inclusion  \begin{equation}\label{eq:unitalization-inclusion}
 \iota_A:A\longrightarrow A^u.
\end{equation}
 {Applying unitalization to $A\to0$ gives the scalar projection
\begin{equation}\label{eq:unitalization-scalar-projection}
 \epsilon_A:A^u\longrightarrow\C.
\end{equation}}
We define the endofunctor  
\[(-)^{\fred}:\Fun(\nCalg_{\sepa},\cC)\to \Fun(\nCalg_{\sepa},\cC)\ , \quad 
 F\longmapsto F^{\fred}:=
       \bigl(\operatorname{Ran}_i(i^*F)\bigr)^{\red}.
\]
Since unitalization is left adjoint to $i$, the right Kan
extension is $B\mapsto F(B^u)$, and its augmentation is induced
by {the scalar projection \eqref{eq:unitalization-scalar-projection}}. Thus we can calculate the reduced construction by the usual formula
\[
 F^{\fred}(B)\simeq\Fib_0\bigl(F(B^u)\to F(\C)\bigr).
\]
Applying reduction to the unit $F\to\operatorname{Ran}_i(i^*F)$
gives a natural transformation $F^{\red}\to F^{\fred}$. If $F$ is
reduced, we compose it with the inverse of the reduction map
\[
 F\xrightarrow{\simeq}F^{\red}\longrightarrow F^{\fred}.
\]
If $\cC$ is additive, {the splitting \eqref{eq:reduction-natural-splitting}} gives a projection
for every $F$. We then use the composite
\[
 F\xrightarrow{\mathrm{pr}} F^{\red}\longrightarrow F^{\fred}.
\]
If $F$ sends the unitalization sequences to fibre sequences,
then $F$ is reduced and 
the map
$F\to F^{\fred}$ is then an equivalence.
  If $F$
also inverts corner inclusions, the composite
\begin{equation}\label{hrtegertger}
 F(-)\longrightarrow F(-\otimes\bbK)
                 \longrightarrow F^{\fred}(-\otimes\bbK)=F^{\sfred}(-)
\end{equation}is a natural equivalence as well. Here $F^{\sfred}$ is the
stabilized reduction defined in \eqref{gewrfrefwrf}.

Let $F:\nCalg_{\sepa}\to\cC$ be a functor to an additive $\infty$-category admitting fibres of split epimorphisms.
\begin{lem}\label{lem:stabilized-reduction-properties}
If $F$ is homotopy invariant, inverts corner inclusions, and preserves finite products, then $F^{\sfred}$ has these properties.
\end{lem}
\begin{proof}
Unitalization preserves homotopies, and corner inclusions between stabilized algebras are homotopy equivalences. Thus $F^{\sfred}$ is homotopy invariant and inverts corner inclusions.

For additivity, it suffices to show that $F^{\fred}$ preserves finite direct sums, since stabilization does so. Let $A$ and $B$ be in $\nCalg_{\sepa}$. Their inclusions and projections induce homomorphisms
\begin{equation}\label{eq:reduced-additivity-inclusions-projections}
 i_A:A^u\longrightarrow(A\oplus B)^u,\qquad
 p_A:(A\oplus B)^u\longrightarrow A^u,
\end{equation}
and similarly $i_B,p_B$. Write $\eta_A,\eta_B,\eta$ for the unit maps of $A^u,B^u,(A\oplus B)^u$, respectively, and $\epsilon_A,\epsilon_B,\epsilon$ for their scalar projections to $\C$. We have
\begin{equation}\label{eq:reduced-additivity-cross-composites}
 p_Ai_A=\id_{A^u},\quad p_Bi_B=\id_{B^u},\quad
 p_Ai_B=\eta_A\epsilon_B,\quad p_Bi_A=\eta_B\epsilon_A.
\end{equation}
We next prove the identity
\begin{equation}\label{eq:reduced-additivity-rotation-identity}
 F(i_A)F(p_A)+F(i_B)F(p_B)
 =\id_{F((A\oplus B)^u)}+F(\eta\epsilon)
 \quad\text{in }\ho(\cC).
\end{equation}
We view $A$ and $B$ as orthogonal ideals in $(A\oplus B)^u$. For $s$ in $[0,\pi/2]$, set
\begin{equation}\label{eq:reduced-additivity-rotation-homotopy}
 R_s:=\begin{pmatrix}\cos s&-\sin s\\\sin s&\cos s\end{pmatrix},\qquad
 H_s(a,b,\lambda):=
 \begin{pmatrix}a&0\\0&0\end{pmatrix}
 +R_s\begin{pmatrix}b&0\\0&0\end{pmatrix}R_s^*+\lambda1_2.
\end{equation}
Orthogonality of the two ideals ensures that the formula defines a homotopy of homomorphisms from $(A\oplus B)^u$ to $M_2((A\oplus B)^u)$. Its endpoints are
\begin{equation}\label{eq:reduced-additivity-rotation-endpoints}
 H_0=\operatorname{diag}(\id_{(A\oplus B)^u},\eta\epsilon),\qquad
 H_{\pi/2}=\operatorname{diag}(i_Ap_A,i_Bp_B).
\end{equation}
 {By \eqref{eq:additive-completion-universal-property}, $F$ factors through an additive functor on $(L_K\nCalg_{\sepa,h})^{\add}$. Consequently, after identifying matrix targets by corner inclusions, $F$ sends block sums to sums of morphisms. Applying $F$ to the homotopy with endpoints \eqref{eq:reduced-additivity-rotation-endpoints} therefore gives \eqref{eq:reduced-additivity-rotation-identity}.} 
 
 We now use this identity in $\ho(\cC)$. 
For $X$ in $\nCalg_{\sepa}$, the scalar projection and its unit section give the natural splitting
\begin{equation}\label{eq:reduced-additivity-scalar-splitting}
 F(X^u)\simeq F^{\fred}(X)\oplus F(\C).
\end{equation}
For $X$ in $\{A,B\}$, write $\bar i_X:F^{\fred}(X)\to F^{\fred}(A\oplus B)$ and $\bar p_X:F^{\fred}(A\oplus B)\to F^{\fred}(X)$ for the induced maps on the reduced summands. Since $i_X$ and $p_X$ preserve scalar projections and unit sections, their matrices with respect to \eqref{eq:reduced-additivity-scalar-splitting} are
\begin{equation}\label{eq:reduced-additivity-block-matrices}
 F(i_X)=\begin{pmatrix}\bar i_X&0\\0&\id_{F(\C)}\end{pmatrix},\qquad
 F(p_X)=\begin{pmatrix}\bar p_X&0\\0&\id_{F(\C)}\end{pmatrix},\qquad
 F(\eta\epsilon)=\begin{pmatrix}0&0\\0&\id_{F(\C)}\end{pmatrix}.
\end{equation}
Equation~\eqref{eq:reduced-additivity-rotation-identity} therefore reads
\begin{equation}\label{eq:reduced-additivity-block-identity}
 \begin{pmatrix}
 \bar i_A\bar p_A+\bar i_B\bar p_B&0\\0&2\id_{F(\C)}
 \end{pmatrix}
 =\begin{pmatrix}
 \id_{F^{\fred}(A\oplus B)}&0\\0&2\id_{F(\C)}
 \end{pmatrix}.
\end{equation}
Together with \eqref{eq:reduced-additivity-cross-composites}, this gives the biproduct identities in $\ho(\cC)$. Hence the canonical map
\begin{equation}\label{eq:reduced-additivity-biproduct}
 F^{\fred}(A\oplus B)\xrightarrow{\simeq}
 F^{\fred}(A)\oplus F^{\fred}(B)
\end{equation}
is an equivalence. Since $F^{\fred}(0)\simeq0$, this proves additivity.
\end{proof}

This completes the description of the stabilized reduction which goes into the definition of the reduced representation group bifunctor \begin{equation}\label{eq:reduced-representation-group-bifunctor}
  {\cT_{-}(-):= ({\cM}_{-}(-)^{\group})^{\sfred}} :(\nCalg_{\sepa})^{\op}\times \nCalg_{\sepa}\to \CGrp(\Ani)\  ,
\end{equation}
where the stabilized reduction is applied to the second variable.
Specializing the first variable at $D$ we obtain 
 $\cT_D$ from  Definition~\ref{def:intro-T}.
Applying stabilized reduction to \eqref{eq:mapping-group-additive-completion} gives
\begin{equation}\label{gwerfwergwerfw}
 \cT_{-}(-)\simeq
 \Bigl(\Hom_{(L_K\nCalg_{\sepa,h})^{\add}}
 \bigl(L^{\add}_{h,K}(-),L^{\add}_{h,K}(-)\bigr)\Bigr)^{\sfred},
\end{equation}
where stabilized reduction acts on the second variable of the mapping-group functor.

The functor $\Omega^\infty E(D,-)$ sends unitalization
sequences to fibre sequences and inverts corner inclusions.
Applying stabilized reduction in the  second variable to $\eta$ in \eqref{eq:intro-first-comparison}
and then using the inverse of  \eqref{hrtegertger} for the target functor gives the natural transformation
\begin{equation}\label{eq:intro-tau}
 \tau:\cT
 \xrightarrow{\eta^{\sfred}}
 \bigl(\Omega^\infty E(-,-)\bigr)^{\sfred}
 \xrightarrow{\simeq}\Omega^\infty E(-,-).
\end{equation}
Specializing the first variable at $D$ we obtain the transformation $\tau_{D}$ in \eqref{gwerrwegw1}.


 \begin{ddd}\label{def:intro-semiprojective}
A  morphism $u:P\to Q$ in $\nCalg_{\sepa}$ is \emph{semiprojective} if, for every
  $C$  in $\nCalg$, increasing sequence of closed ideals
$(J_k)_{k\in \nat}$ of $C$ with $J:=\overline{\bigcup_{k\in \nat}J_k}$, and homomorphism
$\phi:Q\to C/J$, there are $k$ in $\nat$ and a lift in the diagram
\[
 \xymatrix@C=4em@R=2.5em{
 P\ar[r]^u\ar@{..>}[d]&Q\ar[d]^\phi\\
 C/J_k\ar[r]&C/J.
 }
\]
An algebra $D$ is \emph{semiprojective} if $\id_D$ is semiprojective.
\end{ddd}
{The notion of a semiprojective morphism is due to Blackadar \cite[Definition~2.10]{blackadar-shape}. For background on semiprojective $C^*$-algebras and their application to lifting problems, see Loring's book \cite{loring-lifting}, especially Section~14.1 on stable relations.}

\begin{ddd}\label{def:intro-quasidiagonal}
A separable $C^*$-algebra $D$ is \emph{quasidiagonal} if it has a
faithful representation $\pi:D\to\mathcal L(H)$ on a separable
Hilbert space and a sequence of finite-rank projections $(p_k)_{k\in\nat}$ converging strongly
to $1$ such that for every $a$ in $D$
\[
\lim_{k\to\infty}\|[p_k,\pi(a)]\|=0.
\]
\end{ddd}
{Halmos studied quasidiagonality for operators on Hilbert space in \cite[Section~4]{halmos-ten-problems}. Javier \cite{javier-quasidiagonal} subsequently treated the corresponding property of $C^*$-algebras. We refer to \cite[Chapter~7]{brown-ozawa} for a systematic account.}

By construction, the dependence of the 
  comparison $\eta$ and $\tau$   in the first variable factorizes over $L_{h,K}^{\add}$.
  In particular, the property that $\tau_{D}$ is an equivalence only depends on $L_{h,K}^{\add}(D)$ in $(L_{K}\nCalg_{\sepa,h})^{\add}$.

\begin{ddd}\label{def:tau-admissible}
An object $D$ in $(L_{K}\nCalg_{\sepa,h})^{\add}$   is \emph{$\tau$-admissible} if $\tau_D$ is an equivalence.
Write
\[
 (L_K\nCalg_{\sepa,h})^{\add}_{\tau}\subseteq (L_K\nCalg_{\sepa,h})^{\add}
\]
for the full subcategory of these objects. We call a separable
$C^*$-algebra $\tau$-admissible if its image belongs to this
subcategory.
\end{ddd}

The $\tau$-admissible objects are closed under
finite direct sums and retracts. 

\begin{ddd}\label{def:additive-nqs}
We define $(L_K\nCalg_{\sepa,h})^{\add}_{\mathrm{nqs}}$ as the
smallest full subcategory of $(L_K\nCalg_{\sepa,h})^{\add}$
that contains the images of nuclear, quasidiagonal, semiprojective
separable $C^*$-algebras and is closed under equivalences,
finite products, and retracts.
\end{ddd}

Recall the full subcategory of $\eta$-admissible objects from Definition~\ref{def:eta-admissible}.
The following result shows that nuclear, quasidiagonal, semiprojective algebras are $\eta$-admissible, and that $\eta$-admissibility implies $\tau$-admissibility.

\begin{prop}\label{prop:intro-semiprojective}
We have the containments
\[
 (L_K\nCalg_{\sepa,h})^{\add}_{\mathrm{nqs}}
 \subseteq (L_K\nCalg_{\sepa,h})^{\add}_{\eta}
 \subseteq (L_K\nCalg_{\sepa,h})^{\add}_{\tau}.
\]
\end{prop}

The proof is given in
Section~\ref{sec:comparison-components}.

{The $\tau$-admissible class has the same closure properties as
the $\eta$-admissible class (see Definition~\ref{def:eta-admissible} and
Examples~\ref{ex:additive-free-products1} and~\ref{ex:additive-free-products2}).}
It therefore contains algebras that do not satisfy these three assumptions simultaneously, {including e.g. $C^*_{\max}(F_r)$ for $r\geq 2$
(Example~\ref{ex:additive-free-products3})}.

\begin{ex}\label{ex:intro-oinfinity}
The Cuntz algebra $\mathcal O_\infty$ is nuclear and
semiprojective \cite[Corollary~4.6]{enders}. Its image in the
additive completion is neither $\eta$-admissible nor
$\tau$-admissible. Indeed, every homomorphism from
$\mathcal O_\infty$ to either $\bbK$ or $\bbK^u\otimes\bbK$
vanishes. A nonzero homomorphism would send the unit to a
properly infinite projection, and both targets are stably finite.
Consequently,
\[
 {\cM}_{\mathcal O_\infty}(\C)^{\group}\simeq 0,
 \qquad
 \cT_{\mathcal O_\infty}(\C)\simeq 0.
\]
The unit inclusion $\C\to\mathcal O_\infty$ is an
$E$-theory equivalence, so both $\eta_{\mathcal O_\infty,\C}$
and $\tau_{\mathcal O_\infty,\C}$ have a non-zero target  
$\Omega^\infty KU$. \hB
\end{ex}

The functor $\cT_D$ is homotopy invariant and inverts corner inclusions.
The following examples show that exactness need not hold
without the assumptions of Proposition~\ref{prop:intro-semiprojective}.

\begin{ex}\label{ex:intro-nonexact}
{We show that $\cT_{\mathcal O_\infty}$ fails exactness on the
Cuntz--Toeplitz extension}
\[
 0\longrightarrow\bbK\longrightarrow\mathcal E_2
   \longrightarrow\mathcal O_2\longrightarrow0,
\]
where $\mathcal E_2$ is generated by two isometries with orthogonal
ranges. Its $K_0$-group is $\Z$, generated by $[1_{\mathcal E_2}]$
\cite{cuntz-k}.

{By Example~\ref{ex:intro-oinfinity}, every homomorphism
$\mathcal O_\infty\to(B\otimes\bbK)^u\otimes\bbK$ has zero scalar quotient
and hence factors through $B\otimes\bbK\otimes\bbK$. Thus}
\[
 \cT_{\mathcal O_\infty}(B)
 \simeq {\cM}_{\mathcal O_\infty}(B)^{\group},
 \qquad\cT_{\mathcal O_\infty}(\bbK)\simeq0.
\]
The properly infinite unit of $\mathcal E_2$ gives a unital
homomorphism $\mathcal O_\infty\to\mathcal E_2$. Its class survives
group completion: evaluation on the unit maps it to the nonzero
class $[1_{\mathcal E_2}]$ in $K_0(\mathcal E_2)$.
On the other hand, $\cT_{\mathcal O_\infty}$ is homotopy
invariant and additive, inverts corner inclusions, and takes values in
commutative groups in anima. For
$A:=\mathcal O_2\otimes\bbK$, the identity and its block sum
with itself are homotopic: both maps are full
and have the same, zero, $KK$-class, so the uniqueness theorem
\cite[Theorem~B]{gabe-classification} makes them asymptotically
unitarily equivalent. Asymptotic unitary equivalence implies homotopy here; see \cite[Remark~8.7]{gabe-classification}. Thus
$\id_A+\id_A\simeq\id_A$ in $L_K\nCalg_{\sepa,h}$.
Applying the additive functor $\cT_{\mathcal O_\infty}$ gives
\[
 \id_{\cT_{\mathcal O_\infty}(A)}
 +\id_{\cT_{\mathcal O_\infty}(A)}
 \simeq\id_{\cT_{\mathcal O_\infty}(A)}.
\]
Since the target is $\CGrp(\Ani)$, cancellation gives
$\id_{\cT_{\mathcal O_\infty}(A)}\simeq0$ and hence
$\cT_{\mathcal O_\infty}(\mathcal O_2)\simeq0$.
The resulting sequence on component groups
\[
 0=\pi_0\cT_{\mathcal O_\infty}(\bbK)
 \longrightarrow\pi_0\cT_{\mathcal O_\infty}(\mathcal E_2)
 \longrightarrow\pi_0\cT_{\mathcal O_\infty}(\mathcal O_2)=0
\]
fails to be exact at its middle term, since
$\pi_0\cT_{\mathcal O_\infty}(\mathcal E_2)\ne0$.  \hB
\end{ex}


\begin{ex}\label{ex:uhf-unresolved-nonexact}
Let $D:=M_{2^\infty}$ be the UHF algebra of type $2^\infty$.
It is separable, nuclear and quasidiagonal, but not semiprojective
\cite[Example~6.1]{phillips-sorensen-thiel}.
We show that $\cT_D$ is not exact.

The algebra $D$ has no nonzero finite-dimensional
representations. Hence every homomorphism $D\to\bbK$ is zero.
For every $B$, a homomorphism
$D\to(B\otimes\bbK)^{u}\otimes\bbK$ has zero scalar quotient and
therefore factors through $B\otimes\bbK\otimes\bbK$. Consequently,
\[
 \cT_D(B)\simeq {\cM}_D(B)^{\group}.
\]
Put $SD:=C_0((0,1),D)$ and $CD:=C_0([0,1),D)$. The suspension
$SD\otimes\bbK$ has no nonzero projections, so every
homomorphism from the unital algebra $D$ to it is zero.
The cone $CD$ is contractible. Thus
\[
 \cT_D(SD)\simeq0,\qquad\cT_D(CD)\simeq0.
\]

On the other hand, let $p:S^2\to M_2(\C)$ be the rank-one
Bott projection. The continuous family of homomorphisms
\[
 D\longrightarrow D\otimes M_2(\C)\subseteq D\otimes\bbK,
       \qquad a\longmapsto a\otimes p(z)\quad(z\in S^2)
\]
gives, after group completion and translation by its value at
a base point, a class in $\pi_2\cT_D(D)$. Evaluation on $1_D$
induces a map $\cT_D(D)\to\Omega^\infty K^{\topp}(D)$ sending
this class to the Bott class corresponding to
$[1_D]\ne0$ in $K_0(D)$. In particular,
\[
 \pi_2\cT_D(D)\ne0.
\]
If $\cT_D$ were exact, the cone extension
\[
 0\longrightarrow SD\longrightarrow CD
       \xrightarrow{\operatorname{ev}_0}D\longrightarrow0
\]
would imply
$\cT_D(SD)\simeq\Omega\cT_D(D)$, contradicting the preceding
vanishing and nonvanishing statements. \hB
\end{ex}

\subsection{The resolved representation functor}
{The additive category $\mathcal U$ admits sequential colimits by
Proposition~\ref{prop:universal-unreduced-category}, so it is idempotent
complete and admits fibres of split epimorphisms. We can therefore apply
the stabilized reduced construction of Section~\ref{sec:reduced-construction}
to the functor $j$ from Definition~\ref{def:universal-admissible} to obtain}
\begin{equation}\label{gerfwergwergwer}j^{\sfred}:\nCalg_{\sepa}\to \cU\ .
\end{equation}
\begin{lem}\label{lem:stabilized-reduction-additively-continuous}
{The functor $j^{\sfred}$ is additively continuous.}
\end{lem}
\begin{proof}
{Lemma~\ref{lem:stabilized-reduction-properties} gives homotopy invariance, inversion of corner inclusions, and preservation of finite products. Unitalization, stabilization, and $j$ preserve countable filtered colimits, and fibres of split epimorphisms commute with these colimits in the additive category $\mathcal U$. Thus $j^{\sfred}$ also preserves countable filtered colimits.}
\end{proof}   
{The functor $j^{\sfred}$ enters the definition of the resolved stable reduced representation group bifunctor from Definition~\ref{kohprtgerthertgertgertge}:}
\begin{equation}\label{fqwedweqfqwed}
{\cR}_{-}(-):=\Hom_{\cU}(j(-),j^{\sfred}(-)):(\nCalg_{\sepa})^{\op}\times\nCalg_{\sepa}\to\CGrp(\Ani)\ .
\end{equation}
{Formula \eqref{gwerfwergwerfw} defines $\cT$ by applying stabilized reduction to the mapping-group bifunctor in its second variable, whereas in \eqref{fqwedweqfqwed} it is applied to $j$ inside the second entry of $\Hom_{\cU}$. These constructions agree naturally because the mapping-group bifunctor preserves limits in its second variable, in particular the fibres of split epimorphisms used in stabilized reduction.}
{Applying stabilized reduction to the second variable in the map between $\Hom$-bifunctors induced by $j_{\mathrm{add}}$ from \eqref{gergewrfewrgwerge} yields a natural transformation of bifunctors} 
\begin{equation}\label{eq:T-to-R}
 {\rho:\cT \longrightarrow {\cR} :(\nCalg_{\sepa})^{\op}\times\nCalg_{\sepa}
       \longrightarrow\CGrp(\Ani).}
\end{equation}
{Both sides of \eqref{eq:T-to-R} factor through the additive completion in the first variable.}
\begin{ddd}\label{def:rho-admissible}
An object $D$ in $(L_K\nCalg_{\sepa,h})^{\add}$ is \emph{$\rho$-admissible} if
\begin{equation}\label{eq:rho-admissible-comparison}
 \rho_D:\cT_D\longrightarrow {\cR}_D
 \quad\text{in }\Fun(\nCalg_{\sepa},\CGrp(\Ani))
\end{equation}
is an equivalence. We write
\[
 (L_K\nCalg_{\sepa,h})^{\add}_{\rho}
 \subseteq (L_K\nCalg_{\sepa,h})^{\add}
\]
for the full subcategory of these objects. We call a separable $C^*$-algebra $\rho$-admissible if its image belongs to this subcategory. 
\end{ddd}

We now introduce shape systems
to compute the  values ${\cR}_{D}(B)$ explicitly.

\begin{ddd}\label{def:intro-shape-system}
A \emph{shape system} $(D_n)$ for a separable
$C^*$-algebra $D$ is a sequence of separable $C^*$-algebras
\[
 D_0\xrightarrow{\alpha_{0,1}}D_1
 \xrightarrow{\alpha_{1,2}}D_2\longrightarrow\cdots,
 \qquad\colim_nD_n\cong D,
\]
whose transition maps are semiprojective. We write
$\alpha_{nm}:D_n\to D_m$ and $d_n:D_n\to D$ for its
structure maps.
\end{ddd}

Every separable $C^*$-algebra has a shape system
\cite[Theorem~4.3]{blackadar-shape}. 

\begin{rem}\label{lophehgtrheht} We recall
 the construction of the universal
homotopy invariant and countable filtered colimit-preserving functor
$a:\nCalg_{\sepa}\to\cA_{\sepa}$ from \cite[Proposition~3.54]{budu}.
{The ordinary categorical analogues of these constructions were developed by Dadarlat \cite[Section~3]{dadarlat-shape}.}
 Let
\begin{equation}\label{eq:R-ind-inclusion}
 y:\nCalg_{\sepa,h}\longrightarrow
       \Ind_{\aleph_1}(\nCalg_{\sepa,h})
\end{equation}
be the canonical inclusion. {We abbreviate}
\begin{equation}\label{eq:R-homotopy-Yoneda}
{y_h:=y\circ L_h:\nCalg_{\sepa}\longrightarrow
\Ind_{\aleph_1}(\nCalg_{\sepa,h}).}
\end{equation}
By \cite[Theorem~3.49]{budu},
there is a right Bousfield localization
\begin{equation}\label{eq:R-shape-localization}
 L:\cA_{\sepa}\ \leftrightarrows\
 \Ind_{\aleph_1}(\nCalg_{\sepa,h}):P,
\end{equation}
inverting the maps
$\colim_n {y_h(A_n)}\to {y_h(\colim_nA_n)}$ for sequences of
$C^*$-algebras, where 
  $L$ is the fully faithful inclusion. The universal functor is then given by
\begin{equation}\label{gewrfwergwergw}
 a:={P\circ y_h}:\nCalg_{\sepa}\longrightarrow\cA_{\sepa}.
\end{equation}
For a shape system $(D_n)$ of $D$, the counit of the adjunction \eqref{eq:R-shape-localization} at ${y_h(D)}$ is \begin{equation}\label{gerfwergwerfewrfwe}
 \colim_n {y_h(D_n)}\simeq L(a(D))\simeq L(P({y_h(D)}))
 \xrightarrow{\varepsilon_{{y_h(D)}}} {y_h(D)}.
\end{equation}
The identification with the colimit follows from
\cite[Theorem~3.49 and Proposition~3.58]{budu}, and
$L(a(D))\simeq L(P({y_h(D)}))$ follows from \eqref{gewrfwergwergw}.
\hB
\end{rem}

Recall the category of additive presheaves
\[
 \mathcal P=
 \Fun^\oplus\!\left(((L_K\nCalg_{\sepa,h})^{\add})^{\op},
                   \CGrp(\Ani)\right)
\]and the functor  \[
{\mathsf L:\mathcal P\to \mathcal P_{\mathrm{cont}}}
\]
in   \eqref{gerferfwerfwer}   to its full subcategory $\mathcal P_{\mathrm{cont}}$ consisting of
presheaves $F$ such that $F\circ L_{h,K}^{\add}$ sends countable
filtered colimits of $C^*$-algebras to limits, with $L_{h,K}^{\add}$
as in \eqref{eq:additive-completion-category}. 
Let $D$ be in $\nCalg_{\sepa}$ and let $(D_n)$ be a shape system for
$D$.
 \begin{lem}\label{lem:unreduced-presheaf-reflection}
 For  $F$ in $\mathcal P$ we have
\begin{equation}\label{eq:unreduced-presheaf-reflection}
 (\mathsf LF)(L_{h,K}^{\add}(D))\simeq\lim_n F(L_{h,K}^{\add}(D_n))
 \end{equation}
\end{lem}

\begin{proof}
  
  Let $
 L_K^{\add}:=L^{\add}\circ L_K 
$, the functors $y$, $P$ and $L$ be as in \eqref{eq:R-ind-inclusion} and \eqref{eq:R-shape-localization}, and
set \[
{Q:=F^{\sim}\circ L\circ P\circ y:(\nCalg_{\sepa,h})^{\op}\to \CGrp(\Ani)\ ,}
\]
where $F^{\sim}:\Ind_{\aleph_1}(\nCalg_{\sepa,h})^{\op}\to \CGrp(\Ani)$ is the unique functor
sending countable filtered colimits to limits with $F^{\sim}\circ y\simeq F\circ L_K^{\add}$.
{The counit $\varepsilon$ of \eqref{eq:R-shape-localization}, evaluated in \eqref{gerfwergwerfewrfwe}, induces the natural transformation}
\begin{equation}\label{eq:shape-reflection-components1}
 {F\circ L_K^{\add}\simeq F^{\sim}\circ y
 \xrightarrow{F^{\sim}(\varepsilon_y)}F^{\sim}\circ L\circ P\circ y=Q.}
\end{equation}
We compute $Q(L_h(D))$ using the shape system $(D_n)$ by
\begin{equation}\label{eq:shape-reflection-components}
  Q(L_h(D))= F^{\sim}(L(P({y_h(D)})))\simeq\lim_n F(L_{h,K}^{\add}(D_n)),
\end{equation}
employing 
\eqref{gerfwergwerfewrfwe} and the continuity of $F^{\sim}$.

{We proceed in four steps.}
\begin{enumerate}
  \item We verify that $Q$ descends to $\bar Q$ in $\mathcal P_{\mathrm{cont}}$
with
\begin{equation}\label{eq:shape-reflection-descent}
 Q\simeq\bar Q\circ L_K^{\add} :(\nCalg_{\sepa,h})^{\op} \to \CGrp(\Ani).
\end{equation}
\item We show that there exists a unique     natural transformation
 \[
{F\to\bar Q:((L_{K}\nCalg_{\sepa,h})^{\add})^{\op}\to \CGrp(\Ani)}
\]
whose pull-back along $L_{K}^{\add}$ is \eqref{eq:shape-reflection-components1}.
 Precomposition with this  transformation gives a map
\begin{equation}\label{eq:shape-reflection-descended-universal}
 \operatorname{Nat}(\bar Q,G)\longrightarrow\operatorname{Nat}(F,G)
  \end{equation} for all $G$ in $\mathcal P_{\mathrm{cont}}$.
To show that \eqref{eq:shape-reflection-descended-universal}
is an equivalence, we argue that it suffices to show that the
map \begin{equation}\label{eq:shape-reflection-universal}
 \operatorname{Nat}(Q,G\circ L_K^{\add})\longrightarrow
 \operatorname{Nat}(F\circ L_K^{\add},G\circ L_K^{\add}).
\end{equation} obtained by precomposition with $L_K^{\add}$ is an equivalence.
 
\item We prove that \eqref{eq:shape-reflection-universal} is an
equivalence by identifying it with the adjunction equivalence obtained from \eqref{eq:R-shape-localization}  \begin{equation}\label{eq:shape-adjunction-equivalence}
 \operatorname{Nat}(F^{\sim}\circ L,G')
 \simeq
 \operatorname{Nat}(F^{\sim},G'\circ P)
\end{equation}
for  all functors $G':\cA_{\sepa}^{\op}\to\CGrp(\Ani)$.
\item We identify $\bar Q$ with $\mathsf LF$ and deduce the
formula  \eqref{eq:unreduced-presheaf-reflection} from \eqref{eq:shape-reflection-descent} and
\eqref{eq:shape-reflection-components}.
\end{enumerate}

\noindent\textit{Step 1.}
The functor $Q\circ L_{h}$ sends countable filtered colimits   to
limits because $a=P\circ y\circ L_{h}$ and $L$ preserve these colimits and $F^{\sim}$
sends them to limits. It remains to check preservation of finite sums and inversion of corner inclusions.
Finite direct sums and finite matrix
amplifications of shape systems again give shape
systems. Hence
\[ Q(L_h(D\oplus D'))
 \simeq\lim_nF(L_{h,K}^{\add}(D_n\oplus D'_n))
 \simeq Q(L_h(D))\times Q(L_h(D')),
 \]
and finite corner embeddings induce equivalences
$Q(L_h(D))\simeq Q(L_h(D\otimes M_r(\C)))$.
The continuity of $Q\circ L_h$ shown above and $D\otimes\bbK\cong\colim_rD\otimes M_r$
then give
\[
 Q(L_h(D\otimes\bbK))
 \simeq\lim_r Q(L_h(D\otimes M_r(\C)))
 \simeq Q(L_h(D)).
\]
We can now conclude that there is a unique  
  $\bar Q$ in
$\mathcal P_{\mathrm{cont}}$ satisfying  \eqref{eq:shape-reflection-descent}.

\noindent\textit{Step 2.}
  The functor $L_K^{\add}$ is the
localization at corner inclusions followed by additive completion. Consequently, 
 precomposition
with $L_K^{\add}$ {is} fully faithful on additive presheaves. Thus
\eqref{eq:shape-reflection-components1} descends to a map $F\to\bar Q$.
Furthermore, by the universal property of $L^{\add}_{K}$, to show that \eqref{eq:shape-reflection-descended-universal}
is an equivalence, it suffices to show that the map obtained by
precomposition with $L_K^{\add}$ is an equivalence.

\noindent\textit{Step 3.}
Since $G\in\mathcal P_{\mathrm{cont}}$, the functor
$G\circ L_{h,K}^{\add}$ sends countable filtered colimits  
to limits. {By \cite[Proposition~3.54]{budu}, there exists a unique functor $G':\cA_{\sepa}^{\op}\to\CGrp(\Ani)$ sending countable filtered colimits to limits such that}
$G  \circ  L_{h,K}^{\add} \simeq  G'\circ a  $. 
 Thus
\begin{equation}\label{eq:shape-adjunction-target}
 G\circ L_K^{\add}
 \simeq G'\circ P\circ y.
\end{equation}
Using \eqref{eq:shape-adjunction-target} we identify the left-hand side of
\eqref{eq:shape-adjunction-equivalence} with the left-hand side of
\eqref{eq:shape-reflection-universal} by two precompositions:
\begin{equation}\label{eq:shape-adjunction-left-restriction}
\begin{aligned}
 \operatorname{Nat}(F^{\sim}\circ L,G')
 &\xrightarrow{P^*}
 \operatorname{Nat}(F^{\sim}\circ L\circ P,G'\circ P)\\
 &\xrightarrow{y^*}\operatorname{Nat}(Q,G\circ L_K^{\add}).
\end{aligned}
\end{equation}
 To identify the
right-hand sides, we precompose only with $y$:
\begin{equation}\label{eq:shape-adjunction-right-restriction}
 \operatorname{Nat}(F^{\sim},G'\circ P)
 \xrightarrow{y^*}
 \operatorname{Nat}(F\circ L_K^{\add},G\circ L_K^{\add}).
\end{equation}
We now  verify that the arrows in
\eqref{eq:shape-adjunction-left-restriction} and
\eqref{eq:shape-adjunction-right-restriction} are equivalences
and that they carry \eqref{eq:shape-adjunction-equivalence} to
\eqref{eq:shape-reflection-universal}.

The first arrow of \eqref{eq:shape-adjunction-left-restriction}
is an equivalence because $P$ is a localization. The arrows given
by $y^*$ in \eqref{eq:shape-adjunction-left-restriction}
and \eqref{eq:shape-adjunction-right-restriction} are equivalences
by the universal property of the ind-completion. Indeed,
$F^{\sim}$ and $G'\circ P$ send countable
filtered colimits to limits. The same holds for
$F^{\sim}\circ L\circ P$ because $L$ is a left adjoint
and $P$ preserves countable filtered colimits by
\cite[Theorem~3.49]{budu}.

The adjunction map \eqref{eq:shape-adjunction-equivalence} first
precomposes the functors with $P$ and then composes the resulting
transformation with the image of the counit of $L\dashv P$ under
the contravariant functor $F^{\sim}$. Restricting this image along $y$
gives \eqref{eq:shape-reflection-components1}, so the adjunction map
induces precisely \eqref{eq:shape-reflection-universal}.

This proves that \eqref{eq:shape-reflection-descended-universal}
is an equivalence.

\noindent\textit{Step 4.}
Since \eqref{eq:shape-reflection-descended-universal} is an equivalence
for every $G$ in $\mathcal P_{\mathrm{cont}}$, the map $F\to\bar Q$
exhibits $\bar Q$ as $\mathsf LF$. Evaluating
\eqref{eq:shape-reflection-descent} at $L_h(D)$ and using
\eqref{eq:shape-reflection-components}, we obtain
\[
 (\mathsf LF)(L_{h,K}^{\add}(D))
 \simeq\bar Q(L_{h,K}^{\add}(D))
 \simeq Q(L_h(D))
 \simeq\lim_n F(L_{h,K}^{\add}(D_n)).
\]
\end{proof}

\begin{lem}\label{lem:universal-unreduced-mapping}
For every shape system $(D_n)$ of $D$, there are equivalences
of functors
\begin{equation}\label{eq:universal-unreduced-mapping}
 \Hom_{\mathcal U}(j(D),j(-))
 \simeq\lim_nM_{D_n}(-)^{\group}:
 \nCalg_{\sepa}\longrightarrow\CGrp(\Ani),
\end{equation}
and
\begin{equation}\label{eq:R-shape-formula}
 {\cR}_D\simeq\lim_n\cT_{D_n}:
 \nCalg_{\sepa}\longrightarrow\CGrp(\Ani).
\end{equation}
\end{lem}

\begin{proof}
For $B$ in $\nCalg_{\sepa}$, we have natural equivalences 
\[
 \begin{aligned}
 \Hom_{\mathcal U}(j(D),j(B))
 &\stackrel{\eqref{gerwfgwergwerg}}{\simeq}
 \Hom_{\mathcal P_{\mathrm{cont}}}(\mathsf Lh_{L_{h,K}^{\add}(D)},\mathsf Lh_{L_{h,K}^{\add}(B)})\\
 &\stackrel{\eqref{gerferfwerfwer}}{\simeq}\Hom_{\mathcal P}(h_{L_{h,K}^{\add}(D)},\mathsf Lh_{L_{h,K}^{\add}(B)})\\
 &\stackrel{ \eqref{eq:additive-yoneda-presheaf}}{\simeq}(\mathsf Lh_{L_{h,K}^{\add}(B)})(L_{h,K}^{\add}(D))\\
 &{\overset{\eqref{eq:unreduced-presheaf-reflection}}{\simeq}}
 \lim_n h_{L_{h,K}^{\add}(B)}(L_{h,K}^{\add}(D_n))\\
 &{\overset{\eqref{eq:mapping-group-additive-completion}}{\simeq}}
 \lim_nM_{D_n}(B)^{\group}.
 \end{aligned}
\]
Since limits commute with 
 fibres, we can apply the stabilized reduction to
\eqref{eq:universal-unreduced-mapping}  to obtain \eqref{eq:R-shape-formula}.
\end{proof}

{Let $D$ be in $\nCalg_{\sepa}$.}
\begin{kor}\label{cor:semiprojective-rho-admissible}
{If $D$ is semiprojective, then $D$ is $\rho$-admissible.}
\end{kor}
\begin{proof}
{The constant system $(D)$ is a shape system for $D$.
Under \eqref{eq:R-shape-formula} for this system, the map
$\rho_D:\cT_D\to {\cR}_D$ identifies with the identity of $\cT_D$.}
\end{proof}

For a contravariant functor $F:(\nCalg_{\sepa})^{\op}\to\cA$,
we say that $F$ is additively continuous if the opposite functor
$F^{\op}:\nCalg_{\sepa}\to\cA^{\op}$ is additively continuous in the sense of
Definition~\ref{def:admissible}.

{Let $B$ be in $\nCalg_{\sepa}$.}
\begin{prop}\label{rem:R-derived-functor}
{The functor} ${\cR}_{-}(B)$ is
additively continuous, and the transformation ${\rho_{-,B}:\cT_{-}(B)\longrightarrow {\cR}_{-}(B)}$ from \eqref{eq:T-to-R}
presents its target  
 as the universal additively continuous functor under $\cT_{-}(B)$.  \end{prop}
The universal property states that  for every additively continuous functor
$F:(\nCalg_{\sepa})^{\op}\to\CGrp(\Ani)$,
precomposition induces an equivalence  
\[
 \operatorname{Nat}({\cR}_{-}(B),F)
 \xrightarrow{\simeq}
 \operatorname{Nat}(\cT_{-}(B),F).
\]

\begin{proof}
Apply $h_-\circ L_{h,K}^{\add}:{\nCalg_{\sepa}}\to \cP$ to the scalar quotient
$(B\otimes\bbK)^u\to\C$ and its canonical section. We obtain a split
epimorphism in $\mathcal P$. {Let $T$ be its fibre in $\mathcal P$.} By
Lemmas~\ref{lem:adjoin-colimits-representable-unit} and
\ref{lem:reflection-split-fibres}, the map on fibres induced by
$j_{\mathrm{add}}$ identifies with $\upsilon_T\colon T\to\incl(\mathsf LT)$.
Definitions~\ref{def:intro-T} and~\ref{kohprtgerthertgertgertge} identify
its precomposition with $L_{h,K}^{\add}$ with the comparison
\begin{equation}\label{eq:R-derived-unit-comparison}
 \bigl(T\xrightarrow{\upsilon_T}\incl(\mathsf LT)\bigr)
 \circ L_{h,K}^{\add}
 \simeq
 {\bigl(\cT_{-}(B)\xrightarrow{\rho_{-,B}}{\cR}_{-}(B)\bigr).}
\end{equation}
By \eqref{eq:additive-completion-universal-property}, precomposition
with $L_{h,K}^{\add}$ is fully faithful on additive presheaves and
identifies $\mathcal P_{\mathrm{cont}}$ with the contravariant additively continuous functors
with values in $\CGrp(\Ani)$. Hence the universal property of
$\upsilon_T$ gives the asserted universal property of
\eqref{eq:T-to-R} at $B$.
\end{proof}

The functor $e_{\mathcal U}$ {from \eqref{eq:universal-E-factorization}} induces the natural transformation of bifunctors
\begin{equation}\label{eq:R-unreduced-E-comparison}
 \Hom_{\mathcal U}(j(-),j(-))
 \xrightarrow{e_{\mathcal U}} \Omega^\infty E(-,-)
 \colon (\nCalg_{\sepa})^{\op}\times\nCalg_{\sepa}
 \longrightarrow \CGrp(\Ani).
\end{equation}
The mapping group functor $\Hom_{\mathcal U}(j(D),-)$ preserves
fibres, so stabilized reduction of the source is ${\cR}_D$.
We apply stabilized reduction in the second variable to
\eqref{eq:R-unreduced-E-comparison} and compose with the inverse of
\eqref{hrtegertger} for $F=\Omega^\infty E(D,-)$ to obtain
\begin{equation}\label{eq:intro-alpha}
 \alpha_D:{\cR}_D\longrightarrow
 \bigl(\Omega^\infty E(D,-)\bigr)^{\sfred}
 \xrightarrow{\simeq}\Omega^\infty E(D,-).
\end{equation}
This transformation is natural in both variables.
{For future reference, note that 
\begin{equation}\label{eq:tau-alpha-rho}
 \tau_D\simeq\alpha_D\circ\rho_D\ .
 \end{equation}}{Since $e_{\sepa}$ preserves countable filtered colimits by \cite[Corollary~3.45(2) and Propositions~3.67--3.68]{budu}, a shape system $(D_n)$ for $D$ gives an equivalence of functors}
\begin{equation}\label{eq:E-shape-limit-comparison}
 {\Omega^\infty E(D,-)\simeq\lim_n\Omega^\infty E(D_n,-):
 \nCalg_{\sepa}\longrightarrow\CGrp(\Ani).}
\end{equation}
{Under \eqref{eq:R-shape-formula} and \eqref{eq:E-shape-limit-comparison}, $\alpha_D$ identifies with $\lim_n\tau_{D_n}$.}
Equivalently,
Proposition~\ref{rem:R-derived-functor} characterizes
$\alpha_{-,B}$ as the unique extension of $\tau_{-,B}$,
since $\Omega^\infty E(-,B)$ is {additively continuous}.

\begin{ddd}\label{def:alpha-admissible}
{An object $V$ in $\mathcal U$ is \emph{$\alpha$-admissible} if the comparison induced by $e_{\mathcal U}$,
\begin{equation}\label{eq:alpha-admissible-comparison}
 \Hom_{\mathcal U}(V,j^{\sfred}(B))
 \longrightarrow
 \Hom_{\EE_{\sepa}}(e_{\mathcal U}(V),e_{\sepa}(B)),
\end{equation}
is an equivalence for every $B$ in $\nCalg_{\sepa}$. We write
\[
 \mathcal U_{\alpha}\subseteq\mathcal U
\]
for the full subcategory of these objects. We call a separable $C^*$-algebra $D$ $\alpha$-admissible if $j(D)$ belongs to this subcategory; equivalently, if $\alpha_D$ in \eqref{eq:intro-alpha} is an equivalence.}
\end{ddd}
{The $\alpha$-admissible objects are closed under equivalences, finite products, retracts, and countable filtered colimits.}

{The following is one of the main technical results of the present paper. It provides  a sufficient condition for the resolved representation functor to be homological in the sense of Definition~\ref{def:homological}. Let $D$ be in $\nCalg_{\sepa}$.}
\begin{theorem}\label{thm:intro-exactness}
{If $D\in\mathcal U_{\mathrm{nq}}$, then ${\cR}_D$ is homological.}
\end{theorem}
{Since ${\cR}_D$ is reduced, homotopy invariant, and inverts corner inclusions by construction, it remains to verify exactness and Schochet exactness (Definition~\ref{def:functor-properties}(\ref{def:functor-exactness}) and~(\ref{hjjtjrtzh})). For nuclear quasidiagonal $D$ (Definition~\ref{def:intro-quasidiagonal}), these properties are established in Theorem~\ref{thm:general-resolved-exactness}. The factorization of ${\cR}_{-}$ over $j$ and the closure properties defining $\mathcal U_{\mathrm{nq}}$ give the general case. The detailed proof is given in Section~\ref{sec:resolved-exactness}. Together with the results of Subsection~\ref{sec:general-universal-class}, this theorem yields Theorem~\ref{thm:overview-one-splitting} in the introduction.}

\subsection{\texorpdfstring{{Closure properties of $\mathcal U_{\mathrm{nq}}$}}{Closure properties of U_nq}}
{We now investigate the closure properties of $\mathcal U_{\mathrm{nq}}$ and $\mathcal U_{\mathrm{nqs}}$ from Definition~\ref{def:U-nq}. We use them to give examples of algebras in these classes which do not themselves satisfy all the nuclearity, quasidiagonality, and semiprojectivity properties imposed on the generators.}

\begin{lem}\label{lem:finite-dimensional-U-nqs}
{Every finite-dimensional $C^*$-algebra belongs to $\mathcal U_{\mathrm{nqs}}$.}
\end{lem}
\begin{proof}
{Finite-dimensional $C^*$-algebras are nuclear, quasidiagonal, and semiprojective, so the assertion follows from Definition~\ref{def:U-nq}.}
\end{proof}

\begin{kor}\label{cor:AF-U-nqs}
{Every separable AF algebra belongs to $\mathcal U_{\mathrm{nqs}}$.}
\end{kor}
\begin{proof}
{Every separable AF algebra is a sequential inductive limit of finite-dimensional $C^*$-algebras, so Lemma~\ref{lem:finite-dimensional-U-nqs} and closure under countable filtered colimits apply.}
\end{proof}

\begin{ex}\label{UHF}
Typical examples of separable AF algebras are:
\begin{enumerate}
\item {the compact operators $\bbK$ and their unitalization $\bbK^u$;}
\item {the  direct sums $\bigoplus_{n\in\nat}{M}_{k_n}(\C)$, for sequences $(k_n)$ of positive integers;}
\item {UHF algebras, including the CAR algebra $M_{2^\infty}\cong\colim_n{M}_{2^n}(\C)$;}
\item {$C(X)$ for compact metrizable totally disconnected spaces $X$, such as the Cantor set;}
\item $C^{*}_{\max}(G)$   for countable locally finite groups.
{Indeed, we write $G=\bigcup_nG_n$ for an increasing sequence of finite subgroups, and} \begin{equation}\label{eq:locally-finite-group-AF}
C^*_{\max}(G)\cong\colim_n C^*(G_n).
 \end{equation}{These examples include} 
the finitary symmetric group $S_\infty=\bigcup_n S_n$, the finitary alternating group $A_\infty=\bigcup_n A_n$, and $\bigoplus_n\Z/2$.  
\end{enumerate}
{By Corollary~\ref{cor:AF-U-nqs}, all these algebras belong to $\mathcal U_{\mathrm{nqs}}$.}
Note that $AF$-algebras are nuclear and quasidiagonal. They are semiprojective if and only if they are finite-dimensional {{\cite[II.8.3.17(ix)]{blackadar-operator-algebras}}.}
\hB
\end{ex}

\begin{ex}\label{ex:commutative-nqs}
{Every separable commutative $C^*$-algebra belongs to $\mathcal U_{\mathrm{nq}}$, since it is nuclear and quasidiagonal.}

{The class $\mathcal U_{\mathrm{nqs}}$ also contains the following nuclear quasidiagonal semiprojective algebras:}
\begin{enumerate}
\item {$C(X)$ for finite graphs $X$, including $C([0,1])$, $C(S^1)$, and $C(\bigvee_{j=1}^r S^1)$ for positive integers $r$. More generally, $C(X)$ is semiprojective whenever $X$ is a compact metrizable absolute neighbourhood retract of dimension at most one \cite[Theorem~1.2]{sorensen-thiel}.  
}
\item {$C_0(\R)$, which is nuclear and quasidiagonal by commutativity and semiprojective since its one-point compactification is $S^1$ \cite[Corollary~6.2]{sorensen-thiel}.}
\item {The dimension-drop algebras, for positive integers $p,q$,}
\begin{equation}\label{eq:dimension-drop-example}
{I_{p,q}:=\{f\in C([0,1],{M}_p(\C)\otimes{M}_q(\C)):
 f(0)\in{M}_p(\C)\otimes1,\ 
 f(1)\in1\otimes{M}_q(\C)\}.}
\end{equation}
{Their irreducible representations are finite-dimensional \cite{barlak-raum-dimension-drop}, so they are liminal and hence nuclear and quasidiagonal, as in the proof of Proposition~\ref{prop:liminal-U-nq}. They are semiprojective as one-dimensional NCCW complexes by \cite[Theorem~6.2.2]{eilers-loring-pedersen-nccw}.}
 
{Consequently, the Jiang--Su algebra $\cZ$ also belongs to $\mathcal U_{\mathrm{nqs}}$ since it  it is a sequential inductive limit of prime dimension-drop algebras \cite{jiang-su}. }\end{enumerate}

\hB
\end{ex}

\begin{lem}\label{lem:U-free-product-closure}
{Each of the subcategories $\mathcal U_{\mathrm{nq}}$ and $\mathcal U_{\mathrm{nqs}}$ is closed under:}
\begin{enumerate}
\item {countable free products of separable $C^*$-algebras;}
\item {countable unital free products amalgamated over $\C$ of augmented separable unital $C^*$-algebras.}
\end{enumerate}
\end{lem}
\begin{proof}
\noindent{\emph{1.} By Example~\ref{ex:additive-free-products1}, a finite free product has the same image as the finite direct sum in the additive completion. Applying $j_{\mathrm{add}}$ from \eqref{gergewrfewrgwerge} and using closure under finite products proves the finite case. Countable free products are filtered colimits of finite free products.}

\noindent{\emph{2.} Example~\ref{ex:additive-free-products2} expresses the unital free product of two augmented unital algebras as a retract of their direct sum in the additive completion. Applying $j_{\mathrm{add}}$ and using closure under finite products and retracts proves the finite case by induction. Closure under countable filtered colimits gives the countable case.}
\end{proof}

\begin{ex}\label{ex:U-nq-free-products}
{The class $\mathcal U_{\mathrm{nqs}}$ contains:}
\begin{enumerate}
\item {$\ast_I\C$ for countable sets $I$ by  Lemmas~\ref{lem:finite-dimensional-U-nqs} and~\ref{lem:U-free-product-closure};}
\item {$C^*_{\max}(F_r)$ for $r$ in $\nat_{\geq1}\cup\{\omega\}$, where $F_\omega$ is the free group on countably many generators, by Example~\ref{ex:commutative-nqs} and Lemma~\ref{lem:U-free-product-closure}, as in Example~\ref{ex:additive-free-products3};}
\item {  $C^{*}_{\max}(\mathrm{PSL}_2(\Z))$, since}
\begin{equation}\label{eq:modular-group-free-product}
{C^*_{\max}(\mathrm{PSL}_2(\Z))\cong C^*_{\max}(\Z/2*\Z/3)
\cong\C^2*_{\C}\C^3.}
\end{equation}
{The factors are finite-dimensional and augmented, so Lemmas~\ref{lem:finite-dimensional-U-nqs} and~\ref{lem:U-free-product-closure} apply. The resulting algebra is nonnuclear because $\mathrm{PSL}_2(\Z)$ is nonamenable \cite{lance-nuclear}.}
\end{enumerate}
{All these algebras are quasidiagonal, since unital full free products of separable unital quasidiagonal algebras are quasidiagonal \cite[Corollary~14]{boca-free-products}. For $\ast_I\C$, we apply this result to its unitalization $(\ast_I\C)^u\cong\ast_{i\in I,\C}\C^2$ and use that quasidiagonality passes to subalgebras.}

{The algebra $C^*_{\max}(F_r)$ is nonnuclear for $r\geq2$, including $r=\omega$. Similarly, $\ast_I\C$ is nonnuclear for $|I|\geq3$, whereas it is nuclear for $|I|\leq2$: its unitalization is $C^*_{\max}(\ast_I\Z/2)$, and this group is amenable exactly when $|I|\leq2$ \cite{lance-nuclear}.}

{The algebras $\ast_I\C$ for finite $I$, $C^*_{\max}(F_r)$ for finite $r$, and $C^*_{\max}(\mathrm{PSL}_2(\Z))$ are semiprojective, since finite free products of semiprojective algebras are semiprojective, also in the unital category. In contrast, $\ast_{\nat}\C$ and $C^*_{\max}(F_\omega)$ are not semiprojective. Indeed, semiprojectivity would make the identity factor up to homotopy through a finite free-product stage by {\cite[Proposition~II.8.3.18]{blackadar-operator-algebras}}, contradicting the infinite rank of their $K_0$- and $K_1$-groups, respectively.}
\hB
\end{ex}

{We recall that a $C^*$-algebra $D$ is \emph{liminal}, or \emph{CCR}, if every irreducible representation $\pi:D\to B(H_\pi)$ satisfies $\pi(D)=\bbK(H_\pi)$ {\cite[Definition~IV.1.3.1]{blackadar-operator-algebras}}. For unital $D$, this is equivalent to all irreducible representations being finite-dimensional.}

\begin{prop}\label{prop:liminal-U-nq}
{Every separable liminal $C^*$-algebra belongs to $\mathcal U_{\mathrm{nq}}$.}
\end{prop}
\begin{proof}
{Let $D$ be a separable liminal $C^*$-algebra. It is nuclear, since type I $C^*$-algebras are nuclear; see {\cite[IV.3.1, following Proposition~IV.3.1.3]{blackadar-operator-algebras}}. Each irreducible representation has compact image and is therefore quasidiagonal. A countable separating direct sum of irreducible representations gives a faithful quasidiagonal representation of $D$. Thus $D$ is nuclear and quasidiagonal, and the assertion follows from Definition~\ref{def:U-nq}.}
\end{proof}

\begin{ex}\label{ex:liminal-algebras}
{Proposition~\ref{prop:liminal-U-nq} applies to the following separable liminal algebras:}
\begin{enumerate}
\item {The group $C^*$-algebras of connected, simply connected nilpotent Lie groups, including the real Heisenberg group $H_3(\R)$. These give examples with infinite-dimensional irreducible representations and non-Hausdorff primitive ideal spaces; see \cite{beltita-beltita-ludwig-fourier}.}
\item {All separable continuous-trace algebras. If the Dixmier--Douady class vanishes, such an algebra is strongly Morita equivalent to $C_0(X)$, where $X$ is its spectrum; see {\cite[Definition~IV.1.7.10 and Theorem~IV.1.7.15]{blackadar-operator-algebras}}. This case is already covered by Example~\ref{ex:commutative-nqs} and Morita invariance of membership in $\mathcal U_{\mathrm{nq}}$. Nonzero Dixmier--Douady classes give further examples: a locally trivial bundle of compact operators over $S^3$ with nonzero class in $H^3(S^3;\Z)$ has a separable continuous-trace section algebra that is not Morita equivalent to $C(S^3)$ {\cite[IV.1.7.12 and Theorem~IV.1.7.15]{blackadar-operator-algebras}}.}
\end{enumerate}
{All these algebras therefore belong to $\mathcal U_{\mathrm{nq}}$.}
\hB
\end{ex}

We next observe that  ${\mathcal U_{\mathrm{nq}}}$ is stable under taking maximal tensor products.
Let $A$ and $D$ be in $\nCalg_{\sepa}$. 
\begin{prop}\label{prop:U-nq-max-tensor}
{If $A,D\in\mathcal U_{\mathrm{nq}}$, then $A\otimes_{\max}D\in\mathcal U_{\mathrm{nq}}$.}
{If $A$ is AF and $D\in\mathcal U_{\mathrm{nqs}}$, then $A\otimes D\in\mathcal U_{\mathrm{nqs}}$.}
\end{prop}

\begin{proof}
{For fixed $A$, the functor $D\mapsto j(A\otimes_{\max}D)$ is additively continuous: the maximal tensor product preserves homotopies, finite products, and countable filtered colimits in each variable, and tensoring a corner inclusion with $A$ gives a corner inclusion under the canonical associativity isomorphism.} 
The universal property of $j$ therefore extends this functor to an
additive endofunctor of $\mathcal U$ preserving countable filtered
colimits.

If $A$ and $D$ are nuclear and quasidiagonal, then
$A\otimes_{\max}D$ is nuclear and quasidiagonal. Indeed, nuclearity
identifies the maximal and minimal tensor products, and the tensor
products of the finite-rank projections in the tensor product of  faithful quasidiagonal
representations give a quasidiagonal representation of the minimal
tensor product. Thus for nuclear quasidiagonal $A$, the endofunctor
just constructed preserves $\mathcal U_{\mathrm{nq}}$. Now fix $D$
 in $\mathcal U_{\mathrm{nq}}$. By the preceding argument
and symmetry, the endofunctor induced by $A\mapsto j(A\otimes_{\max}D)$
carries the generators of $\mathcal U_{\mathrm{nq}}$ into
$\mathcal U_{\mathrm{nq}}$. It preserves finite products, retracts,
and countable filtered colimits, so it preserves
$\mathcal U_{\mathrm{nq}}$.

{For finite-dimensional $A$, tensoring with $A$ preserves nuclearity, quasidiagonality, and semiprojectivity: we use the decomposition of $A$ into matrix algebras and preservation of these properties under matrix amplification and finite direct sums. The induced endofunctor of $\mathcal U$ therefore preserves the generators of $\mathcal U_{\mathrm{nqs}}$ and their closure under finite products, retracts, and countable filtered colimits.}

{If $A$ is an {AF algebra}, we write $A\cong\colim_n A_n$ with finite-dimensional $A_n$. For $D\in\mathcal U_{\mathrm{nqs}}$, the finite-dimensional case and closure under countable filtered colimits give}
\[
 {A\otimes D\cong\colim_n(A_n\otimes D)\in\mathcal U_{\mathrm{nqs}}.}
\]
\end{proof}

\begin{rem}\label{rem:U-nqs-tensor-products}
{The proof of Proposition~\ref{prop:U-nq-max-tensor} does not establish closure of $\mathcal U_{\mathrm{nqs}}$ under maximal tensor products. For this, it would suffice to show that the tensor product of two nuclear quasidiagonal semiprojective algebras belongs to $\mathcal U_{\mathrm{nqs}}$. Such a tensor product need not itself be semiprojective: $C(S^1)\otimes C(S^1)\cong C(\mathbb T^2)$ is not semiprojective by \cite[Theorem~1.2]{sorensen-thiel}.
Currently we do not know whether 
 $C(\mathbb T^2)$ belongs to $\mathcal U_{\mathrm{nqs}}$.}
\hB
\end{rem}

\begin{ex}\label{ex:U-nq-tensor-products}
\begin{enumerate}
\item {$C^*_{\max}(\prod_{i=1}^{n}F_{r_i})\in\mathcal U_{\mathrm{nq}}$ for positive integers $n$ and $r_1,\ldots,r_n$ in $\nat_{\geq1}\cup\{\omega\}$. Indeed, the full group $C^*$-algebra of a finite direct product is the maximal tensor product of the full group $C^*$-algebras of the factors. We apply Example~\ref{ex:U-nq-free-products} and Proposition~\ref{prop:U-nq-max-tensor}.}
\item {$M_{2^\infty}\otimes C^*_{\max}(F_2)\in\mathcal U_{\mathrm{nqs}}$. The CAR algebra is AF by Example~\ref{UHF}, so this follows from Example~\ref{ex:U-nq-free-products} and the second assertion of Proposition~\ref{prop:U-nq-max-tensor}.}
\item {$C^*_{\max}(H\times F_r)\in\mathcal U_{\mathrm{nqs}}$ for countable locally finite groups $H$ and integers $r\geq2$. Indeed, $C^*(H)$ is AF by Example~\ref{UHF}, and
$C^*_{\max}(H\times F_r)\cong C^*(H)\otimes C^*_{\max}(F_r)$,
so Proposition~\ref{prop:U-nq-max-tensor} applies. 
}
\end{enumerate}

{In 1., the algebra is nonnuclear unless every $r_i=1$ \cite{lance-nuclear}. The algebras in 2. and 3. are nonnuclear: $C^*_{\max}(F_2)$ is the range of a conditional expectation in 2., and $C^*_{\max}(F_r)$ is a quotient in 3.}

{The algebras in 2. and 3. are quasidiagonal. In 1., this follows when at most one factor is nonabelian: the other factors are nuclear, so the maximal tensor product agrees with the minimal tensor product, which preserves quasidiagonality.}

{In 1., semiprojectivity fails for $n\geq2$, since $C(\mathbb T^2)$ is a retract, and for $n=1$, $r_1=\omega$, by Example~\ref{ex:U-nq-free-products}. It also fails in 2.: the algebra is an inductive limit of ${M}_{2^k}(\C)\otimes C^*_{\max}(F_2)$, whose $K$-theory is finitely generated, so semiprojectivity would force its own $K$-theory to be finitely generated by {\cite[Proposition~II.8.3.18]{blackadar-operator-algebras}}, whereas its $K_0$-group is $\Z[1/2]$. In 3., infinite $H$ prevents semiprojectivity because the infinite-dimensional AF algebra $C^*(H)$ is a retract.}
\hB
\end{ex}

{Let $G$ be a countable group acting continuously on a compact metrizable space $X$. For a clopen subset $V\subseteq X$, we define its setwise stabilizer by}
\begin{equation}\label{eq:clopen-setwise-stabilizer}
 {G_V:=\{g\in G\mid gV=V\}.}
\end{equation}
\begin{prop}\label{prop:U-profinite-crossed-products}
{Assume that the clopen subsets of $X$ with finite $G$-orbits form a basis of its topology, and that $C^*_{\max}(G_V)\in\mathcal U_{\mathrm{nqs}}$ for every nonempty such subset $V$. Then $C(X)\rtimes_{\max}G\in\mathcal U_{\mathrm{nqs}}$.}
\end{prop}
\begin{proof}
{We choose a countable basis of clopen subsets with finite $G$-orbits. The translates of the first $n$ basis elements generate a finite Boolean algebra of clopen subsets. Its nonempty atoms form a finite $G$-invariant partition $X_n$ of $X$. These partitions refine one another and give $X\cong\lim_n X_n$. Choosing one atom $V_{n,j}$ from each orbit, we obtain}
\begin{equation}\label{eq:profinite-action-orbits}
 {X_n\cong\bigsqcup_{j=1}^{r_n}G/G_{V_{n,j}}.}
\end{equation}
{The full crossed product preserves inductive limits. Applying the finite-orbit case of Green's imprimitivity theorem \cite[Theorem~17]{green-imprimitivity} to \eqref{eq:profinite-action-orbits} and choosing coset representatives gives}
\begin{equation}\label{eq:profinite-action-crossed-product}
 {C(X)\rtimes_{\max}G
 \cong\colim_n\bigl(C(X_n)\rtimes_{\max}G\bigr)
 \cong\colim_n\bigoplus_{j=1}^{r_n}
 {M}_{[G:G_{V_{n,j}}]}(\C)\otimes_{\max}C^*_{\max}(G_{V_{n,j}}).}
\end{equation}
{Each atom $V_{n,j}$ is clopen with finite $G$-orbit. The claim therefore follows from the assumption, Proposition~\ref{prop:U-nq-max-tensor}, and closure under finite direct sums and countable filtered colimits.}
\end{proof}

\begin{ex}\label{ex:U-nq-profinite-crossed-product}
{We apply Proposition~\ref{prop:U-profinite-crossed-products} with $G=F_2$ and $X=\widehat{F_2}$, the profinite completion, with its translation action. The cosets of open normal subgroups form a clopen basis with finite $F_2$-orbits. The setwise stabilizer of any clopen subset with finite orbit has finite index in $F_2$, and hence is a finitely generated free group. Its full group algebra belongs to $\mathcal U_{\mathrm{nqs}}$ by Example~\ref{ex:U-nq-free-products}. We obtain $C(\widehat{F_2})\rtimes_{\max}F_2\in\mathcal U_{\mathrm{nqs}}$.}
\hB
\end{ex}

\begin{ex}\label{ex:U-nq-irrational-rotation}
For $\theta$ in $\R\setminus\mathbb Q$, let $\Z$ act on $\mathbb T$
by the rotation $z\mapsto e^{2\pi i\theta}z$. The irrational
rotation algebra $C(\mathbb T)\rtimes_\theta\Z$ is nuclear because
$C(\mathbb T)$ is nuclear and $\Z$ is amenable. It is quasidiagonal
by the embedding into an AF algebra of Pimsner and
Voiculescu~\cite{pimsner-voiculescu-rotation}. {Thus $C(\mathbb T)\rtimes_\theta\Z\in\mathcal U_{\mathrm{nq}}$}.

{More precisely, by Elliott--Evans \cite[Theorem~4]{elliott-evans-rotation}, the irrational rotation algebra is a sequential inductive limit of finite direct sums of matrix algebras over $C(S^1)$. It therefore belongs to $\mathcal U_{\mathrm{nqs}}$ by Example~\ref{ex:commutative-nqs}, Proposition~\ref{prop:U-nq-max-tensor}, and closure under finite direct sums and countable filtered colimits.}
\hB
\end{ex}

\begin{ex}\label{ex:rho-free-product-square}
{The nuclear quasidiagonal algebra
$C^{*}(\mathbb Z^2)\cong C(\mathbb T^2)$ belongs to
$\mathcal U_{\mathrm{nq}}$, but is not $\rho$-admissible.
Thus membership in $\mathcal U_{\mathrm{nq}}$ does not imply
$\rho$-admissibility. 
}

{Let $x,y$ be the classes in $K_1(C(\mathbb T^2))$ of the
two coordinate unitaries. Their product, identified with a degree-zero
class by Bott periodicity, is the reduced Bott class $b$.
The index pairing gives a homomorphism}
 \begin{equation}\label{eq:rho-torus-bott-pairing}
 \langle b,-\rangle:E_0(C(\mathbb T^2),\C)\longrightarrow\mathbb Z
\end{equation}which takes the value $1$ on the Dirac class.
For every homomorphism
$\varphi:C(\mathbb T^2)\to\bbK^u\otimes\bbK$, we have
$K_0(\varphi)(b)=0$.
{Indeed, the intersection of its commutative image with the compact ideal $\bbK\otimes\bbK$ is a commutative algebra of compact operators, while its image in the scalar quotient $\bbK$ is finite-dimensional, since it is supported on the finite-rank projection induced by $\varphi(1)$. The spectrum of the commutative image is therefore countable and compact, hence zero-dimensional. Thus the image is AF and its $K_1$ vanishes.}
The quotient map from $C(\mathbb T^2)$ to this commutative image
induces a ring homomorphism on $K$-theory. It sends $x$ and $y$ to zero,
and hence also their product $b$.

{Since the scalar projection has a section, the fibre sequence defining $\cT$ identifies}
\begin{equation}\label{eq:rho-torus-component-kernel}
 {\pi_0\cT_{C(\mathbb T^2)}(\C)
 \cong\ker\!\left(\pi_0{\cM}_{C(\mathbb T^2)}(\bbK^u)^{\group}
 \longrightarrow\pi_0{\cM}_{C(\mathbb T^2)}(\C)^{\group}\right).}
\end{equation}
{Under \eqref{eq:rho-torus-component-kernel}, every element of $\pi_0\cT_{C(\mathbb T^2)}(\C)$ is a formal difference of classes of homomorphisms $\varphi:C(\mathbb T^2)\to\bbK^u\otimes\bbK$.}
Consequently, the image  $ \pi_0\tau_{C(\mathbb T^2),\C}$
is annihilated by \eqref{eq:rho-torus-bott-pairing}.
In contrast, Proposition~\ref{prop:intro-splitting} implies that
$\pi_0\alpha_{C(\mathbb T^2),\C}$ is surjective, since
$C(\mathbb T^2)$ is nuclear and quasidiagonal.
Using \eqref{eq:tau-alpha-rho}, we conclude that
$\pi_0\rho_{C(\mathbb T^2),\C}$ is not surjective.

{
For integers $k,m\geq2$, Example~\ref{ex:U-nq-free-products} and
Proposition~\ref{prop:U-nq-max-tensor} give}
{\begin{equation}\label{eq:rho-free-group-products-nq}
 C^{*}_{\max}(F_k\times F_m)
 \cong C^{*}_{\max}(F_k)\otimes_{\max}C^{*}_{\max}(F_m)
 \in\mathcal U_{\mathrm{nq}}.
\end{equation}}{Choosing one generator in each free factor gives split group
homomorphisms \[
{\mathbb Z^2\to F_k\times F_m\to\mathbb Z^2.}
\]
Thus $C^{*}(\mathbb Z^2)$ is a retract of
$C^{*}_{\max}(F_k\times F_m)$. Since $\rho$-admissibility passes to
retracts by naturality of \eqref{eq:T-to-R}, none of these algebras
is $\rho$-admissible. 
}

{This also shows that $\rho$-admissibility does not pass to
$C^{*}$-subalgebras. Let $A$ be the universal unital $C^{*}$-algebra
generated by a contraction. By
\cite[Theorem~6.10, the proof of Theorem~6.11, and Remark~6.14]{courtney-sherman-contraction},
we have an embedding}
{\begin{equation}\label{eq:rho-universal-contraction-embedding}
 C^{*}_{\max}(F_2\times F_2)
 \cong C^{*}(F_2)\otimes_{\max}C^{*}(F_2)
 \hookrightarrow A\otimes_{\max}A,
\end{equation}}{and $A\otimes_{\max}A$ is quasidiagonal.
Scaling both universal contractions to zero gives a homotopy equivalence
$A\otimes_{\max}A\simeq\C$.
The constant shape system for the semiprojective algebra $\C$, together
with \eqref{eq:R-shape-formula}, shows that $\C$ is $\rho$-admissible.
Homotopy invariance therefore implies that $A\otimes_{\max}A$ is
$\rho$-admissible.} \hB
\end{ex}

\begin{kor}\label{cor:U-nqs-remainder}
 We have  
$
 \mathcal U_{\mathrm{nqs}}\subseteq\mathcal U_\alpha
$.
\end{kor}
\begin{proof}
{For a separable nuclear quasidiagonal semiprojective
algebra $D$, Corollary~\ref{cor:semiprojective-rho-admissible}
gives $\rho$-admissibility, and Theorem~\ref{thm:intro-exactness}
shows that ${\cR}_D$ is homological.
Proposition~\ref{prop:rho-homological-alpha} therefore implies that
$D$ is $\alpha$-admissible.}

{By Definition~\ref{def:alpha-admissible} and the closure properties established there, $\mathcal U_\alpha$ contains the generators of $\mathcal U_{\mathrm{nqs}}$ and is closed under finite products, retracts, and countable filtered colimits. Consequently, $\mathcal U_{\mathrm{nqs}}\subseteq\mathcal U_\alpha$.}
\end{proof}

\begin{prop}\label{ex:U-nq-oinfinity}
{The Cuntz algebra $\mathcal O_\infty$ does not belong to $\mathcal U_{\mathrm{nq}}$. In particular, $\mathcal U_{\mathrm{nq}}$ is properly contained in $\mathcal U$.}
\end{prop}
\begin{proof}
{The algebra $\mathcal O_\infty$ is semiprojective, so its constant shape system gives
${\cR}_{\mathcal O_\infty}\simeq\cT_{\mathcal O_\infty}$.
Example~\ref{ex:intro-nonexact} shows that this functor is
not exact, and hence not homological. By
Theorem~\ref{thm:intro-exactness}, we therefore have
$\mathcal O_\infty\notin \mathcal U_{\mathrm{nq}}$.}
\end{proof}

\begin{rem}\label{rem:remainder-candidates}
{The preceding examples suggest the following algebras and families in $\mathcal U_{\mathrm{nq}}$ for which membership in $\mathcal U_{\mathrm{nqs}}$ has not been established here in general:}
\begin{enumerate}
\item {$C(\mathbb T^d)\cong C^*_{\max}(\Z^d)$ for integers $d\geq2$; see Examples~\ref{ex:U-nq-tensor-products} and~\ref{ex:rho-free-product-square};}
\item {$C^*_{\max}(\prod_{i=1}^n F_{r_i})$ for $n\geq2$, with at least one $r_i\geq2$, allowing $r_i=\omega$; see Example~\ref{ex:U-nq-tensor-products}. This includes $C^*_{\max}(F_2\times F_2)$.}
\item {General separable commutative algebras $C_0(X)$, beyond the cases covered by Corollary~\ref{cor:AF-U-nqs} and Example~\ref{ex:commutative-nqs}. These extend the torus examples in 1.;}
\item {Group $C^*$-algebras of connected, simply connected nilpotent Lie groups, such as $C^*(H_3(\R))$; see Example~\ref{ex:liminal-algebras};}
\item {Separable continuous-trace algebras, in particular those with nonzero Dixmier--Douady class over $S^3$; see Example~\ref{ex:liminal-algebras}. For zero Dixmier--Douady class, Morita equivalence reduces the question to the commutative case in 3.}
\end{enumerate}
{For all these algebras, ${\cR}_D$ is homological. The members not already known to belong to $\mathcal U_{\mathrm{nqs}}$ provide candidates for investigating a nonzero remainder ${\widehat{\cR}}_D$ in subsequent work, since Corollary~\ref{cor:U-nqs-remainder} does not settle its vanishing in these cases. Nonmembership in $\mathcal U_{\mathrm{nqs}}$ would not by itself imply such nonvanishing.}

\hB
\end{rem}

\begin{rem}\label{rem:U-nq-nqs}{By Definition~\ref{def:U-nq} and
Corollary~\ref{cor:U-nqs-remainder}, we have}
\begin{equation}\label{eq:admissible-known-containments}
 \mathcal U_{\mathrm{nqs}}\subseteq
 \mathcal U_{\mathrm{nq}}\cap\mathcal U_{\alpha}.
\end{equation}{We do not know whether
$\mathcal U_{\mathrm{nqs}}=\mathcal U_{\mathrm{nq}}$ or
$\mathcal U_{\mathrm{nq}}=\mathcal U_{\alpha}$.
For the latter equality, we have no counterexample to either inclusion.
The inclusion $\mathcal U_{\mathrm{nq}}\subseteq\mathcal U_{\alpha}$
asks whether the remainder vanishes throughout $\mathcal U_{\mathrm{nq}}$.}

{{The equality $\mathcal U_{\mathrm{nqs}}=\mathcal U_{\mathrm{nq}}$ would imply this vanishing, since $\mathcal U_{\mathrm{nqs}}\subseteq\mathcal U_\alpha$ by \eqref{eq:admissible-known-containments}. The equality $\mathcal U_{\mathrm{nqs}}=\mathcal U_{\mathrm{nq}}$ would follow if every}
separable nuclear quasidiagonal algebra were a sequential inductive
limit of nuclear quasidiagonal semiprojective algebras.
This is a version of Blackadar's question about strong shape systems,
that is, presentations as sequential inductive limits of separable
semiprojective algebras; see
\cite[Problem~LXXXVIII]{schafhauser-tikuisis-white-problems}
for the general and nuclear versions.
{Algebras admitting strong shape systems are closed under shape domination in the sense of \cite[Section~2 and Theorem~4.4]{thiel-semiprojective-limits}. This closure result concerns arbitrary semiprojective stages and does not assert that the stages can be chosen nuclear and quasidiagonal.}}

{More generally, we have no example for which ${\cR}_D$ is
homological and ${\widehat{\cR}}_D\not\simeq0$. We therefore ask whether}
{\begin{equation}\label{eq:homological-implies-alpha-question}
 {\cR}_D\text{ homological}\quad\Longrightarrow\quad
 D\in\mathcal U_{\alpha}.
\end{equation}}{An affirmative answer would imply $\mathcal U_{\mathrm{nq}}\subseteq\mathcal U_\alpha$ by Theorem~\ref{thm:intro-exactness}. Proposition~\ref{prop:rho-homological-alpha} proves \eqref{eq:homological-implies-alpha-question} under the additional assumption that $D$ is $\rho$-admissible. Example~\ref{ex:rho-free-product-square} shows that this assumption is not automatic for $D\in\mathcal U_{\mathrm{nq}}$: neither $C(\mathbb T^2)$ nor $C^*_{\max}(F_k\times F_m)$ for $k,m\geq2$ is $\rho$-admissible. 
}
\hB
\end{rem}

\begin{prob}
In view of Remark \ref{ku-str} the functors ${\cM}^{\group}$, $\cT$ and ${\cR}$ have refinements to $ku$-module valued functors and the comparison maps $\eta$, $\tau$ and $\alpha$ are $ku$-linear. Furthermore,  the splitting in Theorem \ref{thm:overview-one-splitting} is a splitting of $ku$-modules. 
{If ${\cR}_D$ is homological, then the Bott element
$\beta\in\pi_2ku$ induces an equivalence}
{\begin{equation}\label{eq:intro-bott-periodicity}
 {\cR}_D\xrightarrow[\simeq]{\beta}\Omega^2{\cR}_D.
\end{equation}}

{We ask under which conditions on $D$ the functor
${\cR}_D[\beta^{-1}]$ is homological. This is not always the case:
$\cT_{\mathcal O_\infty}\simeq {\cR}_{\mathcal O_\infty}$, and the failure
of exactness in Example~\ref{ex:intro-nonexact} persists after Bott
inversion. Indeed, evaluation on the unit induces a $ku$-linear map}
{\begin{equation}\label{eq:intro-bott-detection}
 \cT_{\mathcal O_\infty}(\mathcal E_2)
 \longrightarrow\Omega^\infty K^{\topp}(\mathcal E_2)
\end{equation}}
{which sends the class of the homomorphism
$\mathcal O_\infty\to\mathcal E_2$ to
$[1_{\mathcal E_2}]\ne0$ in $K_0(\mathcal E_2)$.
Since its target is already Bott periodic, \eqref{eq:intro-bott-detection}
factors through Bott inversion. Thus this class survives, whereas
the values at $\bbK$ and $\mathcal O_2$ remain zero.}
\hB
 \end{prob}

\subsection{\texorpdfstring{{The universal class and the Yoneda retraction}}{The universal class and the Yoneda retraction}}
\label{sec:general-universal-class}

{Let $D$ be in $\nCalg_{\sepa}$.}
\begin{kor}\label{cor:R-E-factorization}
{{If ${\cR}_D$ is homological, then there is a unique}
left-exact functor ${\bar{\cR}}_D$ completing the commutative diagram
\[
 \xymatrix@C=4em@R=3em{
 \nCalg_{\sepa}\ar[r]^{e_{\sepa}}\ar[dr]_{{\cR}_D}
     &\EE_{\sepa}\ar@{..>}[d]^{{\bar{\cR}}_D}\\
     &\CGrp(\Ani).
 }
\]}
\end{kor}
\begin{proof}
{The factorization follows from the universal property of $e_{\sepa}$ in Definition~\ref{def:universal-homological-functor}.}
\end{proof}

{The following is an application of the Yoneda lemma to the canonical point in ${\cR}_{D}(D)$.}
\begin{prop}\label{prop:intro-splitting}
{{If ${\cR}_D$ is homological, then there is a natural transformation}
\[
 \beta_D:\Omega^\infty E(D,-)\longrightarrow {\cR}_D
 \quad\text{in }\Fun(\nCalg_{\sepa},\CGrp(\Ani))
\]
such that $\alpha_D\beta_D\simeq\id$. Consequently,
\[
 {\cR}_D\simeq \Omega^\infty E(D,-)\oplus{\widehat{\cR}}_D,
 \qquad{\widehat{\cR}}_D:=\Fib_0(\alpha_D).
\]}
\end{prop}
{By Theorem~\ref{thm:intro-exactness}, $D\in\mathcal U_{\mathrm{nq}}$ implies that ${\cR}_D$ is homological, so Corollary~\ref{cor:R-E-factorization} and Proposition~\ref{prop:intro-splitting} apply in particular in this case. 
}
%

 \begin{proof}[{Proof of Proposition~\ref{prop:intro-splitting}}]
{We use the factorization ${\bar{\cR}}_D$ from Corollary~\ref{cor:R-E-factorization}.} {Choose a shape system $(D_n)$ for $D$.} Choose a rank-one projection $p$ in $\bbK$. The structure maps
$d_n:D_n\to D$ of the shape system give a cone of  homomorphisms
 
\[
 \xymatrix@C=4em@R=2em{
 \vdots\ar[d]&\\
 D_n\ar[dd]^{\alpha_{n,n+1}}\ar[dr]^{\gamma_n}&\\
 &(D\otimes\bbK)^u\otimes\bbK\\
 D_{n+1}\ar[d]\ar[ur]_{\gamma_{n+1}}&\\
 \vdots&
 }
\]
where
\[
 \gamma_n(x):=\iota_{D\otimes\bbK}(d_n(x)\otimes p)\otimes p
 \qquad(\forall x\in D_n).
\]
 Their composites with  $(D\otimes\bbK)^{u}\otimes\bbK\to \C\otimes \bbK$ vanish so that we get a point  
    \begin{equation}\label{grewfrewgg}
 u:*\longrightarrow\lim_n\cT_{D_n}(D)
       \simeq {\cR}_D(D)\simeq{\bar{\cR}}_D(e_{\sepa}(D)).
 \end{equation} 
 The  Yoneda lemma applied to $u$ gives a natural transformation
\[
 \bar\beta_D:\Hom_{\EE_{\sepa}}(e_{\sepa}(D),-)
       \longrightarrow{\bar{\cR}}_D :\EE_{\sepa}\to \CGrp(\Ani).
\]
 Restriction along $e_{\sepa}$ gives the  asserted  natural transformation
\[
\beta_D:{\Omega^\infty E(D,-)}\to {\cR}_{D}:\nCalg_{\sepa}\to \CGrp(\Ani)\ .
\]%
Unfolding the definition of $\alpha_{D}$ in \eqref{eq:intro-alpha} we have an equivalence
\[
\alpha_{D,D}u \simeq (e_{\sepa}(d_n))_n
\]%
under the equivalence
\[
\Omega^{\infty}E(D,D)\simeq \lim_{n} \Hom_{\EE_{\sepa}}(e_{\sepa}(D_{n}),e_{\sepa}(D))\ .
\]
Since $(e_{\sepa}(d_n))_n $ also presents $\id_{e_{\sepa}(D)}$ under this equivalence we get
\begin{equation}\label{eq:general-universal-class}
 \alpha_{D,D}u\simeq\id_{e_{\sepa}(D)}.
\end{equation}
 Using the Yoneda lemma, $\bar\alpha_D\bar\beta_D\id_{e_{\sepa}(D)}\simeq \alpha_{D,D}u\simeq \id_{e_{\sepa}(D)}$
 and the universal property of $e_{\sepa}$, we conclude that $\alpha_{D}\beta_{D}\simeq \id$.
{The section $\beta_D$ splits the fibre sequence ${\widehat{\cR}}_D\to {\cR}_D\xrightarrow{\alpha_D}\Omega^\infty E(D,-)$, giving the asserted decomposition.}
%
%
\end{proof}

\section{Exactness properties of ${\cR}_{D}$}
\label{sec:resolved-exactness}

\subsection{Statement of the exactness theorem}

In this section we analyse the exactness properties {of the functor} 
${\cR}_D$ of Definition~\ref{kohprtgerthertgertgertge}. {For a shape system $(D_n)$ of $D$, formula~\eqref{eq:R-shape-formula} gives ${\cR}_D\simeq\lim_n\cT_{D_n}$, where $\cT_{D_n}$ is obtained from the representation monoid functor ${\cM}_{D_n}$ of Definition~\ref{def:mapping-monoid} by group completion and stabilized reduction.} The functor ${\cM}_{D}$ sends Schochet fibrant cartesian squares {(see Definition~\ref{def:functor-properties}(\ref{hjjtjrtzh}) and \cite[Section~3]{orbita})} to cartesian squares for every $D$ in $\nCalg_{\sepa}$. {Group completion does not preserve pullbacks in general, so Schochet exactness need not pass from ${\cM}_D$ to ${\cM}_D^{\group}$.} Furthermore, for general $D$ the functor ${\cM}_{D}$ does not send exact sequences to fibre sequences. {The main goal of this section is to show that, under suitable assumptions on $D$, the limit in \eqref{eq:R-shape-formula} is exact and Schochet exact.}

Let $D$ be in $\nCalg_{\sepa}$. {Consider a cartesian square in $\nCalg_{\sepa}$}
\begin{equation}\label{rtgerthrthr}
 {\xymatrix@C=3em@R=2.5em{
 B_{00}\ar[r]\ar[d]&B_{01}\ar[d]^q\\
 B_{10}\ar[r]&B_{11}.
 }}
\end{equation}
\begin{theorem}\label{thm:general-resolved-exactness}
{If $D$ is nuclear and quasidiagonal and $q$ is surjective, then ${\cR}_D$ sends \eqref{rtgerthrthr} to a cartesian square in $\CGrp(\Ani)$.}
\end{theorem}
{We prove Theorem~\ref{thm:general-resolved-exactness} in Section~\ref{sec:resolved-cartesian-proof}.}

\begin{kor}\label{kophertgrtgerg}
If $D$ is nuclear and quasidiagonal, then ${\cR}_{D}$ is exact and Schochet exact.
\end{kor}
\begin{proof}
Exactness and Schochet exactness {require} that ${\cR}_{D}$  preserves {the cartesianness} of   squares \eqref{rtgerthrthr} 
 with $B_{10}=0$ and $q$ surjective,  or with $q$ a Schochet fibration.
 Since in both cases $q$ is surjective, and since ${\cR}_{D}(0)\simeq 0${, the assertion follows} from Theorem \ref{thm:general-resolved-exactness}. 
\end{proof}

\begin{proof}[Proof of Theorem~\ref{thm:intro-exactness} assuming Theorem \ref{thm:general-resolved-exactness}]
{By Definition~\ref{kohprtgerthertgertgertge}, the contravariant
functor $D\mapsto {\cR}_D$ factors through $j$ as
\begin{equation}\label{eq:resolved-source-factorization}
 (\nCalg_{\sepa})^{\op}\xrightarrow{j^{\op}}
 \mathcal U^{\op}\longrightarrow
 \Fun(\nCalg_{\sepa},\CGrp(\Ani)),
 \qquad V\longmapsto\Hom_{\mathcal U}(j(V),j^{\sfred}(-)).
\end{equation}
{Here $V$ is a separable $C^*$-algebra, and the assignment in \eqref{eq:resolved-source-factorization} describes the entire composite. The intermediate functor on $\mathcal U^{\op}$ sends countable filtered colimits in}
$\mathcal U$ to limits and finite products to finite products,
since finite products in $\mathcal U$ are also coproducts.
It also preserves retracts.}

{The subcategory $\Fun^{\mathrm{hom}}(\nCalg_{\sepa},\CGrp(\Ani))$ of homological functors  is
closed under limits and retracts. Consequently, the objects
of $\mathcal U$ whose associated functors are homological form
a full subcategory closed under equivalences, retracts, finite
products, and countable filtered colimits. For nuclear
quasidiagonal $D$, the functor ${\cR}_D$ is reduced, homotopy
invariant by construction and inverts corner inclusions, and
Corollary~\ref{kophertgrtgerg} supplies exactness and Schochet
exactness. Thus this subcategory contains the generators of
$\mathcal U_{\mathrm{nq}}$, and hence all of
$\mathcal U_{\mathrm{nq}}$. The assumption $D\in\mathcal U_{\mathrm{nq}}$ now
implies that ${\cR}_D$ is homological.}
\end{proof}

\subsection{Scalar complements and relative lifting}
\label{sec:general-exactness}

{We use semiprojectivity and nuclearity to straighten and lift maps and homotopies. In the subsequent subsections, we apply these results to prove exactness of ${\cR}_D$.}

\begin{lem}[Straightening lemma]\label{lem:semiprojective-straightening}
{Let $u:P\to Q$ be a semiprojective morphism in $\nCalg_{\sepa}$.}
\begin{enumerate}
\item For every finite subset $F$ of $P$, straightening accuracy $\varepsilon$ in $(0,\infty)$, and norm bound $M$ in $(0,\infty)$, there are a finite subset $G$ of $Q$ and a tolerance  $\delta$
in $(0,\infty)$ such
that the following holds for every {$C$ in $\nCalg$}.
If a bounded linear adjoint-preserving map $f:Q\to C$ satisfies
\[
 \|f\|\leq M,\qquad
 \|f(xy)-f(x)f(y)\|<\delta\quad{(\forall x,y\in G)},
\]
then there is a homomorphism $h:P\to C$ with
\[
 \|h(a)-f(u(a))\|<\varepsilon\quad{(\forall a\in F)}.
\]
\item There are a finite subset $G$ of $Q$ and tolerance $\delta$ in
$(0,\infty)$ such
that the following holds for every {$C$ in $\nCalg$}.
If homomorphisms $h_0,h_1:Q\to C$ satisfy
\[
 \|h_0(x)-h_1(x)\|<\delta\quad{(\forall x\in G)},
\]
then there is a homomorphism $H:P\to C([0,1],C)$ with
\[
 \operatorname{ev}_0H=h_0u,\qquad
 \operatorname{ev}_1H=h_1u.
\]
\end{enumerate}
\end{lem}

\begin{rem}\label{rem:straightening-relative-parameters}
{Let $K$ be a compact Hausdorff space and let $L\subseteq K$ be a closed subspace. Both assertions of Lemma~\ref{lem:semiprojective-straightening} hold for families parametrized by $K$, with the same finite subsets and tolerances and with all estimates in the uniform norm. They also admit the following relative forms.}
{In part~1 of Lemma~\ref{lem:semiprojective-straightening}, let $f:Q\to C(K,C)$ satisfy the stated conditions. If $a\mapsto f(a)|_L$ is a homomorphism from $Q$ to $C(L,C)$, we can choose $h:P\to C(K,C)$ such that}
\begin{equation}\label{eq:straightening-relative-boundary}
 {h(a)|_L=f(u(a))|_L\qquad(\forall a\in P).}
\end{equation}
{In part~2, let $h_0,h_1:Q\to C(K,C)$ satisfy the stated closeness condition and agree after restriction to $L$. We can then choose $H:P\to C([0,1]\times K,C)$ with the stated endpoints and}
\begin{equation}\label{eq:straightening-relative-constant-homotopy}
 {H(a)|_{[0,1]\times L}
 =1_{C([0,1])}\otimes\bigl(h_0(u(a))|_L\bigr)
 \qquad(\forall a\in P).}
\end{equation}
\end{rem}

\par\noindent{\fbox{\begin{minipage}{\dimexpr\linewidth-2\fboxsep-2\fboxrule\relax}
\color{black}
\noindent\emph{Schematic form of (1).}
\[
\xymatrix@C=1.2em@R=3em{
 F\subseteq
 &P\ar[d]_u\ar@{.>}[r]^h
 &C\ar@{=}[d]
 &\|h(a)-f(u(a))\|<\varepsilon\quad{(\forall a\in F)}
 \\
 G\subseteq
 &Q\ar[r]^f
 &C
 &{\|f(xy)-f(x)f(y)\|<\delta\quad(\forall x,y\in G).}
}
\]{The map $f$ is bounded linear and adjoint-preserving, with $\|f\|\leq M$.}
The subset $G$ and tolerance $\delta$ are chosen from
$u,F,\varepsilon,M$, before choosing $C$ and $f$.
\end{minipage}}}\par
\medskip

\par\noindent{\fbox{\begin{minipage}{\dimexpr\linewidth-2\fboxsep-2\fboxrule\relax}
\color{black}
\noindent\emph{Schematic form of (2).}
\[
\xymatrix@C=1.2em@R=3em{
 &P\ar[d]_u\ar@{.>}[r]^-H
 &C([0,1],C)\ar[d]^{(\operatorname{ev}_0,\operatorname{ev}_1)}
 & 
 \\
 G\subseteq
 &Q\ar[r]_-{(h_0,h_1)}
 &C\oplus C
 &\|h_0(x)-h_1(x)\|<\delta\quad{(\forall x\in G)}.
}
\]
The dotted homomorphism $H$ is a
homotopy from $h_0u$ to $h_1u$. The subset $G$ and tolerance
$\delta$ depend only on $u$.
\end{minipage}}}\par
\medskip

\begin{proof}
For the first assertion, we fix $F,\varepsilon,M$ and suppose
that no pair $(G,\delta)$ works. We then choose an increasing sequence  of  finite subsets
$G_j\subseteq Q$ with dense union and numbers $\delta_j$ in $(0,\infty)$ tending to zero.
For each $j$ in $\nat$, there are {$C_j$ in $\nCalg_{\sepa}$} and a bounded
linear adjoint-preserving map $f_j:Q\to C_j$ with
\[
 \|f_j\|\leq M,\qquad
 \|f_j(xy)-f_j(x)f_j(y)\|<\delta_j
 \qquad{(\forall x,y\in G_j)},
\]
but no homomorphism $h:P\to C_j$ satisfies
\[
 \|h(a)-f_j(u(a))\|<\varepsilon
 \qquad{(\forall a\in F)}.
\]
The family of  maps $(f_j)_{j\in \nat}$
induces a homomorphism
\[
 Q\longrightarrow{\prod_j}C_j/{\bigoplus_j}C_j,
 \qquad x\longmapsto[{(f_j(x))_{j}}].
\]
We now apply semiprojectivity of $u$ to the increasing family of ideals
${(\bigoplus_{j\leq N}}C_j)_{N}$. 
{We obtain a  factorization indicated by the dotted arrow in the diagram
\[
 \xymatrix@C=4em@R=3em{
 P\ar[r]^u\ar@/^-2cm/@{-->}[dd]_{h_{i}}\ar@{..>}[d]_{\widetilde f}
   &Q\ar[d]\\
 ({\prod_j}C_j)/{\bigoplus_{j\leq N}}C_j\ar[r]\ar[d]
   &({\prod_j}C_j)/{\bigoplus_j}C_j\\ C_{i}&
 }.
\] 
 For $i\geq N+1$, projecting further to the components   gives
homomorphisms $h_i:P\to C_i$ satisfying
\[
 \lim_{i\to\infty}\|h_i(a)-f_i(u(a))\|=0
 \qquad(\forall a\in P).
\]
Since $F$ is finite, this contradicts the prescribed approximation
failure for all sufficiently large $j$.}

For the second assertion, suppose again that no pair $(G,\delta)$
works. We then obtain {$C_j$ in $\nCalg_{\sepa}$} and homomorphisms
$h_{0,j},h_{1,j}:Q\to C_j$ such that
\[
 \lim_{j\to \infty}\|h_{0,j}(a)-h_{1,j}(a)\|=0
 \qquad{(\forall a\in Q)},
\]
but $h_{0,j}u$ and $h_{1,j}u$ are not homotopic for any $j$.
The linear interpolations
\[
 f_j(a)(t):=(1-t)h_{0,j}(a)+t h_{1,j}(a)
\]
have the prescribed endpoints 
and their   multiplicativity defects vanish at $t=0,1$ and tend
to zero in norm as $j\to\infty$.

We  form the pullback
\[
 \xymatrix@C=4em@R=3em{
 \mathcal B\ar[r]\ar[d]
   &{\prod_j} C([0,1],C_j)
       \ar[d]^{(\operatorname{ev}_0,\operatorname{ev}_1)}\\
 Q\ar[r]_-{{(h_{0,j},h_{1,j})_{j}}}
   &{\prod_j}(C_j\oplus C_j).
 }
\]
Consider the closed ideals
\[
 J:={\bigoplus_j} C_0((0,1),C_j)\subseteq\mathcal B,\qquad
 J_N:={\bigoplus_{j\leq N}}C_0((0,1),C_j)\subseteq J,
 \qquad \overline{{\bigcup_{N\in \nat}}J_N}=J.
\]
Hence we have a homomorphism
\[
 \gamma:Q\longrightarrow\mathcal B/J,\qquad
 a\longmapsto[({(f_j(a))_{j}},a)].
\]

Semiprojectivity of $u$ supplies an $N$ and a lift $\tilde \gamma$ as in  the commutative square
\[
 \xymatrix@C=4em@R=3em{
 P \ar@/^-2cm/@{-->}[dd]_{H_{i}}\ar[r]^u\ar@{.>}[d]_{\widetilde\gamma}
   & Q\ar[d]^\gamma\ar@/^2cm/@{-->}[dd]_{(h_{i,0},h_{i,1})} \\
 \mathcal B/J_N\ar[r]\ar[d]&\mathcal B/J\ar[d]\\ C([0,1],C_i)\ar[r]^{(\ev_{0},\ev_{1})}&C_{i}\oplus C_{i}
 }.
\] 
For $i>N$, the $i$th path
coordinate  gives  the left lower vertical arrow and therefore $H_{i}$. Its existence 
  {contradicts} the choice of $h_{0,j},h_{1,j}$.

{For the relative assertion in part~1, repeat the product-quotient construction with the pullback retaining the homomorphisms $a\mapsto f_j(a)|_L$. The error ideal is}
\[
 {\bigoplus_j\{b\in C(K,C_j):b|_L=0\},}
\]
{with its finite partial sums as the increasing ideals. Since the multiplicativity defects vanish on $L$, semiprojectivity gives lifts satisfying \eqref{eq:straightening-relative-boundary}. For part~2, the linear interpolations are constant on $L$ because $h_0|_L=h_1|_L$. They are therefore homomorphisms on $(\{0,1\}\times K)\cup([0,1]\times L)$. The same construction retains the endpoints and this constant boundary value, giving \eqref{eq:straightening-relative-constant-homotopy}.}
\end{proof}
{Let $(D_n)_{n\in\nat}$ be a sequence in $\nCalg_{\sepa}$
and put $D:=\colim_nD_n$. Write
$\alpha_{nm}:D_n\to D_m$ and $d_n:D_n\to D$ for its structure maps.}

\begin{lem}[Approximate cpc extension lemma]\label{lem:general-cp}
{Assume that $D$ is nuclear. For every $\ell$ in $\nat$,
finite subset $F$ of $D_\ell$, and accuracy $\varepsilon$ in $(0,\infty)$,
there is $m$ in $\nat$ with $m\geq\ell$ such that for every
$C$ in ${\nCalg}$ and morphism $\theta:D_m\to C$,
there is a cpc map $\phi:D\to C$ with
\[
 \|\theta\alpha_{\ell m}(x)-\phi d_\ell(x)\|<\varepsilon
 \qquad(\forall x\in F).
\]
}
\end{lem}

\par\noindent{\fbox{\begin{minipage}{\dimexpr\linewidth-2\fboxsep-2\fboxrule\relax}
\color{black}
\[
\xymatrix@C=1.2em@R=3em{
 F\subseteq
 &D_\ell\ar[d]_{\alpha_{\ell m}}
       \ar@<0.6ex>[r]^{\theta\alpha_{\ell m}}
       \ar@{.>}@<-0.6ex>[r]_{\phi d_\ell}
 &C\ar@{=}[d]
 &{\begin{gathered}
    \|\theta\alpha_{\ell m}(x)-\phi d_\ell(x)\|<\varepsilon\\
    {(\forall x\in F)}
   \end{gathered}}
 \\
 &D_m\ar[d]_{d_m}\ar[r]^\theta
 &C\ar@{=}[d]
 &{\theta\text{ homomorphism} }
 \\
 &D\ar@{.>}[r]^\phi
 &C
 &\phi\text{ cpc}.
}
\]
The stage $m$ is chosen from $\ell,F,\varepsilon$, before choosing
the target $C$ and the {morphism} $\theta$. The dotted arrows
record the cpc approximation supplied by the lemma.
\end{minipage}}}\par
\medskip

\begin{proof}
{Fix $\ell$ in $\nat$, a finite subset $F$ of $D_\ell$, and $\varepsilon$ in $(0,\infty)$ and suppose that no $m$ in $\nat$ with $m\geq\ell$ works.
Choose an increasing sequence $(m_j)_{j\in\nat}$ in $\nat$ with $m_j\geq\ell$
and $m_j\to\infty$.
For each $j$ in $\nat$, there are $C_j$ in ${\nCalg}$ and a
morphism $\theta_j:D_{m_j}\to C_j$ such that every cpc map
$\phi:D\to C_j$ satisfies
\begin{equation}\label{goerktpogkerpogkertgertg}
 \max_{x\in F}
 \|\theta_j\alpha_{\ell m_j}(x)-\phi d_\ell(x)\|
 \geq\varepsilon.
\end{equation}
For each fixed $n$ in $\nat$, the maps $\theta_j\alpha_{nm_j}$ are defined
for all sufficiently large $j$ and induce a homomorphism
\[
 D_n\longrightarrow\prod_jC_j/\bigoplus_jC_j,\qquad
 x\longmapsto[(\theta_j\alpha_{nm_j}(x))_{j\in\nat}].
\]
The finitely many undefined coordinates may be set to zero.
These homomorphisms are compatible with the transition maps,
so the universal property of $D=\colim_nD_n$ gives a homomorphism
\[
 \psi:D\longrightarrow\prod_jC_j/\bigoplus_jC_j.
\]
{Since $D$ is separable and nuclear, the Choi--Effros lifting theorem \cite{choi-effros-lifting} gives a cpc lift of $\psi$}
\[
 \Phi:D\longrightarrow\prod_jC_j.
\]
Its coordinate maps $(\phi_j:D\to C_j)_{j\in\nat}$ are cpc and satisfy
\begin{equation*}\label{gwegerfwerferfw}
 {\lim_{j\to\infty}\max_{x\in F}
 \|\theta_j\alpha_{\ell m_j}(x)-\phi_jd_\ell(x)\|=0.}
\end{equation*}  This contradicts \eqref{goerktpogkerpogkertgertg}.}
\end{proof}

{Let $D$ be in $\nCalg_{\sepa}$ and choose a
shape system $(D_n)_{n\in\nat}$ for $D$. For morphisms,
$\simeq$ denotes homotopy.}
For $A$  in $\nCalg$,  $B$ in $\Calg_{\sepa}$ (unital!), and $\tau:A\to \bbK$ we call
\begin{equation}\label{fqqwedweqfqwe}
\tau_{B}:=1_{B}\otimes \tau:A\to B\otimes \bbK
\end{equation}  the scalar representation associated to $\tau$.

\begin{prop}[Scalar complement proposition]\label{prop:general-complement}
If $D$ is nuclear and quasidiagonal, then
for every $n$ in $\nat$ there is $m$ in $\nat$ with $m\geq n$
such that for every {$B$ in $\Calg_{\sepa}$} and
{morphism} $\theta:D_m\to B\otimes\bbK$  there exist
{morphisms} $\eta:D_n\to B\otimes\bbK$ and $\tau:D_n\to\bbK$ with
\[
 \theta\alpha_{nm}\oplus\eta\simeq\tau_B. 
\]
\end{prop}

\par\noindent{\fbox{\begin{minipage}{\dimexpr\linewidth-2\fboxsep-2\fboxrule\relax}
\color{black}
\[
\xymatrix@C=2em@R=3em{
 D_n\ar[d]_{\alpha_{nm}}
    \ar@{.>}@<0.7ex>[r]^{\theta\alpha_{nm}\oplus\eta}
    \ar@{.>}@<-0.7ex>[r]_{\tau_B}
    \ar@{.>}@/^2.5pc/[rr]^\tau
 &B\otimes\bbK\ar@{=}[d]
 &\bbK\ar[l]_-{1_B\otimes-}
 &\theta\alpha_{nm}\oplus\eta\simeq\tau_B
 \\
 D_m\ar[r]^\theta
 &B\otimes\bbK
}
\]
The upper horizontal maps are joined by the homotopy supplied by
the proposition and $\tau_{B}$ is scalar.  The stage $m$ is chosen from $n$, before choosing $B$
and $\theta$.
\end{minipage}}}\par
\medskip

\begin{proof}
Fix $n$ in $\nat$. We introduce the approximation accuracies as parameters;
the estimates in Step~4 will determine their choice and hence the
stage $m$.

\smallskip
\noindent\emph{1. Homotopy and straightening controls.}
By the second part of
Lemma~\ref{lem:semiprojective-straightening}, applied to the {semiprojective morphism}
$\alpha_{n,n+1}$, there are a finite self-adjoint subset $F$ of the
unit ball of $D_{n+1}$ and a tolerance $\delta$ in $(0,\infty)$ such that,
{for every $C$ in $\nCalg_{\sepa}$,} two {morphisms}
$D_{n+1}\to C$ which are $\delta$-close on $F$ become homotopic
after composition with $\alpha_{n,n+1}$.

We take $\kappa$ in $(0,\delta)$.  
Its precise choice  is deferred to Step~4. For each $\kappa$, apply the first part
of Lemma~\ref{lem:semiprojective-straightening} to the {semiprojective morphism} $\alpha_{n+1,n+2}$, with norm bound $1$ and
straightening accuracy $\kappa$ on $F$. This gives a finite
self-adjoint subset $G$ of the unit ball of $D_{n+2}$ and
a tolerance $\gamma$ in $(0,\infty)$
such that a contractive adjoint-preserving linear map $f:D_{n+2}\to C$ satisfying
\[
 \|f(xy)-f(x)f(y)\|<\gamma\qquad{(\forall x,y\in G)}
\]
admits a {morphism} $\widetilde f:D_{n+1}\to C$ with
\[
 \|\widetilde f(x)-f\alpha_{n+1,n+2}(x)\|<\kappa
 \qquad{(\forall x\in F)}.
\]
Enlarge $G$ to contain $\alpha_{n+1,n+2}(F)$ and put
$E:=G\cup\{xy:x,y\in G\}$. In particular, $E$ contains the products
$x^*x$ and $xx^*$ needed below. All these controls are independent
of the target algebra $C$, but $G$, $E$ and $\gamma$ depend on
$\kappa$.

\par\noindent{\fbox{\begin{minipage}{\dimexpr\linewidth-2\fboxsep-2\fboxrule\relax}
\color{black}
\begin{center}
\small
\[
\xymatrix@C=1.2em@R=3.5em{
 & D_n\ar[d]_{\alpha_{n,n+1}}
      \ar[r]^{h_i\alpha_{n,n+1}}
 & C
 & h_0\alpha_{n,n+1}\simeq h_1\alpha_{n,n+1}
 \\
 F\subseteq
 & D_{n+1}\ar[d]_{\alpha_{n+1,n+2}}
      \ar@<0.6ex>[r]^{h_i}
      \ar@{.>}@<-0.6ex>[r]_{\widetilde f}
 & C
 & {\begin{gathered}
     \|h_0(x)-h_1(x)\|<\delta\quad{(\forall x\in F)},\\
     \|\widetilde f(x)-f\alpha_{n+1,n+2}(x)\|<\kappa\quad{(\forall x\in F)}
   \end{gathered}}
 \\
 G\subseteq E\subseteq
 & D_{n+2}\ar[r]^f
 & C
 & \|f(xy)-f(x)f(y)\|<\gamma\quad{(\forall x,y\in G)}.
}
\]
\end{center}
{Here $i$ is in $\{0,1\}$.} The maps $h_0,h_1$ illustrate the homotopy
control; the contractive adjoint-preserving linear map $f$ and the dotted {morphism} $\widetilde f$
illustrate straightening. Once $\kappa$ is chosen in Step~4, it
determines $G$, $E$ and $\gamma$; the homotopy control $F$, $\delta$
is already fixed.
\end{minipage}}}\par
\medskip

\smallskip
\noindent\emph{2. Choice of $m$ and finite-corner approximation.}
For each $\kappa$ in $(0,\delta)$, we introduce the accuracy $\varepsilon$ in $(0,\infty)$ for the application of Lemma~\ref{lem:general-cp} at stage $n+2$ to the finite subset $E$ of $D_{n+2}$. It must be small enough that the linear adjoint-preserving contractive maps constructed
below satisfy the multiplicativity control $\gamma$, and that their
{norm differences used to obtain the homotopies in Step~4 satisfy the final homotopy control}. We defer its
choice, together with that of $\kappa$, to Step~4. The estimates
below explain why these requirements can be met.

Apply Lemma~\ref{lem:general-cp} at stage $n+2$, on $E$, with this
accuracy $\varepsilon$. For fixed values of the controls, this
fixes a stage $m$ in $\nat$ with $m\geq n+2$ which works uniformly for all target
algebras.

Having fixed $m$, let {$B$ be in $\Calg_{\sepa}$} and let
$\theta:D_m\to B\otimes\bbK$ be a {morphism}. Lemma~\ref{lem:general-cp} supplies
a cpc map $\phi:D\to B\otimes\bbK$ with approximation accuracy
$\varepsilon$ on $E$. Compressing by a sufficiently large
standard matrix corner gives a cpc   map $\phi_r:D\to M_r(B)$ for some $r$ in $\nat_{\geq1}$.
Choose the corner large enough that this compression changes $\phi$
on $d_{n+2}(E)$ by less than $\varepsilon$. Since $E$
contains $G$ and its products and $\theta\alpha_{n+2,m}$ is
multiplicative, the multiplicativity defects of $\phi_r d_{n+2}$
on $G$ tend to zero {as $\varepsilon\to0$}. Thus, for sufficiently small
$\varepsilon$, they satisfy the control $\gamma$.

Since cpc maps are contractive, linear and adjoint-preserving, Step~1 therefore gives a {morphism} $\rho:D_{n+1}\to M_r(B)$
approximating $\phi_r d_{n+1}$ on $F$ within $\kappa$.
For $x$ in $F$ the triangle inequality gives
\begin{equation}\label{eq:complement-first-approximation}
 \begin{split}
 \|\rho(x)-\theta\alpha_{n+1,m}(x)\|
 &< \kappa+2\varepsilon.
 \end{split}
\end{equation}
where matrix corners are viewed inside $B\otimes\bbK$.

\par\noindent{\fbox{\begin{minipage}{\dimexpr\linewidth-2\fboxsep-2\fboxrule\relax}
\color{black}
\begin{center}
\resizebox{\linewidth}{!}{$
\xymatrix@C=1.2em@R=3.5em{
 & D_n\ar[d]_{\alpha_{n,n+1}}
      \ar@<0.6ex>[r]^{\rho\alpha_{n,n+1}}
       \ar@<-0.6ex>[r]_{\theta\alpha_{n,m}}
 & B\otimes\bbK
 & {\begin{gathered}
     \rho\alpha_{n,n+1}\simeq\theta\alpha_{nm}\\
     \text{once }\kappa+2\varepsilon<\delta
   \end{gathered}}
 \\
 F\subseteq
 & D_{n+1}\ar[d]_{\alpha_{n+1,n+2}}
      \ar@<0.6ex>[r]^{\theta\alpha_{n+1,m}}
      \ar@{.>}@<-0.6ex>[r]_{\rho}
 & B\otimes\bbK
 & {\begin{gathered}
     \|\rho(x)-\phi_r d_{n+1}(x)\|<\kappa\quad{(\forall x\in F)},\\
     \|\rho(x)-\theta\alpha_{n+1,m}(x)\|
                  <\kappa+2\varepsilon\quad{(\forall x\in F)}
   \end{gathered}}
 \\
 G\subseteq E\subseteq
 & D_{n+2}\ar[d]_{\alpha_{n+2,m}}
      \ar@<0.6ex>[r]^{\theta\alpha_{n+2,m}}
      \ar@<-0.6ex>[r]_{\phi_r d_{n+2}}
 & B\otimes\bbK
 & {\begin{gathered}
     \|\phi d_{n+2}(x)-\theta\alpha_{n+2,m}(x)\|
                    <\varepsilon\quad{(\forall x\in E)},\\
     \|\phi_r d_{n+2}(x)-\theta\alpha_{n+2,m}(x)\|
                    <2\varepsilon\quad{(\forall x\in E)},\\
     \|\phi_r d_{n+2}(xy)-\phi_r d_{n+2}(x)\phi_r d_{n+2}(y)\|
                    <\gamma\quad{(\forall x,y\in G)}
   \end{gathered}}
 \\
 & D_m\ar[d]_{d_m}\ar[r]^\theta
 & B\otimes\bbK
 & \theta\text{ a {morphism}}
 \\
 d_{n+2}(E)\subseteq
 & D\ar@<0.6ex>[r]^{\phi}
      \ar@<-0.6ex>[r]_{\phi_r}
 & B\otimes\bbK
 & \|\phi_r(a)-\phi(a)\|<\varepsilon
                       \quad{(\forall a\in d_{n+2}(E))}.
}
$}
\end{center}
The maps $\phi_r d_{n+2}$ are cpc, and the dotted map $\rho$ is
obtained by straightening. The estimates for $\rho$ use
$\alpha_{n+1,n+2}(F)\subseteq E$. The top row records the homotopy
that follows when the accuracies are chosen in Step~4.
\end{minipage}}}\par
\medskip

\smallskip
\noindent\emph{3. Absorption and complementary compression.}
 Fix a faithful essential quasidiagonal representation
$\pi:D^{u}\to\mathcal L(H)$ of infinite multiplicity, with finite-rank
{projections $(p_k)_{k\in\nat}$ converging strongly to $1$
and satisfying $[p_k,\pi(a)]\to0$ as $k\to \infty$ for all $a$ in $D^u$.} Since $D^{u}$ is nuclear, $\pi_B:=1_B\otimes\pi$ is
{unitally absorbing by Kasparov's absorption theorem \cite[Theorem~6]{kasparov-absorption}; see \cite[Theorem~2.1(2)]{thomsen-absorbing} for the characterization using ucp maps, and also    \cite[Lemma~5.11]{willett}  for similar arguments as below.}
The cpc map $\phi_r$ has the ucp extension
\[
 \phi_r^{u}:D^{u}\longrightarrow M_r(B)=\mathcal L_B(B^r),
 \qquad
 \phi_r^{u}(a+\lambda1):=\phi_r(a)+\lambda1_{M_r(B)}.
\]
By the definition of unitally absorbing, applied to this map on the
Hilbert $B$-module $B^r$, there is {a family of adjointable isometries}
\[
 {(v_j:B^r\longrightarrow B\otimes H)_{j\in\nat}}
\]
such that
\[
 {\lim_{j\to\infty}\|v_j^*\pi_B(a)v_j-\phi_r^{u}(a)\|=0}
 \qquad{(\forall a\in D^{u})}.
\]
{Since $d_{n+2}(E)$ is finite and lies in $D$, choose $j$
so large that the inequality
\[
 \|v_j^*\pi_B(a)v_j-\phi_r(a)\|<\varepsilon
 \qquad(\forall a\in d_{n+2}(E))
\]
holds, and put $v:=v_j$.}
 The strict inequalities on this finite set leave room for the
following truncation and normalization without violating these estimates.

Write $r_k:=1_B\otimes p_k$. Since $B^r$ is finitely generated,
$r_kv\to v$ in norm {as $k\to\infty$}. For sufficiently large $k$ we may thus replace
$v$ by the isometry
\[
 w:=r_kv(v^*r_kv)^{-1/2},\qquad q:=ww^*\leq r_k.
\]
Its range has the finite projective orthogonal complement
$V:=(r_k-q)(B\otimes H)$. Thus $w$ determines the unitary
\[
 U:B^r\oplus V\xrightarrow{\cong}r_k(B\otimes H),
 \qquad U(\xi,\zeta):=w\xi+\zeta.
\]
Choose $k$ large enough that both
\begin{equation}\label{eq:complement-normalized-compression}
 {\|w^*\pi_B(a)w-\phi_r(a)\|<\varepsilon
 \qquad(\forall a\in d_{n+2}(E))}
\end{equation}
and
\begin{equation}\label{eq:complement-quasidiagonal-bound}
 {\|[r_k,\pi_B(a)]\|<\varepsilon
 \qquad(\forall a\in d_{n+2}(G))}
\end{equation}
hold. 

With respect to the decomposition given by $U$, the compression of
$\pi_B$ to $r_k(B\otimes H)$ has two diagonal blocks: the compression
$w^*\pi_B(-)w$ approximating $\phi_r$, and the complementary cpc map
\[
 b:D\to\mathcal L_B(V),\qquad
 b(a):=(r_k-q)\pi_B(a)(r_k-q).
\]
The full compression is obtained from the scalar cpc map
\[
 \tau:D\to M_s(\C),\qquad \tau(a):=p_k\pi(a)p_k,
 \qquad s:=\operatorname{rank}(p_k).
\]
Write $\tau_B(a):=1_B\otimes\tau(a)$, viewed on $r_k(B\otimes H)$.
We show that its off-diagonal blocks are small, so that
$U(\phi_r\oplus b)U^*$ approximates $\tau_B$.

{Since $E$ contains $G$ and its products, the approximation of $\phi_rd_{n+2}$ by the homomorphism $\theta\alpha_{n+2,m}$ on $E$ gives}
\begin{equation}\label{eq:complement-cpc-multiplicativity-defect}
 {\|\phi_rd_{n+2}(xy)-\phi_rd_{n+2}(x)\phi_rd_{n+2}(y)\|
 =O(\varepsilon)\qquad(\forall x,y\in G).}
\end{equation}
{{Estimate~\eqref{eq:complement-normalized-compression} transfers} the multiplicativity-defect bound \eqref{eq:complement-cpc-multiplicativity-defect} to the compression $w^*\pi_B(-)w$: its multiplicativity defects on $d_{n+2}(G)$ are also $O(\varepsilon)$. For $a$ in $d_{n+2}(G)$, the $C^*$-identity therefore gives}
\begin{equation}\label{eq:complement-variance-estimate}
 {\|(1-q)\pi_B(a)w\|^2
 =\|w^*\pi_B(a^*a)w-(w^*\pi_B(a)w)^*(w^*\pi_B(a)w)\|
 =O(\varepsilon).}
\end{equation}
{Since $G$ is self-adjoint, \eqref{eq:complement-variance-estimate} bounds both off-diagonal blocks by $O(\sqrt\varepsilon)$.}

For $a$ in $d_{n+2}(G)$, the four blocks of $U^*\tau_B(a)U$ on
$B^r\oplus V$ are displayed below, together with their estimates.
\begin{center}
\begingroup
\renewcommand{\arraystretch}{2}
\setlength{\extrarowheight}{1ex}
\begin{tabular}{|>{\centering\arraybackslash}m{0.43\linewidth}|
                 >{\centering\arraybackslash}m{0.43\linewidth}|}
\hline
$\begin{gathered}
 w^*\pi_B(a)w\\[1ex]
 \|w^*\pi_B(a)w-\phi_r(a)\|<\varepsilon
\end{gathered}$
&
$\begin{gathered}
 {w^*\pi_B(a)(r_k-q)}\\[1ex]
 {\|w^*\pi_B(a)(r_k-q)\|=O(\sqrt\varepsilon)}
\end{gathered}$
\\ \hline
$\begin{gathered}
 (r_k-q)\pi_B(a)w\\[1ex]
 \|(r_k-q)\pi_B(a)w\|=O(\sqrt\varepsilon)
\end{gathered}$
&
$\begin{gathered}
 {b(a)=(r_k-q)\pi_B(a)(r_k-q)}\\[1ex]
 \left\|\begin{gathered}
   bd_{n+2}(xy)\\[-0.3ex]
   {}-bd_{n+2}(x)bd_{n+2}(y)
 \end{gathered}\right\|=O(\varepsilon)\\
 {(\forall x,y\in G)}
\end{gathered}$
\\ \hline
\end{tabular}
\endgroup
\end{center}

{For $a,c$ in $d_{n+2}(G)$, the complementary compression satisfies}
\begin{equation}\label{eq:complement-compression-defect}
 {b(ac)-b(a)b(c)
 =(r_k-q)\pi_B(a)(1-r_k+q)\pi_B(c)(r_k-q).}
\end{equation}
{Using $1-r_k+q=(1-r_k)+q$, estimate~\eqref{eq:complement-variance-estimate}, {and \eqref{eq:complement-quasidiagonal-bound},} we have}
\begin{equation}\label{eq:complement-defect-factor-bounds}
 {\begin{aligned}
 \|(r_k-q)\pi_B(a)(1-r_k+q)\|&=O(\sqrt\varepsilon),\\
 \|(1-r_k+q)\pi_B(c)(r_k-q)\|&=O(\sqrt\varepsilon).
 \end{aligned}}
\end{equation}
{Since $1-r_k+q$ is a projection, the product of these two operators is the defect in \eqref{eq:complement-compression-defect}, whose norm is therefore $O(\varepsilon)$. All these bounds are uniform in $B$, $\theta$, and the finite-corner dimensions.}

{For the scalar compression $\tau(a)=p_k\pi(a)p_k$, the analogue of \eqref{eq:complement-compression-defect} gives}
\[
 \tau(ac)-\tau(a)\tau(c)
 =p_k\pi(a)(1-p_k)\pi(c)p_k.
\]
Both factors are small by quasidiagonality. Hence the
multiplicativity defects of $bd_{n+2}$ and $\tau d_{n+2}$ on $G$
are less than $\gamma$ for sufficiently small $\varepsilon$.
{Combining \eqref{eq:complement-normalized-compression} and \eqref{eq:complement-variance-estimate}, we obtain}
\begin{equation}\label{eq:complement-block-comparison}
 {\lim_{\varepsilon\to0}\max_{a\in d_{n+2}(G)}
 \|U(\phi_r(a)\oplus b(a))U^*-\tau_B(a)\|=0.}
\end{equation}
{This compares the cpc maps before straightening.}

\smallskip
\noindent\emph{4. Straightening and the required homotopy.}
Apply the straightening controls of Step~1 to $bd_{n+2}$ and $\tau d_{n+2}$.
We obtain {morphisms}
\[
 \eta':D_{n+1}\to\mathcal L_B(V),\qquad
 \tau':D_{n+1}\to M_s(\C)
\]
which approximate $bd_{n+1}$ and $\tau d_{n+1}$ on $F$, respectively,
within $\kappa$. Hence, for $x$ in $F$,
\begin{equation}\label{eq:complement-straightened-comparison}
 {\|U(\rho(x)\oplus\eta'(x))U^*-\tau'_B(x)\|<2\kappa+\|U(\phi_r d_{n+1}(x)\oplus bd_{n+1}(x))U^*-\tau_B d_{n+1}(x)\|.}
\end{equation}
Here one error $\kappa$ comes from replacing the direct sum
$\rho\oplus\eta'$ by its cpc approximation (the norm is the maximum
of the two summand errors); the other comes from replacing
$\tau'_B$ by the scalar compression.

We can now specify the choices. Choose $\kappa$ so that
$2\kappa<\delta$, and take the corresponding $G$, $E$ and $\gamma$
from Step~1. Then choose $\varepsilon$ small enough that all
multiplicativity defects lie below $\gamma$, the {norm in \eqref{eq:complement-block-comparison}} is smaller than $\delta-2\kappa$, and the cpc approximation
error in \eqref{eq:complement-first-approximation} is smaller than
$\delta-\kappa$. Steps~2--3 show that these errors tend to zero {as $\varepsilon\to0$}, so all requirements hold simultaneously. These
choices first fix the controls and then the stage $m$ by
Lemma~\ref{lem:general-cp}, with accuracy $\varepsilon$, uniformly for all target algebras.
Only after fixing this stage do we take {$B$ in $\Calg_{\sepa}$} and a
{morphism} $\theta:D_m\to B\otimes\bbK$ and carry out the
construction above. Both comparisons then lie strictly within
$\delta$ on $F$.

{The homotopy control of Step~1, applied to \eqref{eq:complement-straightened-comparison} and \eqref{eq:complement-first-approximation}, now gives}
\[
 \theta\alpha_{nm}\simeq\rho\alpha_{n,n+1},\qquad
 (\rho\oplus\eta')\alpha_{n,n+1}
 \simeq\tau'_B\alpha_{n,n+1}.
\]
Here the module identification $U$ can be omitted after stabilization:
two isometric embeddings of a finite projective module become
homotopic by rotating between their copies in an orthogonal direct
sum. View $\mathcal L_B(V)$ as the corner
$(r_k-q)M_s(B)(r_k-q)\subseteq B\otimes\bbK$ and set
\[
 \eta:=\eta'\alpha_{n,n+1}.
\]
Adding $\eta$ to the first homotopy and using the second gives
\[
 \theta\alpha_{nm}\oplus\eta
 \simeq\tau'_B\alpha_{n,n+1}.
\]
Thus the scalar {morphism} $\tau'\alpha_{n,n+1}:D_n\to\bbK$
provides the {morphism} required in the statement.
\end{proof}


{Let $u:P\to Q$ be a semiprojective morphism in $\nCalg_{\sepa}$ and let $q:E\twoheadrightarrow F$ be a surjective morphism in $\nCalg$. Let $h:Q\to F$ be a homomorphism and let $\Phi_t:Q\to E$ be an asymptotic homomorphism satisfying $q\Phi_t=h$.}
\begin{lem}\label{lem:relative-asymptotic-straightening}
{The following lifting problem admits a homomorphism}
$\widetilde h$ making the diagram commute:
\[
 \xymatrix@C=3em@R=3em{
 P\ar[r]^u\ar@{.>}[d]_{\widetilde h}
   &Q\ar[d]^h\ar@{~>}[dl]^{\Phi_t}\\
 E\ar@{->>}[r]_q&F.
 }
\]
The wavy diagonal denotes the given asymptotic morphism
$\Phi_t$, with $q\Phi_t=h$. 
\end{lem}
\begin{proof}
Consider the algebra of bounded continuous paths lifting the constant
quotient values:
\[
 \mathcal B:=\{(b,a)\in C_b([1,\infty),E)\oplus Q:
                         q(b(t))=h(a)\text{ for all }t\}.
\]
The ideal $J:=C_0([1,\infty),\ker q)$ is the closure of the increasing
ideals of paths supported in bounded initial intervals. The
asymptotic lift induces a genuine homomorphism
$Q\to\mathcal B/J$ whose second coordinate is the identity.
Semiprojectivity lifts its composition with $u$ modulo one of these
smaller ideals. Evaluation sufficiently far along the path gives a
homomorphism $P\to E$ with quotient $hu$.
\end{proof}

{Let $u:P\to Q$ be a semiprojective morphism  and $q:E\twoheadrightarrow F$ be a surjection  in $\nCalg_{\sepa}$, and let $(K,L)$ be a finite CW pair. Put $Y:=K\times[0,1]$ and $Z:=K\times\{0\}\cup L\times[0,1]$. Let}
\[
 {H:Q\to C(Y,F),\qquad \psi:Q\to C(Z,E)}
\]
{be homomorphisms satisfying $H|_Z=q\psi$.}
\begin{lem}\label{lem:general-relative-lifting}
{The following lifting problem admits a homomorphism}
$\widetilde H$ making the diagram commute:
\[
 \xymatrix@C=3em@R=3em{
 {P}\ar[r]^{{u}}\ar@{.>}[d]_{\widetilde H}
   & \ar@{~~>}[dl]^{\Phi_t} {Q}\ar[d]^{(H,\psi)}\\
 C(Y,E)\ar@{->>}[r]_-{(q,\operatorname{res}_Z)}
   & C(Y,F)\times_{C(Z,F)}C(Z,E).
 }
\]
\end{lem}
\begin{proof}
{The natural morphism
\[
(q,\operatorname{res}_Z):C(Y,E)\longrightarrow
 C(Y,F)\times_{C(Z,F)}C(Z,E)
\]
is surjective. Indeed, the Bartle--Graves theorem
\cite{bartle-graves} gives a continuous section $s:F\to E$ of
$q$. Choose a retraction $r:Y\to Z$ arising from a deformation
retraction. For a pair $(f,g_Z)$ in the pullback, the continuous
function
\[
 g:=s\circ f+\bigl(g_Z-(s\circ f)|_Z\bigr)\circ r:Y\to E
\]
satisfies $qg=f$ and $g|_Z=g_Z$, since the correction takes
values in $\ker q$.}

{Let $(r_t:Y\to Y)_{t\in[0,1]}$ be the deformation from
$\id_Y$ to the composite $Y\xrightarrow{r}Z\hookrightarrow Y$,
fixing $Z$. The morphisms
\[
 a\longmapsto\bigl(H(a)\circ r_t,\psi(a)\bigr)
\]
give a homotopy from $(H,\psi)$ to
$(r^*(q\psi),\psi)$. The latter has the genuine lift
\[
 r^*\psi:{Q}\longrightarrow C(Y,E).
\]
We apply Shulman's homotopy lifting theorem
\cite[Theorem~11]{shulman} to the quotient morphism $(q,\operatorname{res}_Z)$ and
these two homotopic morphisms  \[
{(H,\psi),(r^*(q\psi),\psi):Q\longrightarrow
 C(Y,F)\times_{C(Z,F)}C(Z,E).}
\]
It gives an
asymptotic homomorphism
\[
 \Phi_t:{Q}\longrightarrow C(Y,E),\qquad
 (q,\operatorname{res}_Z)\Phi_t=(H,\psi)\quad(\forall t\geq1).
\]
Applying Lemma~\ref{lem:relative-asymptotic-straightening} to
${u}$ and $(q,\operatorname{res}_Z)$ gives the required morphism
$\widetilde H$ with
$(q,\operatorname{res}_Z)\widetilde H=(H,\psi){u}$.}
\end{proof}

\subsection{Unital coefficients}
\label{sec:unital-coefficients}

{Let $B$ be in $\nCalg_{\sepa}$.}
Let $p$ be a rank-one projection in $\bbK$. Using the unitalization inclusion \eqref{eq:unitalization-inclusion}, we form the natural homomorphism
\[
B\xrightarrow{\iota_{B\otimes \bbK}(\id\otimes p)\otimes p} (B\otimes \bbK)^{u}\otimes \bbK
\]
Its composite with the scalar projection $\epsilon_{B\otimes\bbK}\otimes \id_{\bbK}$ vanishes and therefore induces a map 
\[
\kappa_{D,B}:{\cM}_{D}(B)^{\group}\to  \cT_{D}(B)\ .
\]%
 \begin{lem}\label{lem:group-scalar-reduction} If $B$ is unital, then
for a family $(D_{n})$ in $\nCalg_{\sepa}$ with semiprojective transition maps  
the family $(\pi_{0}\kappa_{D_{n},B})$ induces an
isomorphism of pro-abelian groups
\[
 \bigl(\pi_0({\cM}_{D_n}(B)^{\group})\bigr)_n
 \xrightarrow{\cong}\bigl(\pi_0\cT_{D_n}(B)\bigr)_n.
\]
\end{lem}
\begin{proof}
  We set $\kappa_{n}:=\pi_{0}\kappa_{D_{n},B}$.

We show that the comparison has trivial pro-kernel and
pro-cokernel. For this, it suffices to verify that for every $n$ in $\nat$ there exists $m$ in $\nat$ with $m\geq n+2$ such that
\begin{equation}\label{eq:group-scalar-pro-criteria}
 \begin{aligned}
 \operatorname{im}\left(
 \pi_0\cT_{D_m}(B)\xrightarrow{\alpha_{nm}^*}
 \pi_0\cT_{D_n}(B)\right)
 &\subseteq\operatorname{im}(\kappa_n),\\
 \alpha_{nm}^*(\ker\kappa_m)&=0.
 \end{aligned}
\end{equation}
We prove \eqref{eq:group-scalar-pro-criteria} with $m=n+4$.

{Fix $n$ in $\nat$. Lemma~\ref{lem:semiprojective-straightening}(2), applied to $\alpha_{n,n+1}$, gives a finite subset $F$ of the unit ball of $D_{n+1}$ and a tolerance $\delta$ in $(0,\infty)$ such that, for every $C$ in $\nCalg_{\sepa}$ and homomorphisms $h_0,h_1:D_{n+1}\to C$,
\begin{equation}\label{eq:group-scalar-homotopy-control}
 \|h_0(a)-h_1(a)\|<\delta\ (\forall a\in F)
 \quad\Longrightarrow\quad
 h_0\alpha_{n,n+1}\simeq h_1\alpha_{n,n+1}.
\end{equation}}{The transition $\alpha_{n+1,n+2}$ is semiprojective. Choose $\varepsilon$ in $(0,\infty)$ with $2\varepsilon<\delta$. Lemma~\ref{lem:semiprojective-straightening}(1), applied to $\alpha_{n+1,n+2}$ with norm bound $1$ and straightening accuracy $\varepsilon$ on $F$, gives a finite subset $G$ of $D_{n+2}$ and a tolerance $\gamma$ in $(0,\infty)$. For every $C$ in $\nCalg_{\sepa}$, a cpc map $f:D_{n+2}\to C$ satisfying
\begin{equation}\label{eq:group-scalar-multiplicativity-control}
 \|f(ab)-f(a)f(b)\|<\gamma\qquad(\forall a,b\in G)
\end{equation}
admits a homomorphism $h:D_{n+1}\to C$ with
\begin{equation}\label{eq:group-scalar-straightening-control}
 \|h(a)-f(\alpha_{n+1,n+2}(a))\|<\varepsilon
 \qquad(\forall a\in F).
\end{equation}
We enlarge $G$ to contain $\alpha_{n+1,n+2}(F)$.}

 {The following construction will be called  the outer compression. Let $\varphi:D_{n+4}\to(B\otimes\bbK)^u\otimes\bbK$ be a homomorphism. Choose finite-rank projections $p_r$ increasing strongly to $1$ in the outer $\bbK$-factor. The cpc maps
\begin{equation}\label{eq:group-outer-compression}
 \varphi_r:D_{n+4}\longrightarrow M_{\operatorname{rank}(p_r)}((B\otimes\bbK)^u),
 \qquad \varphi_r(a):=(1\otimes p_r)\varphi(a)(1\otimes p_r)
\end{equation}
converge pointwise in norm to $\varphi$, after inclusion into the stabilized target. Their multiplicativity defects therefore tend to zero as $r\to\infty$.}

 {Apply Lemma~\ref{lem:semiprojective-straightening}(2) to $\alpha_{n+2,n+3}$ to choose a finite control set and a tolerance. Then apply part (1) to $\alpha_{n+3,n+4}$, with norm bound $1$ and straightening accuracy smaller than this tolerance on the chosen set. For sufficiently large $r$, straightening \eqref{eq:group-outer-compression} gives a homomorphism from $D_{n+3}$ into the finite outer corner. The homotopy control gives, after restriction to $D_{n+2}$, a homomorphism $\varphi':D_{n+2}\to M_{\operatorname{rank}(p_r)}((B\otimes\bbK)^u)$ with
\begin{equation}\label{eq:group-outer-corner-reduction}
 \varphi\alpha_{n+2,n+4}\simeq\varphi'
 \quad\text{in }\Hom_{\nCalg_{\sepa,h}}
       (D_{n+2},(B\otimes\bbK)^u\otimes\bbK).
\end{equation}
The control sets and tolerances depend only on the two transitions. For a compact parameter family, we apply the same argument with a continuous-function algebra as target.}

{We now consider the pro-cokernel. {Since the scalar projection has a section, the fibre sequence defining $\cT_{D_{n+4}}(B)$, given by \eqref{hetrgertget} and Definition~\ref{def:intro-T}, splits and gives}
\begin{equation}\label{eq:group-scalar-kernel-identification}
{\pi_0\cT_{D_{n+4}}(B)\cong
\ker\!\left(
\pi_0\bigl({\cM}_{D_{n+4}}((B\otimes\bbK)^u)^{\group}\bigr)
\xrightarrow{(\epsilon_{B\otimes\bbK})_*}
\pi_0\bigl({\cM}_{D_{n+4}}(\C)^{\group}\bigr)
\right).}
\end{equation}
{Using \eqref{eq:group-scalar-kernel-identification}, we write a class as a difference of classes of homomorphisms whose scalar images agree after group completion.} Applying the outer compression construction  to both representatives, we can represent its restriction to stage $n+2$ by $x=[\varphi]-[\psi]$, where
\[
 \varphi,\psi:D_{n+2}\longrightarrow M_r((B\otimes\bbK)^u)
\]
for a common $r$ in $\nat$. Their scalar quotients have equal classes in $\pi_0({\cM}_{D_{n+2}}(\C)^{\group})$. It suffices to show that $\alpha_{n,n+2}^*(x)$ lies in the image of $\kappa_n$.}

{We describe the inner compression for $\varphi$; the same construction applies to $\psi$. Denote the scalar quotient of $\varphi$ by $\bar\varphi:D_{n+2}\to M_r(\C)$.}

We choose finite-rank projections $e_j$ increasing strongly to $1$
in the inner $\bbK$-factor. Since $B$ is unital, the projections $1_B\otimes e_j$ form an approximate unit of $B\otimes\bbK$ and commute with the scalar part of $\varphi$. For every $a$ in $D_{n+2}$, we have
\[
 \varphi(a)-1\otimes\bar\varphi(a)\in M_r(B\otimes\bbK).
\]
{We consider the  cpc maps
\begin{equation}\label{eq:group-finite-compression}
 f_j:D_{n+2}\longrightarrow M_{\operatorname{rank}(e_j)}(B)\otimes M_r(\C),
 \qquad f_j(a):=(1_B\otimes e_j\otimes1_r)\varphi(a)(1_B\otimes e_j\otimes1_r).
\end{equation}}
We can choose $j$ so large such that
\begin{equation}\label{eq:group-compression-defect-control}
 \|f_j(ab)-f_j(a)f_j(b)\|<\gamma
 \qquad(\forall a,b\in G)
\end{equation}
and
\begin{equation}\label{eq:group-compression-comparison-control}
 \bigl\|\varphi(\alpha_{n+1,n+2}(a))-
 \bigl(f_j(\alpha_{n+1,n+2}(a))\oplus((1-1_B\otimes e_j)\otimes\bar\varphi(\alpha_{n+1,n+2}(a)))\bigr)\bigr\|
 <\delta-\varepsilon\quad(\forall a\in F).
\end{equation}
  Applying Lemma~\ref{lem:semiprojective-straightening}(1) with straightening accuracy $\varepsilon$, using \eqref{eq:group-compression-defect-control}, gives a homomorphism
\[
 h':D_{n+1}\longrightarrow M_{\operatorname{rank}(e_j)}(B)\otimes M_r(\C)
\]
such that
\begin{equation}\label{eq:group-straightened-compression-estimate}
 \|h'(a)-f_j(\alpha_{n+1,n+2}(a))\|<\varepsilon
 \qquad(\forall a\in F).
\end{equation}
Combining \eqref{eq:group-straightened-compression-estimate} and \eqref{eq:group-compression-comparison-control} gives
\[
 \bigl\|\varphi(\alpha_{n+1,n+2}(a))-
 \bigl(h'(a)\oplus((1-1_B\otimes e_j)\otimes\bar\varphi(\alpha_{n+1,n+2}(a)))\bigr)\bigr\|
 <\delta\quad(\forall a\in F).
\]
Applying \eqref{eq:group-scalar-homotopy-control} and putting $h:=h'\alpha_{n,n+1}$ gives
\[
 \varphi\alpha_{n,n+2}\simeq
 h\oplus\bigl((1-1_B\otimes e_j)\otimes\bar\varphi\alpha_{n,n+2}\bigr).
\]
The scalar representation $1\otimes\bar\varphi\alpha_{n,n+2}$
has finite part
$(1_B\otimes\bar\varphi\alpha_{n,n+2})^{\oplus\operatorname{rank}(e_j)}$
and the same complementary summand as the displayed homotopy.
Subtracting this scalar correction cancels the complementary
terms and gives
\[
 \alpha_{n,n+2}^*
 \bigl([\varphi]-[1\otimes\bar\varphi]\bigr)
 =
 \kappa_n\left(
 [h]-[(1_B\otimes\bar\varphi\alpha_{n,n+2})^{\oplus\operatorname{rank}(e_j)}]
 \right).
\]
We view the left side in the scalar kernel defining
$\pi_0\cT_{D_n}(B)$. Applying the same construction to $\psi$
gives the analogous formula. Since the scalar quotients of
$\varphi$ and $\psi$ have equal classes, their scalar
corrections cancel when we subtract these formulas.
Thus $\alpha_{n,n+2}^*(x)$ lies in the image of $\kappa_n$,
proving the first assertion of
\eqref{eq:group-scalar-pro-criteria}.

{For the pro-kernel, let $x$ be in $\ker\kappa_{n+4}$. Write $x=[\rho_0]-[\rho_1]$ with $\rho_i:D_{n+4}\to B\otimes\bbK$. Since the scalar fibre injects into $\pi_0({\cM}_{D_{n+4}}((B\otimes\bbK)^u)^{\group})$, there are a common representation $\psi:D_{n+4}\to(B\otimes\bbK)^u\otimes\bbK$ and a homotopy between the images of $\rho_0\oplus\psi$ and $\rho_1\oplus\psi$.}

{We first apply outer compression to the given homotopy and to $\psi$, using the same finite outer corner and retaining the block decomposition at the endpoints. Straightening both compressions along $\alpha_{n+3,n+4}$ gives homomorphisms}
\begin{equation}\label{eq:group-scalar-outer-straightenings}
{H'':D_{n+3}\to C([0,1],M_r((B\otimes\bbK)^u)),\qquad
\psi'':D_{n+3}\to M_s((B\otimes\bbK)^u).}
\end{equation}
{We choose the approximations sufficiently close that $\operatorname{ev}_iH''$ and $\rho_i\alpha_{n+3,n+4}\oplus\psi''$ satisfy the closeness condition of Lemma~\ref{lem:semiprojective-straightening}(2) for $\alpha_{n+2,n+3}$. That lemma gives homotopies}
\begin{equation}\label{eq:group-scalar-outer-endpoint-homotopies}
{(\operatorname{ev}_iH'')\alpha_{n+2,n+3}
\simeq(\rho_i\alpha_{n+3,n+4}\oplus\psi'')\alpha_{n+2,n+3}
\qquad(i=0,1).}
\end{equation}
{Concatenating \eqref{eq:group-scalar-outer-endpoint-homotopies} with $H''\alpha_{n+2,n+3}$ gives a homomorphism $H':D_{n+2}\to C([0,1],M_r((B\otimes\bbK)^u))$ with, on putting $\psi':=\psi''\alpha_{n+2,n+3}$,}
\begin{equation}\label{eq:group-scalar-outer-endpoints}
{\operatorname{ev}_iH'=\rho_i\alpha_{n+2,n+4}\oplus\psi'
\qquad(i=0,1).}
\end{equation}
{Here and below we include the matrix corners into the stabilized targets.}

{We now compress $H'$ and $\psi'$ in the inner $\bbK$-factor, using the same projection $e_j$, compatibly with the direct sums in \eqref{eq:group-scalar-outer-endpoints}. We choose $j$ uniformly over the path parameter, so that both compressions satisfy \eqref{eq:group-scalar-multiplicativity-control}. Straightening with accuracy $\varepsilon$ gives homomorphisms $h':D_{n+1}\to C([0,1],B\otimes\bbK)$ and $\eta':D_{n+1}\to B\otimes\bbK$. We use the same $\eta'$ at both endpoints. By \eqref{eq:group-scalar-outer-endpoints}, the compressed endpoints are the direct sums of the compressions of $\rho_i\alpha_{n+2,n+4}$ and $\psi'$. Since compression of each $\rho_i$ converges in norm, we can choose $j$ so that its error on the image of $F$ is less than $\delta-2\varepsilon$. The two straightening errors then give}
\begin{equation}\label{eq:group-scalar-inner-endpoints}
{\|\operatorname{ev}_ih'(a)-
 (\rho_i\alpha_{n+1,n+4}(a)\oplus\eta'(a))\|<\delta
 \qquad(\forall a\in F,\ i=0,1).}
\end{equation}
{We apply \eqref{eq:group-scalar-homotopy-control} to both endpoints. Concatenating the resulting homotopies with $h'\alpha_{n,n+1}$ gives, with $\eta:=\eta'\alpha_{n,n+1}$,}
\[
 \rho_0\alpha_{n,n+4}\oplus\eta
 \simeq\rho_1\alpha_{n,n+4}\oplus\eta.
\]
Cancellation of $[\eta]$ in the group completion gives
${\alpha_{n,n+4}^*(x)}=0$, which is the second assertion of
\eqref{eq:group-scalar-pro-criteria}.
 \end{proof}

\subsection{Scalar localization}

 In this subsection, for $D$ in $\nCalg_{\sepa}$,   we restrict the functor $B\mapsto {\cM}_{D}(B)$ to the category  $\Calg_{\sepa}$ of unital algebras and unit-preserving morphisms. The units
$\C\to B$ induce a natural transformation, called  scalar extension,
\[\sigma_{D}: \underline{{\cM}_{D}(\C)}\to {\cM}_{D}(-)\ , \qquad 
 \sigma_{D,B}:{\cM}_{D}(\C)\longrightarrow {\cM}_{D}(B),
 \qquad \tau\longmapsto\tau_B,
\] see also \eqref{fqqwedweqfqwe}.
We call the elements of
$\operatorname{im}(\pi_0\sigma_{D,B})$ the scalar components of $\pi_{0}{\cM}_{D}(B)$.

{Let $(D_n)_{n\in\nat}$ be a shape system of $D$ in $\nCalg$ in the
sense of Definition~\ref{def:intro-shape-system}.}
We abbreviate $\sigma_{n,B}:=\sigma_{D_{n},B}$.
For each $B$, let  \[
 {{\cM}_{D_{n}}(B)\to L_n(B):={\cM}_{D_n}(B)
       [\operatorname{im}(\pi_0\sigma_{n,B})^{-1}]}
\]  be the localization in
$\CMon(\Ani)$ which inverts the scalar components of ${\pi_{0}{\cM}_{D_n}(B)}$.
 By the universal property of localization, group completion
factors canonically as
\[
 {\cM}_{D_n}(B)\longrightarrow L_n(B)\longrightarrow {\cM}_{D_n}(B)^{\group}.
\]

\begin{lem}[Scalar localization and group completion]\label{lem:general-scalar-localization}
\begin{enumerate}
\item The pointwise scalar localizations assemble into functors
\[
 L_n:\Calg_{\sepa}\longrightarrow\CMon(\Ani)
\]
with natural transformations ${\cM}_{D_n}\to L_n\to {\cM}_{D_n}^{\group}$.
{Precomposition with the transition maps makes these transformations compatible with the tower maps.}
\item {If $D$ is nuclear and quasidiagonal, then the comparison}
\[
 (L_n)_n\longrightarrow({\cM}_{D_n}^{\group})_n
\]
is a pro-equivalence of functors.
\end{enumerate}
\end{lem}
\begin{proof}
\noindent{\emph{1.}} For any $A$ in $\nCalg_{\sepa}$,
the map ${\cM}_{A}(B)\to {\cM}_{A}(B')$ induced by a unital map $f:B\to B'$ preserves scalar components since 
\[
 {\cM}_{A}(f)\circ\sigma_{A,B}\simeq \sigma_{A,B'}.
\]
 We therefore obtain the functor ${L_{n}:\Calg_{\sepa}\to \CMon(\Ani)}$ as a  pointwise localization at the scalar components  of the  functor  ${\cM}_{D_{n}}(-)$.


Restriction along $\alpha_{nm}$ also preserves scalar components,
since it sends $\tau_B$ to $(\tau\alpha_{nm})_B$. Consequently the
restriction transformations descend to $L_m\to L_n$.
Group completion inverts all components, so the comparison
${\cM}_{D_n}\to {\cM}_{D_n}^{\group}$ factors naturally through $L_n$.

\noindent{\emph{2.}} Let $L_n^\times(B)$ be the union of the invertible
components of $L_n(B)$. Monoid morphisms preserve invertibility,
so these subspaces define a functor $L_n^\times$ with values in
$\CGrp(\Ani)$ and a natural inclusion {$L_n^\times\to L_n$}.

Fix $n$ in $\nat$ and choose $m$ in $\nat$ with $m\geq n$ as in
Proposition~\ref{prop:general-complement}. For every unital $B$
and every representation $\theta:D_m\to B\otimes\bbK$, that
proposition gives
\[
 [\theta\alpha_{nm}]+[\eta]=[\tau_B]
 \qquad\text{in }\pi_0L_n(B).
\]
Since $\tau_{B}$ is scalar, $[\tau_B]$ is invertible in $\pi_{0}L_{n}(B)$. Consequently,   $[\theta\alpha_{nm}]$ is invertible
with inverse $[\eta]-[\tau_B]$. Thus the natural transformation
${\cM}_{D_m}\to L_n$ takes values in $L_n^\times$.
Here $m$ depends only on $n$ and not on $B$. The factorization is therefore
natural in $B$.
 
{The map ${\cM}_{D_m}\to L_n^\times$ inverts the scalar
components, so the universal property of pointwise localization
gives the dotted factorization in the commutative diagram
\begin{equation}\label{eq:scalar-localization-units-factorization}
{\xymatrix@C=4em@R=3em{
 L_m^\times\ar[r]\ar@{^{(}->}[d]
   &L_n^\times\ar@{^{(}->}[d]\\
 L_m\ar[r]\ar@{..>}[ur]
   &L_n.
 }}
\end{equation}
The horizontal arrows are restriction maps. The lower triangle
commutes by the localization property; the upper triangle
commutes because $L_n^\times\to L_n$ is the inclusion of the
invertible components. Uniqueness of these factorizations makes
them compatible with further restriction. Since $m$ depends only
on $n$, the dotted arrows therefore give the reverse map of
pro-systems in
\[
 \xymatrix@C=5em{
 (L_n^\times)_{n\in\nat}\ar@<0.6ex>[r]
   &(L_n)_{n\in\nat}\ar@<0.6ex>[l].
 }
\]
{The triangle with vertices $L_m^\times$, $L_m$, and
$L_n^\times$ in \eqref{eq:scalar-localization-units-factorization}
shows that the composite from $(L_n^\times)_{n\in\nat}$ to
$(L_n)_{n\in\nat}$ and back is represented by the restriction
$L_m^\times\to L_n^\times$. The triangle with vertices $L_m$,
$L_n^\times$, and $L_n$ shows that the other composite is
represented by $L_m\to L_n$. Both composites are therefore
identities of pro-systems, so the inclusion is a pro-equivalence.}}
Applying pointwise group completion preserves this
pro-equivalence. Since $L_n^\times$ is group-like and scalar
localization does not change group completion, we obtain the
chain of equivalences
\begin{equation}\label{eq:scalar-localization-pro-equivalences}
 {(L_n)_n
 \xleftarrow{\sim}(L_n^\times)_n
 \xrightarrow{\sim}\bigl((L_n^\times)^{\group}\bigr)_n
 \xrightarrow{\sim}(L_n^{\group})_n
 \xrightarrow{\sim}({\cM}_{D_n}^{\group})_n}
\end{equation}
in $\Pro(\Fun(\Calg_{\sepa},\CMon(\Ani)))$.
{By naturality of group completion, \eqref{eq:scalar-localization-pro-equivalences} identifies the}
canonical comparison $(L_n)_n\to({\cM}_{D_n}^{\group})_n$ as an
equivalence.
\end{proof}

{Let $S_n$ be the set of scalar morphisms $\tau:D_n\to\bbK$,
and let $\operatorname{Fin}(S_n)$ be the filtered poset of finite
subsets of $S_n$, ordered by inclusion. For $T$ in
$\operatorname{Fin}(S_n)$, let
\[
 L_{n,T}(B):={\cM}_{D_n}(B)[\tau_B^{-1}:\tau\in T]
\]
be the pointwise localization in $\CMon(\Ani)$. For $B\to B'$ and
$T\subseteq T'$, its universal property gives a natural map
$L_{n,T}(B)\to L_{n,T'}(B')$. The compatibility of  these localizations  shows that we obtain a functor
{\[
 L_{n,-}:\operatorname{Fin}(S_n)\longrightarrow
 \Fun(\Calg_{\sepa},\CMon(\Ani)).
\]}as a {pointwise} localization of the  constant functor with value ${\cM}_{D_{n}}(-)$.
 For each
fixed $T$, {choose a finite block sum}
\[
 \tau_T:=\bigoplus_{\tau\in T}\tau:D_n\longrightarrow\bbK.
\]}

\begin{lem}\label{lem:general-scalar-translations}
\begin{enumerate}
\item The canonical comparison is an equivalence
\[
 \colim_{T\in\operatorname{Fin}(S_n)}L_{n,T}
 \xrightarrow{\simeq}L_n
 \qquad\text{in }\Fun(\Calg_{\sepa},\Ani).
\]
\item For each finite subset $T$ of $S_n$, there is a natural equivalence
\begin{equation}\label{eq:scalar-localization-telescope}
 {\colim\bigl({\cM}_{D_n}\xrightarrow{-\oplus\tau_{T,-}}{\cM}_{D_n}
       \xrightarrow{-\oplus\tau_{T,-}}\cdots\bigr)
 \xrightarrow{\simeq}L_{n,T}}
\end{equation}
in $\Fun(\Calg_{\sepa},\Ani)$.
\item For $n,m$ in $\nat$ with $n\leq m$, restriction along
$\alpha_{nm}$ induces the vertical maps in the commutative diagram
\[
 \xymatrix@C=4em@R=3em{
 \colim_{T\in\operatorname{Fin}(S_m)}L_{m,T}
       \ar[r]^-{\simeq}\ar[d]
   &L_m\ar[d]\\
 \colim_{T\in\operatorname{Fin}(S_n)}L_{n,T}
       \ar[r]^-{\simeq}
   &L_n.
 }
\]
\end{enumerate}
\end{lem}
\begin{proof}
\noindent\emph{1.} {The filtered colimit of the finite localizations inverts
exactly the scalar components defining $L_n$. Its universal
property therefore identifies it with $L_n$ as a functor to
$\CMon(\Ani)$. {Since the forgetful functor $\CMon(\Ani)\to\Ani$
preserves filtered colimits, this also proves the asserted equivalence
in anima-valued functors.}}

\noindent\emph{2.} {Fix a finite subset $T$ of $S_n$. Inverting the sum
$\tau_{T,B}$ is equivalent to inverting every $\tau_B$ with
$\tau$ in $T$.  {In the operadic model of Remark~\ref{kheoprhtertg},
the cyclic permutation automorphism of
$\tau_{T,B}\oplus\tau_{T,B}\oplus\tau_{T,B}$ is homotopic to the
identity automorphism. Indeed, a path from the cyclic permutation
matrix to the identity in $SU(3)$ commutes throughout with the triple
block representation, since its three summands are equal.} Thus
\cite[Proposition~5.1]{bekhsy} identifies localization at
$\tau_{T,B}$ with {the telescope in \eqref{eq:scalar-localization-telescope}}.  The formation of the telescopes  is
natural in $B$, since scalar extension commutes with morphisms
in $\Calg_{\sepa}$. Note that the telescope is formed in anima-valued
functors; its translation maps need not preserve the monoid unit.
}

\noindent\emph{3.} {We have a functor $\alpha_{nm}^{*} :\operatorname{Fin}(S_m)\to \operatorname{Fin}(S_n)$ sending a finite subset $T$ of $S_m$
to the finite subset
$
 \{\tau\alpha_{nm}:\tau\in T\}\subseteq S_n
$.
We get a commutative diagram
\[\xymatrix@C=4em@R=3em{
 L_{m,-}
       \ar[r]\ar[d]
   &\underline{L_m}\ar[d]\\
  L_{n,\alpha_{nm}^{*}}
       \ar[r]
   &\underline{L_n}.
 }
\]
{in
\[
 \Fun\bigl(\operatorname{Fin}(S_m),
       \Fun(\Calg_{\sepa},\CMon(\Ani))\bigr).
\]}
{Taking colimits over $\operatorname{Fin}(S_m)$ gives}
\[
 \xymatrix@C=4em@R=3em{
 \colim_{T\in\operatorname{Fin}(S_m)}L_{m,T}
       \ar[r]^{\simeq}\ar[d]
   &L_m\ar[d]\\
 \colim_{T\in\operatorname{Fin}(S_m)}L_{n,\alpha_{nm}^{*}T}
       \ar[r]
   &L_n.
 }
\]
{We then compose with}
\[
\xymatrix{ \colim_{T\in\operatorname{Fin}(S_m)}L_{n,\alpha_{nm}^{*}T}
\ar[dr]\ar[rr]&&\colim_{T\in\operatorname{Fin}(S_n)}L_{n,T}\ar[dl]_{\simeq}\\&L_{n}&}
\]%
to get the asserted square.
}
\end{proof}

\subsection{Cartesian squares preserved by ${\cR}_{D}$}
\label{sec:resolved-cartesian-proof}

\begin{proof}[Proof of Theorem~\ref{thm:general-resolved-exactness}]
\emph{Strategy.}
For a cartesian square $(A_{ij})$ of separable $C^*$-algebras with
surjective right vertical map, we must prove that the comparison
\[
 {\cR}_D(A_{00})\longrightarrow
 {\cR}_D(A_{01})\times_{{\cR}_D(A_{11})}{\cR}_D(A_{10})
\]
is an equivalence in $\Ani$. This also proves the assertion in
$\CGrp(\Ani)$, since its forgetful functor to $\Ani$ creates limits
and detects equivalences. 

{Choose a shape system $(D_n)$ of $D$, with transition maps $\alpha_{nm}:D_n\to D_m$. Formula~\eqref{eq:R-shape-formula} and commutation of limits give}
\begin{equation}\label{eq:exactness-resolved-fibre}
{ {\cR}_D(A)\simeq\Fib_0\left(
 \lim_n {\cM}_{D_n}((A\otimes\bbK)^{u})^{\group}
 \longrightarrow\lim_n {\cM}_{D_n}(\C)^{\group}\right).}
\end{equation}
We will establish the required equivalence in three steps.
\begin{enumerate}
\item For a unital cartesian square $(B_{ij})$ with surjective
right vertical map, we prove that the scalar-localized comparison
\[
 L_n(B_{00})\longrightarrow
 L_n(B_{01})\times_{L_n(B_{11})}L_n(B_{10})
\]
induces an equivalence on inverse limits by verifying the
shifted lifting property.
\item Lemma~\ref{lem:general-scalar-localization} identifies the
limits of the scalar localizations with those of the stagewise
group completions. {The shifted lifting property of \cite[Definitions~2.1 and~2.29]{budu}, verified in Step~1, gives an equivalence on inverse limits. Since limits commute with pullbacks, we conclude that $\lim_nM_{D_n}^{\group}$ sends the unital square to a pullback in $\Ani$.}
\item Apply this to $B_{ij}:=(A_{ij}\otimes\bbK)^{u}$. Stabilization
and unitalization preserve the given cartesian square and the
surjectivity of its right vertical map. Taking fibres over zero of
the compatible maps to $\lim_nM_{D_n}(\C)^{\group}$ then gives the required pullback
for ${\cR}_D$, again because limits commute.
\end{enumerate}

\emph{Step 1. The scalar-localized comparison.}
Consider a unital cartesian square
\begin{equation}\label{gwregwerfwrf} 
 \xymatrix{
 B_{00}\ar[r]\ar[d]&B_{01}\ar[d]^q\\
 B_{10}\ar[r]^f&B_{11}}
 \qquad(q\text{ surjective}).
 \end{equation} Form the pullbacks in $\Ani$
\begin{equation}\label{eq:exactness-comparison-pullbacks}
{ \begin{split}
 P_n&:={\cM}_{D_n}(B_{01})
       \times_{{\cM}_{D_n}(B_{11})}{\cM}_{D_n}(B_{10}),\\
 Q_n&:=L_n(B_{01})\times_{L_n(B_{11})}L_n(B_{10}).
 \end{split}}
\end{equation}
The square induces comparison maps ${\cM}_{D_n}(B_{00})\to P_n$ and
$L_n(B_{00})\to Q_n$. We will show that the latter comparison
becomes an equivalence after taking the limit over $n$.

We first establish the relative lifting property before scalar
localization. {Fix $n$ in $\nat$ and put $m:=n+1$, so that $\alpha_{nm}$ is semiprojective.} Let $L\to K$ be a map of finite anima and consider a
commutative square
\begin{equation}\label{eq:relative-mapping-lifting-problem}
{\xymatrix@C=3em@R=3em{
 L\ar[r]\ar[d]&{\cM}_{D_m}(B_{00})\ar[d]\\
 K\ar[r]&P_m.
 }}
\end{equation}
After restriction along $\alpha_{nm}$, we seek the dotted map in
\begin{equation}\label{eq:relative-mapping-shifted-lift}
{\xymatrix@C=2em@R=3em{
 L\ar[r]\ar[d]
   &{\cM}_{D_m}(B_{00})\ar[r]\ar[d]
   &{\cM}_{D_n}(B_{00})\ar[d]\\
 K\ar[r]\ar@{.>}[urr]
   &P_m\ar[r]&P_n.
 }}
\end{equation}
The lift includes the coherent commutativity data for the two
triangles, compatible with the given outer square. 

To apply Lemma~\ref{lem:general-relative-lifting}, we now choose a
model for this lifting problem. The mapping anima
${\cM}_{D_m}(B_{ij})$ are modelled by the singular complexes of the
topological spaces of homomorphisms
$D_m\to B_{ij}\otimes\bbK$
\cite[Remark~3.6]{budu}. Use a path-space model for their pullback
$P_m$. By \cite[Remark~3.10]{budu}, the square \eqref{eq:relative-mapping-lifting-problem} can be represented
by a strictly commutative diagram of simplicial sets with a finite
simplicial inclusion at the left. Denote its realized finite CW
pair by $(K,L)$.

In this model {the map $K\to P_m$ in \eqref{eq:relative-mapping-lifting-problem}} is represented by homomorphisms
\[
 \begin{split}
 \rho&:D_m\longrightarrow C(K,B_{01}\otimes\bbK),\\
 \lambda&:D_m\longrightarrow C(K,B_{10}\otimes\bbK),\\
 h&:D_m\longrightarrow C(K\times[0,1],B_{11}\otimes\bbK),
 \end{split}
 \qquad
 h|_{K\times\{0\}}=q\rho,\quad
 h|_{K\times\{1\}}=f\lambda.
\]
Let $p:B_{00}\to B_{01}$ be the top horizontal map of
the coefficient square \eqref{gwregwerfwrf}. The map $L\to {\cM}_{D_m}(B_{00})$ in
\eqref{eq:relative-mapping-lifting-problem} is represented by
\[
 \sigma:D_m\longrightarrow C(L,B_{00}\otimes\bbK).
\]
Put $Y:=K\times[0,1]$ and
$Z:=(K\times\{0\})\cup(L\times[0,1])$.
Define $\psi:D_m\to C(Z,B_{01}\otimes\bbK)$ by
\[
 \begin{aligned}
 \psi(a)(x,0)&:=\rho(a)(x)
     &&(\forall x\in K),\\
 \psi(a)(x,u)&:=(p\otimes\id_{\bbK})(\sigma(a)(x))
     &&(\forall (x,u)\in L\times[0,1])
 \end{aligned}
 \qquad(\forall a\in D_m).
\]
The formulas agree on $L\times\{0\}$ because the chosen
model of \eqref{eq:relative-mapping-lifting-problem} commutes
strictly. The same commutativity gives
\[
 h|_Z=(q\otimes\id_{\bbK})\psi.
\]
Thus the lifting problem of
Lemma~\ref{lem:general-relative-lifting} is
\begin{equation}\label{eq:relative-coefficient-lifting-problem}
 \xymatrix@C=1.5em@R=3em{
 D_n\ar[r]^{\alpha_{nm}}\ar@{..>}[d]_{\widetilde h}
   &D_m\ar[d]^{(h,\psi)}\\
 C(Y,B_{01}\otimes\bbK)\ar@{->>}[r]
   &C(Y,B_{11}\otimes\bbK)
      \times_{C(Z,B_{11}\otimes\bbK)}C(Z,B_{01}\otimes\bbK).
 }
\end{equation}
The bottom map sends $g$ to $((q\otimes\id_{\bbK})g,g|_Z)$.
Lemma~\ref{lem:general-relative-lifting} gives $\widetilde h$ with
\[
 \widetilde h|_Z=\psi\alpha_{nm},\qquad
 (q\otimes\id_{\bbK})\widetilde h=h\alpha_{nm}.
\]
Write $\widetilde h_u$ and $h_u$ for evaluation at time $u$.
At time $1$ we have
\[
 (q\otimes\id_{\bbK})\widetilde h_1
 = h_1\alpha_{nm}
 = (f\otimes\id_{\bbK})\lambda\alpha_{nm}.
\]
Thus $(\widetilde h_1,\lambda\alpha_{nm})$ defines a map
$D_n\to C(K,B_{00}\otimes\bbK)$ and hence a lift
$K\to {\cM}_{D_n}(B_{00})$. In the path-space model for $P_n$,
the comparison homotopy is given at time $u$ by
\begin{equation}\label{eq:exactness-path-comparison-homotopy}
{ \bigl(\widetilde h_u,\lambda\alpha_{nm},
    (h_{u+(1-u)v}\alpha_{nm})_{v\in[0,1]}\bigr).}
\end{equation}
{At $u=0$, \eqref{eq:exactness-path-comparison-homotopy} is the restricted map} $K\to P_m\to P_n$.
At $u=1$ it is the image of the lift, with a constant quotient
path. The homotopy is constant on $L$ by the boundary condition
on $Z$, giving the relative lift in
\eqref{eq:relative-mapping-shifted-lift}.

{By Lemma~\ref{lem:general-scalar-translations}, each scalar
localization is a filtered colimit of finite localizations
indexed by $\operatorname{Fin}(S_n)$. For a fixed  
$T$ in $\operatorname{Fin}(S_n)$, the scalar morphism $\tau_T$ gives a single
telescope for each of the four coefficient algebras. Scalar
extension is compatible with  morphisms between coefficient algebras, so for fixed $T$   
these telescope models also compute the corresponding pullback \begin{equation}\label{eq:scalar-localized-pullback-telescope}
{L_{n,T}(B_{01})\times_{L_{n,T}(B_{11})}L_{n,T}(B_{10})
\simeq\colim\bigl(P_n\xrightarrow{-\oplus\tau_T}P_n
\xrightarrow{-\oplus\tau_T}\cdots\bigr).}
\end{equation}
Here we use that
filtered colimits commute with finite limits.  
 {The translation on $P_n$ in \eqref{eq:scalar-localized-pullback-telescope} adds the scalar extensions of $\tau_T$ at the three coefficient algebras, with the constant comparison path.}
Applying this first to the telescopes and then to the finite-localization
diagrams computes $Q_n$ from the analogous telescopes on $P_n$.}

{{More explicitly, fix $n$ in $\nat$.}
Consider a finite relative lifting problem
\[
 \xymatrix@C=3em@R=3em{
 L\ar[r]\ar[d]&{L_{n+1}(B_{00})}\ar[d]\\
 K\ar[r]&{Q_{n+1}}.
 }
\]
{Mapping out of the finite anima $K$ and $L$ preserves filtered colimits. The space of commutative squares with left edge $L\to K$ is a finite limit of these mapping anima, so it also preserves filtered colimits in the right edge. Thus both maps and the specified commutativity homotopy are represented at one finite localization, for a finite subset $T$ of $S_{n+1}$, and then at a common stage $k$ of its telescope models.} 
Thus the
problem is represented by a commutative square
\[
 \xymatrix@C=3em@R=3em{
 L\ar[r]\ar[d]&{\cM}_{D_{{n+1}}}(B_{00})\ar[d]\\
 K\ar[r]&{P_{n+1}},
 }
\]
including its commutativity homotopy.}

{The map from stage $k$ of the telescope to the finite
localization sends $x$ to $x-k\tau_{T,B}$. 
After restriction along ${\alpha_{n,n+1}}$,
the preceding argument supplies the dotted lift in
\eqref{eq:relative-mapping-shifted-lift}. Compose this lift
with ${\cM}_{D_n}(B_{00})\to L_n(B_{00})$ and translate by
$-k(\tau_T{\alpha_{n,n+1}})_{B_{00}}$. {Scalar extension commutes with all morphisms in the coefficient square, so these translations commute with the comparison maps and induce a translation on $Q_n$. Applying them also to the comparison homotopy recovers the original localized lifting problem after restriction.} The resulting map
$K\to L_n(B_{00})$ gives the dotted lift in the original
localized lifting problem after restriction:
\begin{equation}\label{eq:localized-relative-shifted-lift}
 \xymatrix@C=2em@R=3em{
 L\ar[r]\ar[d]
   &{L_{n+1}(B_{00})}\ar[r]\ar[d]
   &L_n(B_{00})\ar[d]\\
 K\ar[r]\ar@{..>}[urr]
   &{Q_{n+1}}\ar[r]&Q_n.
 }
\end{equation}
}

Thus
\[
 (L_n(B_{00})\longrightarrow Q_n)_n
\]
has the shifted lifting property. By \cite[Theorem~2.30]{budu}, its
inverse limit is an equivalence.

\emph{Step 2. Passage to group completions.}
{By Step~1, the first map in \eqref{eq:exactness-scalar-localized-limit}}
\begin{equation}\label{eq:exactness-scalar-localized-limit}
{ \lim_n L_n(B_{00})\xrightarrow{\simeq}\lim_n Q_n
 \simeq
 (\lim_n L_n(B_{01}))\times_{\lim_n L_n(B_{11})}
 (\lim_n L_n(B_{10}))}
\end{equation}
is an equivalence. The second equivalence follows from the
{definition of $Q_n$ in \eqref{eq:exactness-comparison-pullbacks} and commutation of limits with pullbacks.}
By Lemma~\ref{lem:general-scalar-localization}, the maps
$L_n(B_{ij})\to {\cM}_{D_n}(B_{ij})^{\group}$ induce equivalences on
inverse limits, naturally in the coefficient algebra. Replacing
each of the four limits of $L_n$ by the corresponding limit of
group completions gives
\begin{equation}\label{eq:exactness-group-completed-pullback}
{ \lim_nM_{D_n}(B_{00})^{\group}
 \simeq \left(\lim_nM_{D_n}(B_{01})^{\group}\right)
 \times_{\lim_nM_{D_n}(B_{11})^{\group}}
 \left(\lim_nM_{D_n}(B_{10})^{\group}\right).}
\end{equation}

\emph{Step 3. Passage to ${\cR}_D$.}
For a cartesian square $(A_{ij})$ with surjective right vertical
map, stabilization followed by unitalization gives a unital
cartesian square $B_{ij}:=(A_{ij}\otimes\bbK)^{u}$, again with
surjective right vertical map. {Apply \eqref{eq:exactness-group-completed-pullback} and take fibres of the}
compatible scalar quotient maps to $\lim_nM_{D_n}(\C)^{\group}$. {By \eqref{eq:exactness-resolved-fibre} and commutation of limits,}
\[
 {\cR}_D(A_{00})\simeq
 {\cR}_D(A_{01})\times_{{\cR}_D(A_{11})}{\cR}_D(A_{10}).
\]
%
\end{proof}

\section{Comparison of component groups}
\label{sec:general-comparison}

\subsection{Reduction to connected components}
\label{sec:comparison-components}

{Throughout this subsection, let $D$ be in $\nCalg_{\sepa}$ and assume that ${\cR}_D$ is homological. This holds, in particular, if $D\in\mathcal U_{\mathrm{nq}}$ by Theorem~\ref{thm:intro-exactness}. 
We use the factorization ${\bar{\cR}}_D$ described in Corollary \ref{cor:R-E-factorization}.
  {Proposition~\ref{prop:intro-splitting} gives the decomposition}}
\[
{\cR}_D\simeq\Omega^\infty E(D,-)\oplus{\widehat{\cR}}_D.
\]
The following criterion reduces the
vanishing of ${\widehat{\cR}}_D$ to injectivity on connected
components for unital coefficients. 
It generalizes the first
step in \cite[Lemma~10.17]{orbita}.

\begin{lem}\label{lem:general-comparison-pi-zero}
 If 
$
 \pi_0\alpha_{D,C}:\pi_0{\cR}_D(C)\longrightarrow E_0(D,C)
$
is injective for every unital $C$ in $\Calg_{\sepa}$, then {$D$ is $\alpha$-admissible,  $\beta_D$ is an inverse to $\alpha_D$, and ${\widehat{\cR}}_D\simeq0$.}
\end{lem}
\begin{proof}
 We first pass from unital to arbitrary coefficients $B$ in $\nCalg_{\sepa}$. The  
unitalization sequence
\[
 0\longrightarrow B\xrightarrow{\iota_B}B^{u}
   \xrightarrow{{\epsilon_B}}\C\longrightarrow0
\]
is split.
We consider the commutative square
\[
 \xymatrix@C=3em@R=2.5em{
 \pi_0{\cR}_D(B)\ar[r]^{\pi_0\alpha_{D,B}}\ar[d]^{\pi_{0}{\cR}_{D}(\iota_{B})}
   &E_0(D,B)\ar[d]\\
 \pi_0{\cR}_D(B^{u})\ar[r]^{\pi_0\alpha_{D,B^{u}}}
   &E_0(D,B^{u}).
 }
\]
{The left vertical map is injective} by exactness of ${\cR}_{D}$, and the lower horizontal map is injective by assumption.
Hence the upper horizontal map is injective.

%
Together with  $\alpha_D\beta_D\simeq\id$  from {Proposition~\ref{prop:intro-splitting}}   giving surjectivity of $\alpha_{D}$  this proves that $\pi_0\alpha_{D}$
is an isomorphism.
Using exactness and homotopy invariance of the domain and target of $\alpha_{D}$ and applying it to suspensions one concludes that $\pi_{j}\alpha_{D}$ is an isomorphism for every $j$ in $\nat$.
 
 Consequently $\alpha_{D}$ is an equivalence with inverse $\beta_{D}$ and ${\widehat{\cR}}_{D}\simeq 0$.
\end{proof}

The following result provides a tool to verify the assumption of Lemma \ref{lem:general-comparison-pi-zero} without the need to calculate $\alpha_{D}$. It is motivated by 
 the generator argument
in \cite[Lemma~10.17]{orbita} and uses the factorization ${\bar{\cR}}_{D}$ from Corollary \ref{cor:R-E-factorization}. For  $C$ in $\Calg_{\sepa}$, consider the assertion
\begin{equation}\label{eq:general-comparison-generation}
 \pi_0{\cR}_D(C)=
 \left\langle
   (\pi_0{\bar{\cR}}_D)(f)([u])\;\middle|\;
   f\in E_{0}(D,C)
 \right\rangle.
\end{equation}
Here $[u]\in\pi_0{\cR}_D(D)$ is the universal class  \eqref{grewfrewgg} and $\langle \cdots \rangle$ means the subgroup generated by the  list of elements indicated by $\cdots$. 

\begin{lem}\label{lem:general-comparison-generation}
If \eqref{eq:general-comparison-generation} holds for  every unital  $C$ in $\Calg_{\sepa}$, then the assumption of Lemma \ref{lem:general-comparison-pi-zero} is satisfied. \end{lem}
\begin{proof}
Any natural endomorphism $t$ of  the functor $\pi_0{\bar{\cR}}_D$ which fixes $[u]$
{fixes every generator on the right-hand side of} \eqref{eq:general-comparison-generation}.  
For such an endomorphism the assumption therefore implies that  $t_{e_{\sepa}(C)}=\id$ for every $C$ in $\Calg_{\sepa}$.

We apply this to the transformation  $t:=\pi_0(\bar\beta_D\bar\alpha_D)$. The construction of
$\bar\beta_D$ and \eqref{eq:general-universal-class} give
\[
 t_{e_{\sepa}(D)}([u])
 =\pi_0\bar\beta_{D,e_{\sepa}(D)}([\id_{e_{\sepa}(D)}])=[u].
\]
It follows that $\pi_0(\beta_{D,C}\alpha_{D,C})=\id$ and  $\pi_0\alpha_{D,C}$ is injective for all  $C$ in $\Calg_{\sepa}$.
\end{proof}
{{The following technical result will be used} in the proof of Proposition \ref{prop:rho-homological-alpha} below. Let $C$ and $D$ be in $\nCalg_{\sepa}$ and choose a rank-one projection $p$ in $\bbK$. We consider the homomorphisms}
\begin{equation}\label{eq:reduction-rotation-maps}
{\begin{aligned}
f,g:(C\otimes\bbK)^u&\longrightarrow((C\otimes\bbK)^u\otimes\bbK)^u,\\
f(c\otimes k,z)&:=((c\otimes k,z)\otimes p,0),\\
g(c\otimes k,z)&:=((c\otimes p,0)\otimes k,z).
\end{aligned}}
\end{equation}
{Using  the projection  $P:=(1_{(C\otimes\bbK)^u}\otimes p,0)$ in ${((C\otimes\bbK)^u\otimes\bbK)^u}$ we  define the homomorphism}
\begin{equation}\label{eq:reduction-rotation-scalar-correction}
{h:\C\longrightarrow((C\otimes\bbK)^u\otimes\bbK)^u,\qquad h(z):=z(1-P).}
\end{equation}
Let
\begin{equation}\label{eq:reduction-rotation-fibre-inclusion}
{i:\Fib_0\Bigl({\cM}_D((C\otimes\bbK)^u)^{\group}
\xrightarrow{{\cM}_D(\epsilon_{C\otimes\bbK})^{\group}}
{\cM}_D(\C)^{\group}\Bigr)\longrightarrow {\cM}_D((C\otimes\bbK)^u)^{\group}}
\end{equation}
be the canonical map.
\begin{lem}
\label{lem:reduction-rotation}
{There is an equivalence of morphisms of commutative monoids in anima}
\begin{equation}\label{eq:reduction-rotation-monoid-identity}
{{\cM}_D(g)\simeq {\cM}_D(f)+{\cM}_D(h)\circ {\cM}_D(\epsilon_{C\otimes\bbK}).}
\end{equation}
{Consequently, ${\cM}_{D}(f)^{\group}\circ i\simeq {\cM}_{D}(g)^{\group}\circ i$.}
\end{lem}
\begin{proof}
{Define
\begin{equation}\label{eq:reduction-rotation-unital-extension}
{\widetilde f:(C\otimes\bbK)^u\longrightarrow((C\otimes\bbK)^u\otimes\bbK)^u,\qquad
\widetilde f(c\otimes k,z):=((c\otimes k,0)\otimes p,z).}
\end{equation}
Then $\widetilde f=f+h\epsilon_{C\otimes\bbK}$, with orthogonal summands. A rotation between the two corner embeddings $k\mapsto k\otimes p$ and $k\mapsto p\otimes k$, followed by unital extension, gives $\widetilde f\simeq g$. Hence \eqref{eq:reduction-rotation-monoid-identity} follows from additivity. 
The second assertion follows  since ${\cM}_D(\epsilon_{C\otimes\bbK})^{\group}\circ i\simeq 0$.}
\end{proof}

Let $D$ be in $\nCalg_{\sepa}$ and recall that we assume that  ${\cR}_D$ is homological.
\begin{prop}\label{prop:rho-homological-alpha}
{If  $D$ is $\rho$-admissible, then $D$ is $\alpha$-admissible.}
\end{prop}
\begin{proof}
 By Lemmas~\ref{lem:general-comparison-generation} and~\ref{lem:general-comparison-pi-zero}, it suffices to establish \eqref{eq:general-comparison-generation} for every unital $C$ in $\Calg_{\sepa}$.
 

{Let $A$ be in $\nCalg_{\sepa}$. We assign to a homomorphism $\varphi:D\to A\otimes\bbK$ the element}
\[
{(\pi_0{\bar{\cR}}_D)(e_{\sepa}(\varphi))([u])\in\pi_0{\cR}_D(A).}
\]
This assignment is homotopy invariant and respects block sums of homomorphisms and
therefore extends to a homomorphism
\begin{equation}\label{eq:rho-generation-kappa}
{\kappa_A:(\pi_0{\cM}_D(A))^{\group}\longrightarrow\pi_0{\cR}_D(A).}
\end{equation}
{Choose a rank-one projection $p$ in $\bbK$ and define}
\begin{equation}\label{eq:rho-generation-ideal-corner}
{j:=\iota_{C\otimes\bbK}\circ(\id_C\otimes p):C\longrightarrow(C\otimes\bbK)^u.}
\end{equation}
By $\rho$-admissibility, we identify ${\cR}_D$   with   $\cT_D$  via \eqref{eq:rho-admissible-comparison} and hence $\pi_{0}{\cR}_{D}(C)$ with the kernel  of  
the right upper horizontal map  in 
\begin{equation}\label{eq:rho-generation-fibre-diagram}
{
 \xymatrix@C=3em@R=3em{
 \pi_0{\cR}_D(C)\ar@{^{(}->}[r]^{i}\ar[dr]^{{\pi_0{\cR}_D(j)}}\ar[d]^{{\pi_0{\cR}_D(\id_C\otimes p)}}_{\simeq}
  &(\pi_0{\cM}_D((C\otimes\bbK)^{u}))^{\group}
       \ar[d]^{{\kappa_{(C\otimes\bbK)^u}}}\ar[r]^-{\pi_{0}{\cM}_{D}(\epsilon_{C\otimes \bbK})^{\group}}& (\pi_{0}{\cM}_{D}(\C))^{\group}\\
 \pi_0{\cR}_D(C\otimes\bbK)
       \ar@{^{(}->}[r]_-{\pi_0{\cR}_D(\iota_{C\otimes\bbK})}
  &\pi_0{\cR}_D((C\otimes\bbK)^{u})\ar@{^{(}->}[r]&\pi_0 {\cM}_D(((C\otimes\bbK)^u\otimes\bbK)^u)^{\group} 
 }.
}
\end{equation}
We justify the {upper left triangle} using Lemma~\ref{lem:reduction-rotation}. By the definition of $[u]$ in \eqref{grewfrewgg} and the identification by $\rho$, the composite through $\kappa_{(C\otimes\bbK)^u}$ is induced by $f$ in \eqref{eq:reduction-rotation-maps}, whereas $\cT_D(j)$ is induced by $g$. {By the last assertion of the lemma, these composites agree after composition with the lower right horizontal map in \eqref{eq:rho-generation-fibre-diagram}. This map is injective because the scalar projection $\epsilon_{(C\otimes \bbK)^{u}\otimes \bbK}$ has a section. Thus the upper left triangle commutes.} 
{The upper left horizontal map is injective because the fibre sequence defining $\cT_D(C)$ splits. The lower left horizontal map is injective because exactness of ${\cR}_D$ sends the split unitalization sequence to a split fibre sequence. The lower left triangle commutes by functoriality of $\pi_0{\cR}_D$ and the definition of $j$ in \eqref{eq:rho-generation-ideal-corner}.}

{Every class $x$ in $\pi_0{\cR}_D(C)$} 
 is a class in $\pi_0({\cM}_D((C\otimes\bbK)^{u})^{\group})$ with trivial 
image under $\pi_{0}{\cM}_{D}(\epsilon_{C\otimes \bbK})^{\group}$. {Using $e_{\mathrm{add}}$ from \eqref{eq:additive-E-theory-factorization}, we have}
{a unique $y$ in $E_0(D,C)$ such that}
\begin{equation}\label{eq:rho-generation-scalar-kernel-class}
{e_{\add}(x)=e_{\sepa}(j)\circ y.}
\end{equation}
{Indeed, the splitting gives}
\begin{equation}\label{eq:rho-generation-E-splitting}
{E_0(D,(C\otimes\bbK)^u)\cong E_0(D,C)\oplus E_0(D,\C),}
\end{equation}
{with first summand included by postcomposition with $e_{\sepa}(j)$, and $e_{\add}(x)$ has coordinates $(y,0)$. The commutative triangle and \eqref{eq:rho-generation-scalar-kernel-class} give}
\begin{equation}\label{eq:rho-generation-cancellation-calculation}
{\pi_0{\cR}_D(j)(x)=\kappa_{(C\otimes\bbK)^u}(x)=(\pi_0{\bar{\cR}}_D)(e_{\add}(x))([u])=\pi_0{\cR}_D(j)\bigl((\pi_0{\bar{\cR}}_D)(y)([u])\bigr).}
\end{equation}
%
%
%
%
{Since $\pi_0{\cR}_D(j)$ is the composite of the corner equivalence and the lower left injection, it is injective. Thus \eqref{eq:rho-generation-cancellation-calculation} implies $x=(\pi_0{\bar{\cR}}_D)(y)([u])$.}
{Since $x$ was arbitrary} this establishes \eqref{eq:general-comparison-generation}.
\end{proof}

\begin{proof}[Proof of Proposition~\ref{prop:intro-semiprojective}]
 Note that $\tau$ is obtained from $\eta$ by applying stabilized reduction in the second variable. 
Therefore, if $\eta_{D}$ is an equivalence for unital coefficients, then $\tau_{D}$ is an equivalence for all coefficients. 
Hence $\eta$-admissibility implies {$\tau$-admissibility} which  is the 
  second containment.
  
  For the first containment, closure of the $\eta$-admissible
objects under equivalences, finite products, and retracts reduces the
proof to nuclear, quasidiagonal, semiprojective $D$. In this case, we first show that
$\tau_D$ is an equivalence and then identify it with $\eta_{D,B}$ for
unital $B$.

{By Theorem~\ref{thm:intro-exactness}, ${\cR}_D$ is homological. Since $D$ is semiprojective, its constant shape system shows that $D$ is $\rho$-admissible. Proposition~\ref{prop:rho-homological-alpha} therefore implies that $D$ is $\alpha$-admissible. By \eqref{eq:tau-alpha-rho}, $D$ is also $\tau$-admissible.}  
 
 Now we deduce $\eta$-admissibility from $\tau$-admissibility for nuclear, quasidiagonal, semiprojective $D$.
Let $B$ be in $\Calg_{\sepa}$ and put $B_n:=B\otimes M_n(\C)$,
with the upper-left corner inclusions as transition maps. Write
$\epsilon_n:B_n^u\to\C$ for the scalar quotient. 
Since  ${\cM}_D(-)^{\group}$  preserves finite products,
the canonical inclusion  $B_{n}\to B_{n}^{u}$ induces an equivalence
\begin{equation}\label{gerwwfrwerfrgwregre}
 {\cM}_D(B_n)^{\group}\xrightarrow{\simeq}
 \Fib_0\Bigl({\cM}_D(B_n^u)^{\group}
       \xrightarrow{{\cM}_D(\epsilon_n)^{\group}}{\cM}_D(\C)^{\group}\Bigr).\end{equation}
These equivalences commute with the transition maps by
naturality of the ideal inclusions and scalar quotients.

The colimit of the $B_n^u$ is $(B\otimes\bbK)^u$, and the
$\epsilon_n$ induce its scalar quotient
$\epsilon:(B\otimes\bbK)^u\to\C$. Since $D$ is semiprojective,
mapping out of $D$ in $\nCalg_{\sepa,h}$ sends sequential colimits
of separable $C^*$-algebras to colimits of mapping anima
\cite[proof of Corollary~3.79 and \textup{(3.43)}]{budu}.
Applying this to $B_n^u\otimes\bbK$ shows that
$\colim_n {\cM}_D(B_n^u)\to {\cM}_D((B\otimes\bbK)^u)$ is an equivalence.
Group completion preserves filtered colimits, and filtered
colimits commute with fibres in $\CGrp(\Ani)$. Taking the colimit
 of  the equivalences \eqref{gerwwfrwerfrgwregre} therefore gives
\begin{equation}\label{wregwerfwerfwerf}{\cM}_D(B)^{\group}\simeq\colim_n {\cM}_D(B_n)^{\group}
   \xrightarrow{\simeq}
 \Fib_0\Bigl({\cM}_D((B\otimes\bbK)^u)^{\group}
       \xrightarrow{{\cM}_D(\epsilon)^{\group}}{\cM}_D(\C)^{\group}\Bigr)
   \simeq\cT_D(B).
\end{equation} The first equivalence uses that   ${\cM}_{D}^{\group}$ inverts corner inclusions.
The composition of \eqref{wregwerfwerfwerf} with $\tau_{D,B}$ is $\eta_{D,B}$.
Since $\tau_{D,B}$ is an equivalence as shown above, also $\eta_{D,B}$ is an equivalence. 
Since $B$ was arbitrary in $\Calg_{\sepa}$, we conclude that 
 $D$ is {$\eta$-admissible}.
\end{proof}
%
%
\subsection{Topological component groups}
\label{sec:comparison-hausdorff}

{Let $F:\nCalg_{\sepa}\to\Ab$ be a functor. We equip each $F(B)$ with the condensed topology as follows.} Recall  that a light profinite space is an inverse limit of a sequence
of finite discrete spaces {\cite[Lecture~2]{clausen-scholze-analytic-stacks}}; see also \cite[Definition~4.1]{budu}.
For such a space $S$, any class $x$ in $F(C(S,B))$ gives rise to a map
$\hat x:S\to F(B)$ sending $s$  in $S$ to $F(\ev_{s})(x)$ in $F(B)$, where $\ev_{s}:C(S,B)\to B$ is the evaluation at $s$.
 \begin{ddd}\label{def:condensed-topology} {The condensed topology on $F(B)$ is the finest topology making $\hat x:S\to F(B)$ continuous for every light profinite space $S$ and every $x$ in $F(C(S,B))$.}
\end{ddd}
{Addition in $F(B)$ is separately continuous. However, it is not automatic from this definition that the condensed topology is a group topology: passing from condensed objects to their underlying topological spaces need not preserve products; see \cite[Remark~1.8]{scholze-condensed}.} 
{For  our main examples $\pi_0{\cR}_D(-)$ and $E_0(D,-)$, the comparison with the shape topology will show that the condensed topology is also a group topology; see Remark~\ref{rem:condensed-group-topologies}.}

Any map $B\to B'$ induces functorially maps $C(S,B)\to C(S,B')$ for all {light profinite spaces} $S$. This implies that $F$ sends morphisms in $\nCalg_{\sepa}$   to continuous maps.
Moreover, 
every 
 natural transformation between abelian group-valued functors
 has continuous components with respect to the condensed topologies.
 
This equips the functors $\pi_{0}{\cR}_{D}$ and $E_{0}(D,-)$ with the condensed topology. {If ${\cR}_D$ is homological, the section $\beta_D$ exists by Proposition~\ref{prop:intro-splitting}, and the transformations}
\[
\alpha_{D}:\pi_{0}{\cR}_{D}\to E_{0}(D,-)\ , \qquad  \beta_{D}:E_{0}(D,-)\to \pi_{0}{\cR}_{D}
\]%
are continuous.

\begin{rem}\label{kophertgertgerth} {We recall Pimsner's condition, note that it holds for $E$-theory, and give an example showing that it fails for general functors $F$ as above.}

{Let $\nat^{\dagger}$ be the one-point compactification of $\nat$. We say that $F$ satisfies Pimsner's condition if, for every $B$ in $\nCalg_{\sepa}$, a sequence $(x_n)_{n\in\nat}$ in $F(B)$ converges to $x$ in the condensed topology if and only if there is a class in $F(C(\nat^{\dagger},B))$ with evaluations $x_n$ at $n$ and $x$ at $\infty$. For every $D$ in $\nCalg_{\sepa}$, the functor $E_0(D,-)$ satisfies this condition by \cite[Theorem~A]{carrion-schafhauser-topology}.}

For any functor $F$ as above, a class in $F(C({\nat^{\dagger}},B))$
gives a convergent sequence by the definition of the condensed
topology.  {However,} the converse can fail.

Consider the functor to abelian groups given by
\[
 F(B):=\prod_{\nat}B/\bigoplus_{\nat}B,
\]
where the product and sums are  interpreted in $C^{*}$-algebras and we then take the underlying abelian group.
We show that the condensed topology on $F(B)$ is indiscrete.
For $[(b_j)_j]$ in $F(B)$,
define    $(f_j)_j$ in $\prod_{\nat} C({\nat^{\dagger}},B)$
by
\[
 f_j(n):=
 \begin{cases}
 b_j & n\leq j,\\
 0 & n>j,
 \end{cases}
 \qquad f_j(\infty):=0.
\]
Let $[(f_{j})_{j}]$ be its class in $F(C(\nat^{\dagger},B))$. Then
 $\ev_{n}([(f_{j})_{j}])=[(b_{j})_{j}]$
 since   representatives differ   in only finitely many entries.
 Furthermore $\ev_{\infty}([(f_{j})_{j}])=0$.
   Thus every constant sequence in $F(B)$  converges to
zero. Since addition in $F(B)$ is separately continuous, this shows that  the condensed topology on $F(B)$ is indiscrete.

For $B=\C$ we consider the class  $[(n)_{j}]$ in $F(B)$ of   the constant  family with value $n$.  Then $\lim_{n\to \infty} [(n)_{j}]=0$ in the condensed topology and
$\|[(n)_{j}]\|=n$
 in the canonical norm on $F(B)$.  
 For any class $y$ in $F(C({\nat^{\dagger}},\C))$ we have
 $\sup_{n} \|\ev_{n}(y)\|<\infty$.
 Consequently, there is no class in $F(C({\nat^{\dagger}},\C))$ witnessing  $\lim_{n\to \infty} [(n)_{j}]=0$.
  \hB
\end{rem}

Any group-valued functor $F:(\nCalg_{\sepa})^{\op}\to \Ab$ refines to a functor with values in topological abelian groups by equipping the values $F(D)$ with the shape topology
  
\begin{ddd}\label{def:intro-shape-topology}
The \emph{shape topology} on $F(D)$ is the coarsest  
topology making every restriction
\[
 F(u):F(D)\longrightarrow F(D')_{\mathrm{discrete}}
\]
along a semiprojective homomorphism $u:D'\to D$ continuous.
\end{ddd}
 
{The shape topology on $F(D)$ is a group topology with a neighbourhood basis of open subgroups.}
If $D\to C$ is a morphism in $\nCalg_{\sepa}$ and $D'\to D$ is semiprojective, then also $D'\to D\to C$
is semiprojective. This implies that $F$ sends morphisms in ${\nCalg_{\sepa}}$ to continuous maps. The components of  any natural transformation between contravariant functors are continuous.

\begin{rem}\label{okpzehhtre}If $F:(\nCalg_{\sepa})^{\op}\to \Ab$  is homotopy invariant and $(D_{n})$ is a
shape system for $D$ with structure maps $d_{n}:D_{n}\to D$ and $\alpha_{n,n+1}:D_{n}\to D_{n+1}$, then the subgroups $\ker F(d_n)$ form a
neighbourhood basis of zero for the shape topology on $F(D)$.

Indeed, $d_n=d_{n+1}\alpha_{n,n+1}$ is semiprojective,
since semiprojectivity is preserved by postcomposition.
Thus $\ker F(d_n)$ is open in the shape topology.
Conversely, a semiprojective homomorphism $u:D'\to D$ admits
a factorization up to homotopy
\[
 D'\xrightarrow{v}D_n\xrightarrow{d_n}D
\]
for some $n$ in $\nat$; see also the proof of
\cite[Proposition~3.78]{budu}. Homotopy invariance gives
\[
 F(u)=F(v)F(d_n),\qquad
 \ker F(d_n)\subseteq\ker F(u).
\]
Finally, for $m\geq n$, the identity $d_n=d_m\alpha_{nm}$
implies $\ker F(d_m)\subseteq\ker F(d_n)$. Hence every finite
intersection of the defining open kernels contains one of the
subgroups $\ker F(d_n)$, proving the assertion. 
  \hB
\end{rem}

 Using the {contravariant} functoriality in $D$  the values of   $\pi_{0}{\cR}_{D}$ and $E_{0}(D,-)$ also have the shape topology.  
\begin{prop} \label{kophertgertr}For all $D,B$ in $\nCalg_{\sepa}$
the condensed and shape topologies on $E_0(D,B)$ coincide.
  \end{prop}
\begin{proof} This is shown in
\cite[Propositions~4.25 and 4.42]{budu},
see also Carri\'on--Schafhauser
\cite[Theorem~A and Section~4.2]{carrion-schafhauser-topology}. \end{proof}

  Let $D$ and $B$ be in $\nCalg_{\sepa}$.
\begin{prop}\label{prop:R-stage-shape-topology}
{For any shape system $(D_{n})$ of $D$, the shape topology
of Definition~\ref{def:intro-shape-topology} on $\pi_0{\cR}_D(B)$
is the topology induced by
\[
 \pi_0{\cR}_D(B)\longrightarrow
 \lim_n\pi_0\cT_{D_n}(B)_{\mathrm{discrete}},
\]
where the target has the inverse limit topology.
}
\end{prop}

\begin{proof}
{Let $Q:=L\circ P$ be the colocalization endofunctor  of
$\Ind_{\aleph_1}(\nCalg_{\sepa,h})$ obtained from the adjunction   \eqref{eq:R-shape-localization}  
 and write $\epsilon:Q\to\id$ for
its counit. The functor $Q$ preserves countable filtered colimits since 
  $L$ is a left adjoint and $P$ preserves these colimits
by \cite[Theorem~3.49(3)]{budu}. 
}

For fixed $B$, let $\widetilde{\cT}_{-}(B):\Ind_{\aleph_1}(\nCalg_{\sepa,h})^{\op}\to \CGrp(\Ani)$ be the unique extension along $y$ of
$X\mapsto\cT_X(B)$  
sending countable filtered colimits to limits. {For $X$ in $\nCalg_{\sepa}$, we have} \begin{equation}\label{twertwertwerwert}
   \widetilde{\cT}_{{y_h(X)}}(B)\simeq\cT_X(B),
  \qquad
  \widetilde{\cT}_{Q({y_h(X)})}(B)\simeq {\cR}_X(B).
 \end{equation}


{Write
\[
 p_n:{\cR}_D(B)\longrightarrow\cT_{D_n}(B),\qquad
 r_n:{\cR}_D(B)\longrightarrow {\cR}_{D_n}(B)
\]
for the limit projection and the map induced by $d_n:D_n\to D$,
respectively. By Remark~\ref{okpzehhtre}, the kernels of
$\pi_0r_n$ form a neighbourhood basis for the shape topology.
We show that these kernels and the kernels of $\pi_0p_n$
are cofinal.
}

By \eqref{gerfwergwerfewrfwe},
we have an equivalence $\colim_m {y_h(D_m)}\simeq Q({y_h(D)})$. We
let $i_{n}:{y_h(D_n)}\to Q({y_h(D)})$ denote the structure map.
%
%
{The composite
\[
 Q({y_h(D_n)})\xrightarrow{\epsilon_{{y_h(D_n)}}}{y_h(D_n)}
 \xrightarrow{i_n} Q({y_h(D)})
\]
is $Q({y_h(d_n)})$. Indeed, both maps have the same composite
with $\epsilon_{{y_h(D)}}$, and the universal property of
colocalization identifies maps from the colocal object
$Q({y_h(D_n)})$ to $Q({y_h(D)})$ with maps to ${y_h(D)}$.
Applying $\widetilde{\cT}_{-}(B)$ and using \eqref{twertwertwerwert} gives a factorization
\[
 r_n:\quad {\cR}_D(B)\xrightarrow{p_n}\cT_{D_n}(B)
       \longrightarrow {\cR}_{D_n}(B).
\]
Consequently,
\[
 \ker(\pi_0p_n)\subseteq\ker(\pi_0r_n).
\]
}

 {For the reverse inclusion we use \[
 \begin{aligned}
 \colim_m Q({y_h(D_m)})
 &\simeq Q\left(\colim_m {y_h(D_m)}\right) 
 \simeq Q(Q({y_h(D)}))\simeq Q({y_h(D)}).
 \end{aligned}
\]
The object ${y_h(D_n)}$ is compact with respect to countable filtered
colimits in the ind-completion. Hence $i_n$  factors, for some $m\geq n$, as
\[
 {y_h(D_n)}\longrightarrow Q({y_h(D_m)})
       \xrightarrow{Q({y_h(d_m)})}Q({y_h(D)}).
\]
Applying $\widetilde{\cT}_{-}(B)$ now gives
\[
 p_n:\quad {\cR}_D(B)\xrightarrow{r_m}{\cR}_{D_m}(B)
       \longrightarrow\cT_{D_n}(B),
\]
and therefore
\[
 \ker(\pi_0r_m)\subseteq\ker(\pi_0p_n).
\]
Thus the two families of kernels define the same topology.
}
\end{proof}

 \begin{rem}\label{kopgrethertherthge}
 The standard Milnor sequence applied to the limit formula \eqref{eq:R-shape-formula}  for ${\cR}_{D}(B)$ gives
\begin{equation}\label{gwegerfrewfw}
 0\longrightarrow\lim\nolimits_{n\in \nat}^1\pi_1\cT_{D_n}(B)
 \xrightarrow{{\partial_B}}\pi_0{\cR}_D(B)
 \xrightarrow{r}\lim_{n\in \nat}\pi_0\cT_{D_n}(B)
 \longrightarrow0.
\end{equation}
{We use the quotient description}
\begin{equation}\label{eq:Milnor-lim-one-representatives}
{{\lim\nolimits_n^1}\pi_1\cT_{D_n}(B)
\cong\operatorname{coker}\left(1-\mathrm{shift}:\prod_n\pi_1\cT_{D_n}(B)\longrightarrow\prod_n\pi_1\cT_{D_n}(B)\right),}
\end{equation}
{where $\mathrm{shift}((b_n)_n):=(\alpha_{n,n+1}^*b_{n+1})_n$. A representative for an element in the ${\lim\nolimits^{1}}$-term is thus a sequence $(a_n)_n$ with $a_n$ in $\pi_1\cT_{D_n}(B)$. 
We write $[(a_n)_n]$ for its class. 
}

{By Proposition~\ref{prop:R-stage-shape-topology}, the shape topology on $\pi_0{\cR}_D(B)$ is induced by $r$ in \eqref{gwegerfrewfw}, with the prodiscrete topology on its target. The image of ${\partial_B}$ coincides with the closure of zero.}
 \hB
\end{rem}

%
%

\begin{prop} 
\label{prop:R-condensed-shape}
The condensed and shape topologies on $\pi_0{\cR}_D(B)$ coincide.
 \end{prop}
\begin{proof}
{We first show that the identity from the condensed to the shape topology is continuous. Let $S$ be a light profinite space and let $x$ be in $\pi_0{\cR}_D(C(S,B))$. Write $p_n:{\cR}_D(-)\to\cT_{D_n}(-)$ for the limit projections. By Proposition~\ref{prop:R-stage-shape-topology}, it suffices to show that
\begin{equation}\label{eq:R-condensed-stage-evaluation}
 S\xrightarrow{\hat x}\pi_0{\cR}_D(B)
   \xrightarrow{\pi_0p_n}\pi_0\cT_{D_n}(B)
\end{equation}
is locally constant for every $n$ in $\nat$.}

{Fix $n$ in $\nat$. 
The image of $(\pi_0p_{n+1})(x)$ under
 the  injection
\begin{equation}\label{eq:R-condensed-scalar-inclusion}
 \pi_0\cT_{D_{n+1}}(C(S,B))\lhook\joinrel\longrightarrow
 \pi_0\bigl({\cM}_{D_{n+1}}((C(S,B)\otimes\bbK)^u)^{\group}\bigr)
  \end{equation}
 is a difference $[\varphi]-[\psi]$ represented by homomorphisms
\begin{equation}\label{eq:R-condensed-stage-representatives}
 \varphi,\psi:D_{n+1}\longrightarrow
       (C(S,B)\otimes\bbK)^u\otimes\bbK.
\end{equation}}{For $s$ in $S$, let $\varphi_s$ and $\psi_s$ be the evaluations of \eqref{eq:R-condensed-stage-representatives}. By norm continuity and Lemma~\ref{lem:semiprojective-straightening}(2), for $s'$ in a sufficiently small neighbourhood of $s$ we have
\[
 \varphi_s\alpha_{n,n+1}\simeq\varphi_{s'}\alpha_{n,n+1},
 \qquad \psi_s\alpha_{n,n+1}\simeq\psi_{s'}\alpha_{n,n+1}.
\]
By \eqref{eq:R-condensed-scalar-inclusion} and compatibility of the limit projections, \eqref{eq:R-condensed-stage-evaluation} is therefore locally constant.}

We now show that the map
$\pi_{0}{\cR}_{D}(B)_{\mathrm{shape}}\to \pi_{0}{\cR}_{D}(B)_{\mathrm{cond}}  $ is continuous. 
{We have not yet shown that the condensed topology is a group topology. Nevertheless, by separate continuity of addition translations are homeomorphisms. Since the shape topology is first countable, it therefore suffices to show that sequences converging to zero in the shape topology also converge to zero in the condensed topology.} 
We use the profinite space $\nat^{\dagger}$ corepresenting converging sequences. {Let $(x_k)_{k\in\nat}$ be a sequence converging to zero in the shape topology.}
 {Put $x_{k,n}:=\operatorname{pr}_n(r(x_k))$ in $\pi_0\cT_{D_n}(B)$, where $r$ is defined in \eqref{gwegerfrewfw} and $\operatorname{pr}_n$ evaluates the limit at $n$. By Proposition~\ref{prop:R-stage-shape-topology}, for each fixed $n$ in $\nat$, the set $\{k\in \nat \mid x_{k,n}\not=0\}$ is finite.  } Let
\[
 \iota_k:B\longrightarrow C({\nat^{\dagger}},B)\ , \qquad b\mapsto \left(n\mapsto  \left\{\begin{array}{cc} b&n=k\\0 &else  \end{array} \right.\right)
\]
For every $n$ in $\nat $ the sum
\[
 y_n:=\sum_k(\iota_k)_*x_{k,n}
 \in\pi_0\cT_{D_n}(C({\nat^{\dagger}},B))
\]
{has finitely many nonzero terms, and $\alpha_{n,n+1}^{*}y_{n+1}=y_n$. By the Milnor sequence \eqref{gwegerfrewfw}, we choose a class $y$ in $\pi_0{\cR}_D(C(\nat^{\dagger},B))$ with $r(y)=(y_n)_n$. Let $p:\nat^\dagger\to *$ be the unique map, inducing}
\begin{equation}\label{eq:R-condensed-constant-inclusion}
{p^*:B\longrightarrow C(\nat^\dagger,B).}
\end{equation}
{Since $(\ev_\infty)_*y_n=0$ for every $n$, we have $r((\ev_\infty)_*y)=0$. Replacing $y$ by $y-(p^*)_*(\ev_\infty)_*y$ therefore preserves $r(y)=(y_n)_n$ and gives $(\ev_\infty)_*y=0$. Furthermore,}
$x_k-(\operatorname{ev}_k)_*y$   vanishes at every stage.
{The image of $\partial_B$ from \eqref{gwegerfrewfw}}
  is precisely the subgroup of classes vanishing at
every stage. Choose representatives with
\[
 {\partial_B([(a_{k,n})_n])}
   =x_k-(\operatorname{ev}_k)_*y,\qquad
 a_{k,n}=0\quad(n<k).
\]
The latter condition can always be arranged: changing finitely
many coordinates does not change a ${\lim\nolimits^1}$ class, because every
finitely supported sequence lies in the image of
$1-\mathrm{shift}$.

{For each $n$ in $\nat$, we consider the sum $\sum_{k\leq n}(\iota_k)_*a_{k,n}$ in $\pi_1\cT_{D_n}(C(\nat^{\dagger},B))$.}  {Apply the Milnor map in \eqref{gwegerfrewfw} with coefficient algebra}
$C({\nat^{\dagger}},B)$ to this sequence to get
\[
 z:=\partial_{C({\nat^{\dagger}},B)}
       {\left(\left[\left(\sum_{k\leq n}(\iota_k)_*a_{k,n}\right)_n\right]\right)} \in \pi_{0}{\cR}_{D}(C(\nat^{\dagger},B)). 
\]
Naturality gives
\[
 (\operatorname{ev}_k)_*z
   =x_k-(\operatorname{ev}_k)_*y,\qquad
 (\operatorname{ev}_\infty)_*z=0.
\]
The element $y+z$ in $ \pi_{0}{\cR}_{D}(C(\nat^{\dagger},B))  $  witnesses that the condensed limit of  the sequence $(x_{k})_{k}$ is zero.
%
\end{proof}

\begin{rem}\label{rem:condensed-group-topologies}
{Propositions~\ref{kophertgertr} and~\ref{prop:R-condensed-shape} identify the condensed and shape topologies on $E_0(D,B)$ and $\pi_0{\cR}_D(B)$, respectively. Since the shape topology is known to be a group topology 
these two groups are topological abelian groups also for their condensed topologies.}
\hB
\end{rem}

\begin{kor}\label{cor:comparison-shape-continuity}
{Assume that ${\cR}_D$ is homological. The transformations}
\[
{\alpha_D:\pi_0{\cR}_D(-)\longrightarrow E_0(D,-),\qquad
\beta_D:E_0(D,-)\longrightarrow\pi_0{\cR}_D(-)}
\]
{are continuous for the shape topologies.}
\end{kor}
\begin{proof}
{Naturality with respect to coefficient evaluations and Definition~\ref{def:condensed-topology} give continuity of $\alpha_D$ and $\beta_D$ for the condensed topologies. Propositions~\ref{kophertgertr} and~\ref{prop:R-condensed-shape} identify these with the shape topologies, proving the assertion.}
\end{proof}

We let $(-)_{\mathrm{Hd}}$ denote the operation of taking the Hausdorff quotient.

{Let $D$ and $B$  be in $\nCalg_{\sepa}$  and let $(D_n)$ be a shape system for $D$.}
\begin{kor}\label{cor:unital-hausdorff-formula}
 If $B$ is unital, then there is a canonical isomorphism of topological abelian groups
\begin{equation}\label{eq:unital-hausdorff-formula}
 (\pi_0{\cR}_D(B))_{\mathrm{Hd}}
 \cong\lim_n\pi_0\bigl({\cM}_{D_n}(B)^{\group}\bigr)_{\mathrm{discrete}}{.}
\end{equation}
\end{kor}
\begin{proof}
 {The shape-limit formula \eqref{eq:R-shape-formula}, the Milnor sequence \eqref{gwegerfrewfw}, and Proposition~\ref{prop:R-stage-shape-topology} give
\[
 (\pi_0{\cR}_D(B))_{\mathrm{Hd}}
 \cong\lim_n\pi_0\cT_{D_n}(B)
 \cong\lim_n\pi_0\bigl({\cM}_{D_n}(B)^{\group}\bigr).
\]
The second isomorphism is induced by Lemma~\ref{lem:group-scalar-reduction}.}
\end{proof}

\begin{kor}
{Assume that ${\cR}_D$ is homological.} We have induced transformations
\begin{equation}\label{hrtegrthetr}
 \alpha_{D,\mathrm{Hd}}:(\pi_{0}{\cR}_{D}(-))_{\mathrm{Hd}}\to E_{0}(D,-)_{\mathrm{Hd}} \ ,  \quad 
 \beta_{D,\mathrm{Hd}}:E_{0}(D,-)_{\mathrm{Hd}}\to (\pi_{0}{\cR}_{D}(-))_{\mathrm{Hd}}\end{equation}
 such that
 $\alpha_{D,\mathrm{Hd}}\beta_{D,\mathrm{Hd}}\simeq \id$.
\end{kor}
\begin{proof}
{By Corollary~\ref{cor:comparison-shape-continuity}, both transformations have continuous homomorphisms as components. These preserve the closure of zero and hence induce continuous homomorphisms on Hausdorff quotients. The identity follows from $\alpha_D\beta_D\simeq\id$ in Proposition~\ref{prop:intro-splitting}.}
\end{proof}

{The following criterion reduces the Hausdorff comparison to injectivity on inverse limits of stage groups. We verify this criterion in the proof of Theorem~\ref{thm:nuclear-qd-hausdorff-comparison}.}

\begin{kor}\label{prop:hausdorff-comparison-criterion}
{Assume that ${\cR}_D$ is homological.} If the comparison map
\begin{equation}\label{gregwergwefwerf} \lim_n\pi_0\cT_{D_n}(B)
   \longrightarrow\lim_n E_0(D_n,B)
\end{equation} is injective, then the {components} at $B$ of transformations in  \eqref{hrtegrthetr} are  isomorphisms of topological
 groups. \end{kor}
\begin{proof}
{By \eqref{eq:E-shape-limit-comparison} and the Milnor sequence, the map $E_0(D,B)\to\lim_n E_0(D_n,B)$ is surjective. Remark~\ref{okpzehhtre} identifies its kernel with the closure of zero for the shape topology. Thus it induces an isomorphism of topological groups from the Hausdorff quotient to the inverse limit of the discrete stage groups. Together with Proposition~\ref{prop:R-stage-shape-topology} and \eqref{gwegerfrewfw}, we obtain the commutative square}
\begin{equation}\label{eq:hausdorff-stage-comparison-square}
{\xymatrix@C=4em@R=3em{
 (\pi_0{\cR}_D(B))_{\mathrm{Hd}}
 \ar[r]^{\alpha_{D,B,\mathrm{Hd}}}\ar[d]_{\cong}
 &E_0(D,B)_{\mathrm{Hd}}\ar[d]^{\cong}\\
 \lim_n\pi_0\cT_{D_n}(B)
 \ar[r]_{\lim_n\pi_0\tau_{D_n,B}}
 &\lim_n E_0(D_n,B).
}}
\end{equation}
{The square commutes because $\alpha_D$ identifies with $\lim_n\tau_{D_n}$ under \eqref{eq:R-shape-formula} and \eqref{eq:E-shape-limit-comparison}. By the hypothesis, $\alpha_{D,B,\mathrm{Hd}}$ is injective. Since \eqref{hrtegrthetr} gives the continuous section $\beta_{D,B,\mathrm{Hd}}$, the two maps are inverse isomorphisms of topological groups.}
\end{proof}

  In particular, it suffices to establish a
pro-equivalence between the two systems of component groups.

\subsection{The finite-stage cancellation  theorem}


Let $C$ be in $\nCalg_{\sepa}$. The fold map gives a split
exact sequence
\[
 0\longrightarrow qC\longrightarrow C*C
   \xrightarrow{\mathrm{fold}}C\longrightarrow0.
\] defining the Cuntz construction $qC$ of $C$. Composing the inclusion with $(\id,0):C*C\to C$ provides the
structure map $\varepsilon_C:qC\to C$ of the Cuntz construction. 
The Cuntz picture of $KK$-theory asserts that  
\[
\Omega^{\infty }KK(-,-)\simeq {\Hom_{L_K\nCalg_{\sepa,h}}(q-,-)}:(\nCalg_{\sepa})^{\op} \times \nCalg_{\sepa}\to \CGrp(\Ani)\ ,
\]%
\cite[Corollary~11.14]{orbita}.
For every $B$ in $\nCalg_{\sepa}$, {$\varepsilon_C$} induces a
morphism of commutative monoids
\begin{equation}\label{gewerfwerfrewgfwer}
 \varepsilon_C^*:{\cM}_C(B)\longrightarrow
 \Hom_{L_K\nCalg_{\sepa,h}}(qC,B)
 \simeq\Omega^\infty KK(C,B).
\end{equation}
   Since the target of
\eqref{gewerfwerfrewgfwer} is a commutative group in anima,
the map extends uniquely through the group completion
\begin{equation}\label{hrtgertherthert}{\cM}_C(B) \longrightarrow {\cM}_C(B)^{\group}\longrightarrow\Omega^\infty KK(C,B).
\end{equation}
{Let $D$ be in $\nCalg_{\sepa}$ and choose a shape system $(D_n)$ for $D$. We write the induced map on components as}
\begin{equation}\label{eq:q-stage-comparison}
 c_n:\pi_0\bigl({\cM}_{D_n}(B)^{\group}\bigr)
       \longrightarrow KK_0(D_n,B).
\end{equation}
These maps are 
 compatible with restriction along the structure maps of
the shape system:
\[
{c_{n}\circ \pi_{0}M_{\alpha_{n,n+1}}(B)^{\group} =KK_{0}(\alpha_{n,n+1},B) \circ  c_{n+1}\ .}
\]

For a homomorphism $\tau:D_n\to M_s(\C)$ we use the notation $\tau_B$ for the scalar representation
associated to $\tau$ in \eqref{fqqwedweqfqwe}, viewing $M_s(\C)$
as a finite matrix corner of $\bbK$.

\begin{theorem}\label{thm:q-finite-stage-cancellation}
We assume that $D$ is nuclear and quasidiagonal.  Then
for every $n$ in $\nat$ there is $m$ in $\nat$ with $m\geq n+2$ such that for every $B$ in $\Calg_{\sepa}$ and 
 homomorphisms
\[
 \theta_0,\theta_1:D_m\longrightarrow B\otimes\bbK
\]
having the same class in $KK_0(D_m,B)$ there are $s$ in $\nat$ and a homomorphism $\tau:D_n\to M_s(\C)$ such that
\begin{equation}\label{eq:q-common-scalar-homotopy}
 \theta_0\alpha_{nm}\oplus\tau_B
 \simeq\theta_1\alpha_{nm}\oplus\tau_B.
\end{equation}
\end{theorem}

{The cancellation in Theorem~\ref{thm:q-finite-stage-cancellation} is analogous to the stable uniqueness results of Dadarlat--Eilers \cite{dadarlat-eilers-asymptotic} and Szab\'o \cite[Theorem~A]{szabo-uniqueness}, which relate equality of $KK$-classes to asymptotic unitary equivalence after absorption. Here we obtain a homotopy after restriction to an earlier stage, with a common summand coming from a finite-dimensional scalar representation.}

 The following constructions and results prepare the proof of this theorem.
 Choose once and for all a faithful essential unital representation
\[
 \pi:D^{u}\longrightarrow\mathcal L(H)
\]
of infinite multiplicity on a separable infinite-dimensional Hilbert
space. Since $D$ is quasidiagonal, it can be chosen with  a family of finite-rank
projections $(p_j)_{j\in \nat}$ converging strongly to $1$ and satisfying
$\|[p_j,\pi(a)]\|\to0$ as $j\to\infty$ for $a$ in $D^{u}$.
For  $B$ in $\Calg_{\sepa}$, put $H_B:=B\otimes H$ and
$\pi_B:=1_B\otimes\pi$. The operators $1_B\otimes p_j$ are compact
projections and form an approximate unit for $\mathcal K_B(H_B)$.

The next lemma is the place where nuclearity of $D$ is used. 

\begin{lem}\label{lem:q-finite-control-absorption}
Let $\ell$ be in $\nat$, let $F$ be a finite self-adjoint
subset of the unit ball of $D_\ell$, and let $\eta$ be in
$(0,\infty)$.
There is $m$ in $\nat$ with $m\geq\ell$ 
such that every $B$   in $\Calg_{\sepa}$ and
representation
\[
 \rho:D_m\longrightarrow\mathcal L_B(H_B)
\]
{admits an isometry $v$ in $\mathcal L_B(H_B)$ such that, with $P:=vv^*$, we have
\begin{equation}\label{eq:q-absorption-block-estimates}
 \|P\pi_B(d_\ell(a))P-v\rho(\alpha_{\ell m}(a))v^*\|<\eta,
 \quad \|(1-P)\pi_B(d_\ell(a))P\|<\eta
 \quad(\forall a\in F).
\end{equation}}
\end{lem}
\begin{proof}
{Choose $\kappa$ in $(0,\eta)$ with $3\kappa<\eta^2$.}
Lemma~\ref{lem:general-cp}, applied to
$F\cup\{a^*a:a\in F\}$ with accuracy $\kappa/2$, gives $m$ in $\nat$ with $m\geq\ell$ such that every $\rho$ admits a cpc map
\[
 \phi:D\longrightarrow\mathcal L_B(H_B)
\]
with
\begin{equation}\label{eq:q-cpc-extension-approximation}
 \|\phi(d_\ell(x))-\rho(\alpha_{\ell m}(x))\|<\kappa/2
 \qquad(\forall x\in F\cup\{a^*a:a\in F\}).
\end{equation}
The lemma is uniform in its target algebra, so it applies to
$\mathcal L_B(H_B)$ and the choice of $m$ is independent of $B$.

The map $\phi$ has the ucp extension
$\phi^{u}(a+\lambda1):=\phi(a)+\lambda1$ to $D^{u}$.
Since $D^{u}$ is nuclear, $\pi_B$ is unitally absorbing by
Kasparov's absorption theorem \cite[Theorem~6]{kasparov-absorption};
{see \cite[Theorem~2.1(2)]{thomsen-absorbing} for this formulation. The absorption and compression-defect argument is also used in \cite[Lemma~5.11]{willett}.} {We can therefore choose an isometry $v$ in $\mathcal L_B(H_B)$  such that}
\begin{equation}\label{eq:q-absorption-compression-approximation}
 \|v^*\pi_B(d_\ell(x))v-\phi(d_\ell(x))\|<\kappa/2
 \qquad(\forall x\in F\cup\{a^*a:a\in F\}).
\end{equation}
{Combining \eqref{eq:q-cpc-extension-approximation} and \eqref{eq:q-absorption-compression-approximation} gives
\begin{equation}\label{eq:q-absorption-combined-compression}
 \|v^*\pi_B(d_\ell(x))v-\rho(\alpha_{\ell m}(x))\|<\kappa
 \qquad(\forall x\in F\cup\{a^*a:a\in F\}).
\end{equation}
Set $P:=vv^*$. Since $\kappa<\eta$, \eqref{eq:q-absorption-combined-compression} gives the first estimate in \eqref{eq:q-absorption-block-estimates}.}

{For the off-diagonal terms, we use
\[
 \begin{aligned}
 \|(1-P)\pi_B(d_\ell(a))v\|^2
 &=\bigl\|v^*\pi_B(d_\ell(a^*a))v
       -(v^*\pi_B(d_\ell(a))v)^*(v^*\pi_B(d_\ell(a))v)\bigr\|\\
 &<3\kappa<\eta^2
 \qquad(\forall a\in F).
 \end{aligned}
\]
Indeed, \eqref{eq:q-absorption-combined-compression} changes the first term by less than $\kappa$ and the product by less than $2\kappa$ when we replace the compressions by $\rho\alpha_{\ell m}$; the resulting difference is zero by multiplicativity. Multiplying on the right by $v^*$ gives the second estimate in \eqref{eq:q-absorption-block-estimates}.}
\end{proof}
%
%

\begin{rem}\label{rem:stable-operator-homotopy}
{We use the following nonequivariant form of stable operator
homotopy \cite[Definition~2.5 and Theorem~2.9]{gabe-szabo}.
Let $B$ be in $\Calg_{\sepa}$, let $A$ be a separable $C^*$-algebra, and let
\[
 \varphi_0,\varphi_1:A\longrightarrow\mathcal L_B(H_B)
\]
be homomorphisms such that
$\varphi_0(a)-\varphi_1(a)\in\mathcal K_B(H_B)$ for every $a$ in $A$.
If their Cuntz-pair class vanishes in $KK_0(A,B)$, then there are
a homomorphism $\rho:A\to\mathcal L_B(H_B)$ and a
norm-continuous path
\[
 U:[0,1]\longrightarrow
 U\bigl(\mathcal L_B(H_B\oplus H_B)\bigr),\qquad U_0=1,
\]
such that, for every $a$ in $A$ and $t$ in $[0,1]$,
\[
 \begin{split}
 U_1(\varphi_0(a)\oplus\rho(a))U_1^*
    =\varphi_1(a)\oplus\rho(a),\qquad 
 [U_t,\varphi_0(a)\oplus\rho(a)]
    \in\mathcal K_B(H_B\oplus H_B).
 \end{split}
\]}
\hB
\end{rem}

\begin{proof}[Proof of Theorem~\ref{thm:q-finite-stage-cancellation}]
\emph{1. Choose the controls and the resolution index.}
Fix $n$ in $\nat$. We use the second assertion of
Lemma~\ref{lem:semiprojective-straightening} to choose a finite set
$F$ of contractions in $D_{n+1}$ and a tolerance $\delta$ in $(0,\infty)$.
Two homomorphisms that are $\delta$-close on $F$ become homotopic
after restriction to $D_n$.

{Choose $\varepsilon$ in $(0,\infty)$ with $2\varepsilon<\delta$. Lemma~\ref{lem:semiprojective-straightening}(1), applied to $\alpha_{n+1,n+2}$ with straightening accuracy $\varepsilon$ on $F$, supplies a finite self-adjoint subset $E$ of the unit ball of $D_{n+2}$ and a tolerance $\gamma$ in $(0,\infty)$ such that a cpc map with multiplicativity defects less than $\gamma$ on $E$ can be straightened within $\varepsilon$ on $F$. We enlarge $E$ to contain $\alpha_{n+1,n+2}(F)$.}

Choose $\eta$ in $(0,\infty)$ with $2\varepsilon+2\eta<\delta$.
We impose further upper bounds on $\eta$ in Step~4. Apply
Lemma~\ref{lem:q-finite-control-absorption} at stage $n+2$ on $E$
to obtain $m$.

\emph{2. Obtain a homotopy with an infinite common summand.}
Let $\theta_0$ and $\theta_1$ be as in the theorem. We regard them
as representations into $\mathcal L_B(H_B)$ that take values in the
compact operators. They form a Cuntz pair with class
\begin{equation}\label{eq:q-vanishing-cuntz-pair}
 [\theta_0,\theta_1]=[\theta_0]-[\theta_1]=0
 \quad\text{in }KK_0(D_m,B).
\end{equation}

{We apply Remark~\ref{rem:stable-operator-homotopy} to the Cuntz pair in}
\eqref{eq:q-vanishing-cuntz-pair} and obtain $\rho$ and $U_t$.
The formula
\[
 \sigma_t(a):=U_t(\theta_0(a)\oplus\rho(a))U_t^*
\]
defines a point-norm continuous path of representations with
endpoints $\theta_i\oplus\rho$. Along this path we have
\begin{equation}\label{eq:q-compact-background-path}
 \sigma_t(a)-(0\oplus\rho(a))
 \in\mathcal K_B(H_B\oplus H_B).
\end{equation}


\emph{3. Transport the homotopy to the quasidiagonal background.}
{We apply Lemma~\ref{lem:q-finite-control-absorption} to $\rho$ with $\ell=n+2$, $F=E$, and $\eta$ from Step~1. We obtain an isometry $v$ in $\mathcal L_B(H_B)$. Putting
\[
 P:=vv^*,\qquad Q:=1-P,\qquad W:=1\oplus v.
\] we have 
the compression and off-diagonal estimates as in    \eqref{eq:q-absorption-block-estimates}, with $\ell=n+2$ and $a$ in $E$.}
Define
\begin{equation}\label{eq:q-transported-path}
 \Phi_t(a):=W\sigma_t(\alpha_{n+2,m}(a))W^*
            +\bigl(0\oplus Q\pi_B(d_{n+2}(a))Q\bigr).
\end{equation}
The summands are cpc and have orthogonal supports, so
$\Phi_t$ is cpc. The first summand is multiplicative.  Therefore we have
\begin{equation}\label{eq:q-transported-defect}
 \Phi_t(ab)-\Phi_t(a)\Phi_t(b)
       =0\oplus Q\pi_B(d_{n+2}(a))P\pi_B(d_{n+2}(b))Q,
\end{equation}
{Its norm is at most $\eta^2$ for all $a,b$ in $E$, uniformly in $t$, by the off-diagonal estimate in \eqref{eq:q-absorption-block-estimates} and its adjoint, since $E$ is self-adjoint.}
To prepare for finite-rank compression, we separate the compact
perturbation from the error in replacing the common summand.
\begin{equation}\label{eq:q-compact-plus-small}
 \Phi_t(a)-\bigl(0\oplus\pi_B(d_{n+2}(a))\bigr)=K_t(a)+E(a),
\end{equation}
where
\[
 K_t(a):=W\Bigl(\sigma_t(\alpha_{n+2,m}(a))
              -\bigl(0\oplus\rho(\alpha_{n+2,m}(a))\bigr)\Bigr)W^*
\]
and
\[
 \begin{aligned}
 E(a):=0\oplus\bigl(&v\rho(\alpha_{n+2,m}(a))v^*
                  +Q\pi_B(d_{n+2}(a))Q-\pi_B(d_{n+2}(a))\bigr).
 \end{aligned}
\]
For each $a$, the path $t\mapsto K_t(a)$ is norm continuous and
takes values in the compact operators by
\eqref{eq:q-compact-background-path}. {The compression and off-diagonal estimates in \eqref{eq:q-absorption-block-estimates} give
\begin{equation}\label{eq:q-background-error-bound}
 \|E(a)\|\leq2\eta\qquad(\forall a\in E).
\end{equation}}
At the endpoints we obtain
\begin{equation}\label{eq:q-transported-endpoints}
 \Phi_i(a)=\bigl(\theta_i(\alpha_{n+2,m}(a))
                 \oplus\pi_B(d_{n+2}(a))\bigr)+E(a).
\end{equation}


\emph{4. Compress the entire path to finite corners.}
Set
\[
 e_j:=(1_B\otimes p_j)\oplus(1_B\otimes p_j),
 \qquad k_j:=\operatorname{rank}(p_j).
\]
We compress the entire path by setting
\[
 \Psi_j:D_{n+2}\longrightarrow C([0,1],M_{2k_j}(B)),
 \qquad \Psi_j(a)(t):=e_j\Phi_t(a)e_j.
\]
Quasidiagonality gives $[e_j,0\oplus\pi_B(d_{n+2}(a))]\to0$ as $j\to\infty$.
The projections $e_j$ form an approximate unit for the compact
operators. Since $t\mapsto K_t(a)$ has compact image in norm, we
obtain uniform convergence
\[
 \sup_t\|[e_j,K_t(a)]\|\longrightarrow0\qquad(j\to\infty).
\]
For $a$ in $E$, the remaining term satisfies
$\|[e_j,E(a)]\|\le2\|E(a)\|\le4\eta$. Thus
\begin{equation}\label{eq:q-compressed-commutator}
 \limsup_{j\to\infty}\sup_t\|[e_j,\Phi_t(a)]\|\le4\eta\qquad(\forall a\in E).
\end{equation}
The multiplicativity defect is
\[
 \begin{split}
 \Psi_j(ab)(t)-\Psi_j(a)(t)\Psi_j(b)(t)
 ={}e_j\bigl(\Phi_t(ab)-\Phi_t(a)\Phi_t(b)\bigr)e_j
    +e_j\Phi_t(a)(1-e_j)\Phi_t(b)e_j.
 \end{split}
\]
Equations~\eqref{eq:q-transported-defect} and
\eqref{eq:q-compressed-commutator} show that the limiting bound
for this defect tends to zero with $\eta$. {We choose $\eta$ small enough that this bound is strictly less than the tolerance $\gamma$ from Step~1 and $2\varepsilon+2\eta<\delta$. Then we  choose $j$ sufficiently large that the multiplicativity defect of $\Psi_j$ is less than $\gamma$ on $E$.}

We also compress the scalar representation.
\[
 b_j:D_{n+2}\longrightarrow M_{k_j}(\C),\qquad
 b_j(a):=p_j\pi(d_{n+2}(a))p_j.
\]
Its multiplicativity defect tends to zero on every finite set by
quasidiagonality. {Increasing $j$ if necessary, we therefore obtain
\begin{equation}\label{eq:q-compressed-multiplicativity-control}
 {\|\Psi_j(ab)-\Psi_j(a)\Psi_j(b)\|<\gamma,
 \quad \|b_j(ab)-b_j(a)b_j(b)\|<\gamma
 \quad(\forall a,b\in E).}
\end{equation}}

\emph{5. Straighten and attach the endpoint homotopies.}
{In view of the bounds  \eqref{eq:q-compressed-multiplicativity-control}   we can apply Lemma~\ref{lem:semiprojective-straightening}(1) with straightening accuracy $\varepsilon$ and the  tolerance $\gamma$  to get   homomorphisms}
\[
 h_j:D_{n+1}\longrightarrow C([0,1],M_{2k_j}(B)),
 \qquad
 \tau_j:D_{n+1}\longrightarrow M_{k_j}(\C),
\]
{such that
\begin{equation}\label{eq:q-straightening-approximation-estimates}
 \|h_j(x)-\Psi_j(\alpha_{n+1,n+2}(x))\|<\varepsilon,
 \quad \|\tau_j(x)-b_j(\alpha_{n+1,n+2}(x))\|<\varepsilon
 \quad(\forall x\in F).
\end{equation}}
We include the finite corners in the compact operators on the
original Hilbert modules.

{For $x$ in $F$ and $i$ in $\{0,1\}$, substitute $a=\alpha_{n+1,n+2}(x)$ in \eqref{eq:q-transported-endpoints} and compress by $e_j$. By the definitions of $\Psi_j$ and $b_j$, this gives
\begin{equation}\label{eq:q-compressed-endpoint-estimate}
 \begin{aligned}
 &\Bigl\|\Psi_j(\alpha_{n+1,n+2}(x))(i)-\bigl((1_B\otimes p_j)\theta_i(\alpha_{n+1,m}(x))(1_B\otimes p_j)
   \oplus(1_B\otimes b_j(\alpha_{n+1,n+2}(x)))\bigr)\Bigr\|\\
 &\qquad=\|e_jE(\alpha_{n+1,n+2}(x))e_j\|
 \leq\|E(\alpha_{n+1,n+2}(x))\|\leq2\eta.
 \end{aligned}
\end{equation}
{The last inequality follows from \eqref{eq:q-background-error-bound}, since $\alpha_{n+1,n+2}(F)\subseteq E$ by Step~1.}}

{Since $\theta_i$ takes values in the compact operators and $(1_B\otimes p_j)_j$ is an approximate unit, the compression errors
\begin{equation}\label{eq:q-endpoint-compression-errors}
 \zeta_j:=\max_{\substack{x\in F\\i\in\{0,1\}}}
 \bigl\|(1_B\otimes p_j)\theta_i(\alpha_{n+1,m}(x))(1_B\otimes p_j)
       -\theta_i(\alpha_{n+1,m}(x))\bigr\|
\end{equation}
satisfy $\lim_{j\to\infty}\zeta_j=0$. Combining \eqref{eq:q-straightening-approximation-estimates}, \eqref{eq:q-compressed-endpoint-estimate}, and \eqref{eq:q-endpoint-compression-errors} gives}
\[
 \bigl\|{h_j(x)(i)}-
       \bigl(\theta_i(\alpha_{n+1,m}(x))
                    \oplus(1_B\otimes\tau_j(x))\bigr)\bigr\|
 <2\varepsilon+2\eta+\zeta_j.
\]
Increase $j$ so that the right side is smaller than $\delta$.
The choice of endpoint control in Step~1 then gives homotopies
after restriction to $D_n$.

With $\tau:=\tau_j\alpha_{n,n+1}$, we concatenate the endpoint
homotopies with the restriction of $h_j$ to obtain
\[
 \theta_0\alpha_{nm}\oplus\tau_B
 \simeq h_j(0)\alpha_{n,n+1}
 \simeq h_j(1)\alpha_{n,n+1}
 \simeq\theta_1\alpha_{nm}\oplus\tau_B.
\]
This proves \eqref{eq:q-common-scalar-homotopy}.
\end{proof}

Recall the definition of $c_{n}$ in \eqref{eq:q-stage-comparison}.

Let $D$ be in $\nCalg_{\sepa}$. 
\begin{kor} \label{jhiortjohertge}If  $D$ is nuclear and quasidiagonal, then 
for every $n$  in $\nat$ there exists $m$ in $\nat$ with $m\ge n$ such that for every $B$ in $\Calg_{\sepa}$ we have 
\begin{equation}\label{eq:q-pro-kernel-zero}
 \alpha_{nm}^*(\ker c_m)=0
 \quad\text{in }\pi_0{\cM}_{D_n}(B)^{\group}.
\end{equation}
{In particular, for every $B$ in $\Calg_{\sepa}$ the comparison induces an injection
\begin{equation}\label{eq:q-unital-limit-injectivity}
 \lim_n c_n:\lim_n\pi_0{\cM}_{D_n}(B)^{\group}
 \lhook\joinrel\longrightarrow\lim_n KK_0(D_n,B).
\end{equation}}
\end{kor}
\begin{proof}
{Choose $m$ as in Theorem~\ref{thm:q-finite-stage-cancellation}. Every element of $\ker c_m$ is a difference $[\theta_0]-[\theta_1]$ of representation classes with equal images in $KK_0(D_m,B)$. By \eqref{eq:q-common-scalar-homotopy}, their restrictions to $D_n$ become homotopic after adding the same scalar representation. Group completion cancels this summand, proving \eqref{eq:q-pro-kernel-zero}. If $x=(x_n)_n$ lies in the kernel of \eqref{eq:q-unital-limit-injectivity}, then, for each $n$, choosing such an $m$ gives $x_n=\alpha_{nm}^*(x_m)=0$. Hence $x=0$.}
\end{proof}
We now extend from unital to general $B$.
\begin{kor}\label{ijwohterhetgertge}
If $D$ is nuclear and quasidiagonal, then for every $B$ in
$\nCalg_{\sepa}$ the reduced comparison induces an injection
\[
 \lim_n\pi_0\cT_{D_n}(B)
 \longrightarrow\lim_n KK_0(D_n,B).
\]
\end{kor}
\begin{proof}
The scalar quotient
$\epsilon_{B\otimes\bbK}:(B\otimes\bbK)^u\to\C$ has the section
$\lambda\mapsto\lambda1$. Applying ${\cM}_{D_n}$ and group completion
preserves this splitting.
The long exact sequence for the fibre defining $\cT_{D_n}(B)$ contains
\[
 \resizebox{\linewidth}{!}{$\displaystyle
 \pi_1\bigl({\cM}_{D_n}((B\otimes\bbK)^u)^{\group}\bigr)
 \to \pi_1\bigl({\cM}_{D_n}(\C)^{\group}\bigr)
 \xrightarrow{\partial}\pi_0\cT_{D_n}(B)
 \to \pi_0\bigl({\cM}_{D_n}((B\otimes\bbK)^u)^{\group}\bigr)
 \to \pi_0\bigl({\cM}_{D_n}(\C)^{\group}\bigr)
 \to 0.$}
\]
The scalar quotient induces a split surjection on $\pi_1$, so
$\partial=0$. Consequently,
\[
 \pi_0\cT_{D_n}(B)\cong
 \ker\left(
 \pi_0\bigl({\cM}_{D_n}((B\otimes\bbK)^u)^{\group}\bigr)
 \longrightarrow\pi_0\bigl({\cM}_{D_n}(\C)^{\group}\bigr)\right).
\]
{We use the map $j$ from \eqref{eq:rho-generation-ideal-corner}, with $C=B$. Taking limits gives the left vertical injection in the natural square, whose right vertical map is induced by $j$:}
\[
 \xymatrix@C=3em@R=3em{
 \lim_n\pi_0\cT_{D_n}(B)\ar[r]\ar@{^{(}->}[d]
  &\lim_n KK_0(D_n,B)\ar[d]^{{j_*}}\\
 \lim_n\pi_0\bigl({\cM}_{D_n}((B\otimes\bbK)^u)^{\group}\bigr)
       \ar@{^{(}->}[r]
  &\lim_n KK_0(D_n,(B\otimes\bbK)^u).
 }
\]
 Since $(B\otimes\bbK)^u$ is
unital, the lower horizontal
map {is injective} by Corollary~\ref{jhiortjohertge}.  We conclude that
  the upper horizontal map is injective.
\end{proof}

%
%
%

\subsection{Shape topologies and the Hausdorff comparison}
\label{sec:kk-e-topologies}

Recall that for a nuclear $D$ in $\nCalg_{\sepa}$ and any $B$ in $\nCalg_{\sepa}$ the natural comparison
\[
\KK(D,B)\to \EE(D,B)
\]%
is an equivalence \cite[Remark~9.18]{orbita}. {The classical comparison follows from Skandalis's half-exactness result for $K$-nuclear algebras \cite{skandalis-k-nuclearity} and Higson's comparison theorem \cite[Theorem~3.5]{higson-excision}; see also \cite[Theorem~25.6.3]{blackadar-kbook} for a textbook account.} 

%
%
%
{Let $D$ and $B$ be in $\nCalg_{\sepa}$.}
\begin{prop}\label{prop:kk-e-topologies}
{If $D$ is nuclear, the natural comparison induces an isomorphism of topological groups}
\[
 KK_0(D,B)\xrightarrow{\cong}E_0(D,B),
\]
where either side can be equipped with its condensed topology or
its shape topology. In particular, the condensed and
shape topologies on $KK_0(D,B)$ coincide.
\end{prop}

\begin{proof}
We first show that the group isomorphism
\begin{equation}\label{gwerfewrfwerferf}{KK_0(D,B)} \xrightarrow{\cong}E_0(D,B) 
\end{equation}
is a homeomorphism ${KK_0(D,B)}_{\mathrm{cond}}\xrightarrow{\cong}E_0(D,B)_{\mathrm{cond}}$.
Nuclearity of $D$ gives natural isomorphisms
\[
 KK_0(D,C(S,B))\xrightarrow{\cong}E_0(D,C(S,B))
\]
for all light profinite spaces $S$.
{These isomorphisms identify the generating maps in Definition~\ref{def:condensed-topology} under \eqref{gwerfewrfwerferf}. Hence \eqref{gwerfewrfwerferf} is a homeomorphism for the condensed topologies.}

{We next prove continuity from the condensed to the shape topology on $KK_0(D,B)$. Let $S$ be a light profinite space, let $x$ be in $KK_0(D,C(S,B))$, and let $u:C\to D$ be semiprojective. We must show that
\begin{equation}\label{eq:KK-condensed-shape-evaluation}
 S\longrightarrow KK_0(C,B),
 \qquad s\longmapsto u^*(\ev_s)_*(x)
\end{equation}
is locally constant.}

{Consider the coefficient inclusion
\begin{equation}\label{eq:KK-infinite-coefficient-inclusion}
 \iota:B\longrightarrow B^u\otimes\mathcal O_\infty\otimes\bbK,
 \qquad b\longmapsto\iota_B(b)\otimes1\otimes e_{11}.
\end{equation}
The split unitalization sequence and the $KK$-equivalence $\C\to\mathcal O_\infty$ imply that
\begin{equation}\label{eq:KK-infinite-coefficient-injectivity}
 \iota_*:KK_0(C,B)\lhook\joinrel\longrightarrow
 KK_0(C,B^u\otimes\mathcal O_\infty\otimes\bbK)
\end{equation}
is injective. {The parameterized target $C(S,B^u\otimes\mathcal O_\infty\otimes\bbK)$ contains the full properly infinite projection given by the constant function with value $1_{B^u}\otimes1_{\mathcal O_\infty}\otimes e_{11}$. Since $D$ is nuclear, \cite[Theorem~A]{gabe-classification} therefore gives a homomorphism}
\begin{equation}\label{eq:KK-parameterized-homomorphism}
 \varphi:D\longrightarrow C(S,B^u\otimes\mathcal O_\infty\otimes\bbK),
 \qquad [\varphi]=(C(S,\iota))_*(x).
\end{equation}}{For $s'$ in a sufficiently small neighbourhood of $s$, norm continuity and Lemma~\ref{lem:semiprojective-straightening}(2) give $\varphi_su\simeq\varphi_{s'}u$. Thus the values of \eqref{eq:KK-condensed-shape-evaluation} at $s$ and $s'$ have equal images under \eqref{eq:KK-infinite-coefficient-injectivity}, and hence are equal. This proves the required local constancy.}

 Using the commutative squares
 \[
\xymatrix@C=3em@R=2.5em{
 KK_0(D,B)\ar[r]^{}\ar[d]_{{u^*}}
   &E_0(D,B)\ar[d]^{{u^*}}\\
 KK_0(C,B) \ar[r]
   &E_0(C,B).
 }
\]
for all semiprojective maps {$u:C\to D$} we conclude that
the natural map
\[
KK_{0}(D,B)_{\mathrm{shape}}\to E_{0}(D,B)_{\mathrm{shape}}
\]%
is continuous.
{Via \eqref{gwerfewrfwerferf}, we regard the condensed and shape topologies on both groups as topologies on $KK_0(D,B)$, denoted by $\tau_{\mathrm{cond}}^{KK}$, $\tau_{\mathrm{shape}}^{KK}$, $\tau_{\mathrm{cond}}^E$, and $\tau_{\mathrm{shape}}^E$, respectively. The preceding arguments give} 
   \[
 \tau_{\mathrm{cond}}^{KK}
 =\tau_{\mathrm{cond}}^E
 =\tau_{\mathrm{shape}}^E
 \subseteq\tau_{\mathrm{shape}}^{KK}
 \subseteq\tau_{\mathrm{cond}}^{KK}.
\]
Therefore all coincide.
 \end{proof}

{Let $D$ and $B$ be in $\nCalg_{\sepa}$.}
\begin{theorem}\label{thm:nuclear-qd-hausdorff-comparison}
{If $D$ is nuclear and quasidiagonal, then $\alpha_{D,B}$ induces an isomorphism of topological abelian groups}
\begin{equation}\label{eq:nuclear-qd-hausdorff-comparison}
{(\pi_0{\cR}_D(B))_{\mathrm{Hd}}\xrightarrow{\cong}E_0(D,B)_{\mathrm{Hd}},}
\end{equation}
{where the Hausdorff quotients are taken for the shape topologies.}
\end{theorem}

\begin{proof}
{Choose a shape system $(D_n)$ for $D$.}
Proposition~\ref{prop:intro-splitting} gives a section of
$\alpha_{D,B}$, so the induced map on Hausdorff quotients is
surjective. It remains to prove injectivity.
{We consider the commutative square}
\[
{\xymatrix{
\lim_n\Omega^\infty KK(D_n,B)\ar[r]&\lim_n\Omega^\infty E(D_n,B)\\
\Omega^\infty KK(D,B)\ar[r]^{\simeq}\ar[u]^{\simeq}&\Omega^\infty E(D,B)\ar[u]^{\simeq}.
}}
\]

The right vertical map is an equivalence since $E$-theory preserves countable
filtered colimits. {The left vertical equivalence requires a separate argument because $KK$-theory does not preserve all countable filtered colimits. Every sequential inductive system of separable $C^*$-algebras with nuclear colimit is admissible in the sense of \cite[Definition~2.5]{meyer-nest-localisation}, by the Choi--Effros lifting theorem \cite{choi-effros-lifting} and \cite[Section~2.4, discussion following Proposition~2.6]{meyer-nest-localisation}. Hence the separable $KK$-theory functor preserves the colimit of $(D_n)$ by \cite[Lemma~2.19(3)]{bunke-engel-land-kk}. Applying $\Omega^\infty KK(-,B)$ gives the left vertical equivalence.}

The lower horizontal map is an equivalence because $D$ is nuclear.

{{By Proposition~\ref{prop:kk-e-topologies}, the comparison is an isomorphism of topological groups for the shape topologies and therefore induces an isomorphism on Hausdorff quotients. Using the vertical equivalences and the description of the shape topologies by restriction to the stages, we obtain}
\begin{equation}\label{eq:KK-E-Hausdorff-stage-comparison}
{\lim_n KK_0(D_n,B)\cong KK_0(D,B)_{\mathrm{Hd}}\xrightarrow{\cong}E_0(D,B)_{\mathrm{Hd}}\cong\lim_n E_0(D_n,B).}
\end{equation}
Together with Corollary~\ref{ijwohterhetgertge}, this gives injectivity of
$\lim_n\pi_0\cT_{D_n}(B)\to\lim_n E_0(D_n,B)$.
Corollary~\ref{prop:hausdorff-comparison-criterion} therefore gives the asserted isomorphism of topological groups.}
\end{proof}

%
%

\begin{prob}\label{prob:LLP-comparison} 
{For $D$ and $B$ in $\nCalg_{\sepa}$, we write}
\[
{KL(D,B):=\bigl(KK_0(D,B)_{\mathrm{cond}}\bigr)_{\mathrm{Hd}}}
\]
{for the Hausdorff quotient with respect to the condensed topology, also called the canonical or Pimsner topology \cite{dadarlat-topology}. The comparison in Proposition~\ref{prop:kk-e-topologies} leads to the following questions.}
The first two  are motivated by the attempt to understand the role of  Willett's LLP condition for the Hausdorff comparison results. \begin{enumerate}
\item {If $D$ has LLP, is $KK(D,-)$ exact? Equivalently, is the natural comparison}
\begin{equation}\label{eq:LLP-KK-E-comparison}
{KK(D,-)\longrightarrow E(D,-):\nCalg_{\sepa}\longrightarrow\Sp}
\end{equation}
{an equivalence of functors?}

\item {If $D$ has LLP, does the comparison induce an isomorphism of topological
abelian groups
\[
 \bigl(KK_0(D,B)_{\mathrm{shape}}\bigr)_{\mathrm{Hd}}
 \longrightarrow
 \bigl(E_0(D,B)_{\mathrm{shape}}\bigr)_{\mathrm{Hd}}?
\]}
\item { For which   $D$ and $B$ in $\nCalg_{\sepa}$  do the condensed and shape topologies on $KK_0(D,B)$
coincide?  For nuclear $D$, the comparison of $KK$ with $E$ is an equivalence, and Proposition~\ref{prop:kk-e-topologies} identifies the condensed and shape topologies. For general $D$, we do not know whether the condensed topology defining $KL(D,B)$ coincides with the shape topology.}

\end{enumerate}
\end{prob}

\subsection{{Hausdorff comparison}}
\label{sec:hausdorff-source-continuity}


The  $\infty$-category $\cU$ from Definition \ref{def:universal-admissible} admits countable filtered colimits by Proposition \ref{prop:universal-unreduced-category}. 
Hence we can apply the theory of shapes from  \cite[Section~2.2]{budu}.
We have an adjunction $$|-|: \Ind_{\aleph_1}(\mathcal U) \leftrightarrows \cU:y_{\mathcal U}\ ,$$where $|-|$ is called realization and $y_{\cU}$ is the 
 Yoneda embedding. According to  \cite[Definition~2.7]{budu},
a shape of $A$ in $\cU$, if it exists, is an object $S_{\cU}(A)$ in $ \Ind_{\aleph_1}(\mathcal U)$ together with an equivalence
of functors
$$\Hom_{ \Ind_{\aleph_1}(\mathcal U)}(S_{\cU}(A),-)\simeq \Hom_{\cU}(A,|-|): \Ind_{\aleph_1}(\mathcal U)\to \Ani  \ .$$Below we use the criterion   \cite[Proposition~2.22]{budu}: The object $A$ in $\cU$ admits a shape if and only if there exists a  system $(A_{n})_{n\in \nat}$ in $\cU$  with ${\colim_n A_n}\simeq A$ and such that the canonical map is an equivalence
\begin{equation}\label{eq:U-shape-criterion}
{\lim_n \colim_i\Hom_{\cU}(A_{n},X_{i})\stackrel{\simeq}{\to}\Hom_{\cU}(A,\colim_iX_{i})}
\end{equation}
for all countable filtered systems $(X_{i})_{i\in I}$ in $\cU$. In this case ${S_{\cU}(A)}\simeq {\colim_n} y_{\cU}(A_{n})$.

\begin{lem}
\label{lem:U-shapes}
\mbox{}\begin{enumerate}
\item 
{Every object of $\mathcal U$ admits a shape.  \item If $(D_n)$ is a shape system of  $D$ in $\nCalg_{\sepa}$, then}
\begin{equation}\label{eq:U-algebra-shape}
{S_{\mathcal U}(j(D))\simeq\colim_n y_{\mathcal U}(j(D_n)).}
\end{equation}\end{enumerate}
\end{lem}
\begin{proof} In the following we  use the inclusions $\cU\subseteq \mathcal P_{\mathrm{cont}}\subseteq \cP$ from the  proof of Proposition~\ref{prop:universal-unreduced-category} in order to consider objects  of $\cU$ also as objects of ${\mathcal P_{\mathrm{cont}}}$ or $\cP$. Let $(X_i)_{i\in I}$ be a  countable filtered system in $\mathcal U$.
Then by the localization \eqref{gerferfwerfwer} we have \begin{equation}\label{gwrefwerfwref}\colim^{\cU}_{i} X_{i} \simeq \mathsf L(\colim_i^{\cP} X_i)\ .
\end{equation}  Let $D$ be in $\nCalg_{\sepa}$ and $(D_{n})$ be a shape system of $D$.  Using  Yoneda  and the notation \eqref{eq:additive-yoneda-presheaf} at the equivalences marked by $!$
\begin{equation}\label{eq:U-shape-mapping-calculation}
 {\begin{array}{>{\displaystyle}r >{\displaystyle}c >{\displaystyle}l}
\multicolumn{3}{l}{\displaystyle\Hom_{\mathcal U}(j(D),\colim^{\cU}_iX_i)}\\
&\stackrel{\eqref{gerwfgwergwerg},\eqref{gwrefwerfwref}}{\simeq}&\Hom_{\mathcal P_{\mathrm{cont}}}
 (\mathsf Lh_{L_{h,K}^{\add}(D)},\mathsf L(\colim^{\cP}_iX_i))\\
&\stackrel{\eqref{gerferfwerfwer}}{\simeq}&\Hom_{\mathcal P}
 (h_{L_{h,K}^{\add}(D)},\incl\mathsf L(\colim^{\cP}_iX_i))\\
&\stackrel{!}{\simeq}&\bigl(\mathsf L(\colim^{\cP}_iX_i)\bigr)(L_{h,K}^{\add}(D))\\
&\stackrel{\eqref{eq:unreduced-presheaf-reflection},\eqref{gwrefwerfwref}}{\simeq}&\lim_n\colim_i X_i(L_{h,K}^{\add}(D_n))\\
&\stackrel{!}{\simeq}&\lim_n\colim_i\Hom_{\mathcal P}
 (h_{L_{h,K}^{\add}(D_n)},\incl X_i)\\
&\stackrel{\eqref{gerferfwerfwer}}{\simeq}&\lim_n\colim_i\Hom_{\mathcal P_{\mathrm{cont}}}
 (\mathsf Lh_{L_{h,K}^{\add}(D_n)},X_i)\\
&\stackrel{\eqref{gerwfgwergwerg} }{\simeq}&\lim_n\colim_i\Hom_{\mathcal U}(j(D_n),X_i).
\end{array}}
\end{equation}
Since $\colim_{n}j(D_{n})\simeq j(D)$, the equivalence \eqref{eq:U-shape-mapping-calculation} 
verifies {the criterion \eqref{eq:U-shape-criterion}} for the existence of the shape of $j(D)$ and proves
 \eqref{eq:U-algebra-shape}.

 The   additive completion functor ${L^{\add}_{h,K}}$ from \eqref{eq:additive-completion-category} is essentially surjective. Thus the construction of $\mathcal U$ in Proposition~\ref{prop:universal-unreduced-category} generates $\mathcal U$ from the objects $j(D)$ under countable filtered colimits. Objects admitting shapes are closed under these colimits by \cite[Lemma~2.9]{budu}. Hence every object of $\mathcal U$ admits a shape.
\end{proof}

{For a functor $F:\mathcal U^{\op}\to\CGrp(\Ani)$, we equip $\pi_0F(V)$ with the coarsest group topology making restriction along every weakly compact map $W\to V$ continuous when $\pi_0F(W)$ is discrete. This is the shape topology, as in \cite[Definition~4.19]{budu}. If $V=j(D)$, \eqref{eq:U-algebra-shape} and \cite[Corollary~2.23 and Lemma~4.21]{budu} show that the restrictions to $j(D_n)$ define this topology. It therefore agrees with Definition~\ref{def:intro-shape-topology} for the functor $F\circ j$, by Remark~\ref{okpzehhtre}.}

\begin{prop}
\label{prop:hausdorff-source-continuity}
{Let $F:\mathcal U^{\op}\to\CGrp(\Ani)$ send countable filtered colimits to limits. For every system $(V_i)_{i\in\nat}$ in $\mathcal U$, with $V:=\colim_iV_i$, restriction induces an isomorphism of topological abelian groups}
\begin{equation}\label{eq:hausdorff-source-continuity}
{(\pi_0F(V))_{\mathrm{Hd}}\xrightarrow{\cong}
\lim_i(\pi_0F(V_i))_{\mathrm{Hd}}.}
\end{equation}
\end{prop}
\begin{proof}
{We use the diagonal argument of \cite[Lemma~4.35 and Proposition~4.36]{budu}. We explain why it applies to $F$. If a sequence $(W_n)$ represents the shape of $W$, then $W\simeq\colim_nW_n$, so $F(W)\simeq\lim_nF(W_n)$. The Milnor exact sequence gives}
\begin{equation}\label{eq:hausdorff-source-milnor}
{0\longrightarrow{\lim\nolimits_n^1}\pi_1F(W_n)
\xrightarrow{\partial}\pi_0F(W)
\xrightarrow{r}\lim_n\pi_0F(W_n)\longrightarrow0.}
\end{equation}
{The shape topology on $\pi_0F(W)$ is induced by $r$, where the target has the inverse-limit topology of discrete groups. Thus $r$ in \eqref{eq:hausdorff-source-milnor} realizes the Hausdorff quotient, and the image of $\partial$ is $\ker( r)=\overline{\{0\}}$. We obtain}
\begin{equation}\label{eq:hausdorff-shape-stage-formula}
{(\pi_0F(W))_{\mathrm{Hd}}\cong
\lim_n\pi_0F(W_n)_{\mathrm{discrete}}.}
\end{equation}
{By Lemma~\ref{lem:U-shapes}, all $V_i$ admit shapes. {By the construction in the proof of \cite[Proposition~4.36]{budu}, we choose a diagram}
\begin{equation}\label{eq:compatible-U-shape-presentations}
{(i,n)\longmapsto V_{i,n}:\nat\times\nat\longrightarrow\cU,
\qquad \colim_n y_{\cU}(V_{i,n})\simeq S_{\cU}(V_i),}
\end{equation}
{with these equivalences natural in $i$. The diagonal $(V_{n,n})_n$ of \eqref{eq:compatible-U-shape-presentations} represents the shape of $V$.} Applying \eqref{eq:hausdorff-shape-stage-formula} gives}
\begin{equation}\label{eq:hausdorff-source-diagonal}
{\begin{aligned}
(\pi_0F(V))_{\mathrm{Hd}}
&\cong\lim_n\pi_0F(V_{n,n})_{\mathrm{discrete}}\\
&\cong\lim_i\lim_n\pi_0F(V_{i,n})_{\mathrm{discrete}}\\
&\cong\lim_i(\pi_0F(V_i))_{\mathrm{Hd}}.
\end{aligned}}
\end{equation}
{The middle isomorphism uses cofinality of the diagonal and commutation of limits. All identifications use restriction maps and the inverse-limit topologies, so \eqref{eq:hausdorff-source-diagonal} is the claimed natural topological isomorphism.}
\end{proof}

 {Fix $B$ in $\nCalg_{\sepa}$. We apply Proposition~\ref{prop:hausdorff-source-continuity} to the functors}
\begin{equation}\label{eq:hausdorff-U-functors}
 {\begin{aligned}
F_B,G_B:\mathcal U^{\op}&\longrightarrow\CGrp(\Ani),\\
F_B(V)&:=\Hom_{\mathcal U}(V,j^{\sfred}(B)),\\
G_B(V)&:=\Hom_{\EE_{\sepa}}(e_{\mathcal U}(V),e_{\sepa}(B)).
\end{aligned}}
\end{equation}
 {Both functors send countable filtered colimits to limits: for $G_B$ we use preservation of these colimits by $e_{\mathcal U}$ from \eqref{eq:universal-E-factorization}. Applying $e_{\mathcal U}$ and using $e_{\mathcal U}(j^{\sfred}(B))\simeq e_{\sepa}(B)$ gives a natural transformation}
\begin{equation}\label{eq:hausdorff-U-comparison}
 \alpha^{\mathcal U}_{-,B}:F_B\longrightarrow G_B:\cU^{\op} \to
 {\CGrp(\Ani)\ .} 
\end{equation}
 {At $V=j(D)$, \eqref{eq:hausdorff-U-comparison} is \eqref{eq:intro-alpha}.}

\begin{kor}\label{cor:hausdorff-comparison-colimits}
The full subcategory of objects $V$ of $\mathcal U$ for which \eqref{eq:hausdorff-U-comparison} induces an isomorphism on Hausdorff component groups for every $B$ is closed under retracts, finite products, and countable filtered colimits. 
\end{kor}
\begin{proof}
{Since $\mathcal U$ is additive, finite products are also coproducts. The functors $F_B$ and $G_B$ therefore send finite products in $\mathcal U$ to products, compatibly with the shape topologies.  Since Hausdorff quotients commute with finite products we get the closedness under finite products.}
If $V$ is a retract of $W$, naturality of \eqref{eq:hausdorff-U-comparison} exhibits the Hausdorff comparison for $V$ as a retract of that for $W$ in the category of morphisms of topological abelian groups. A retract of an isomorphism is an isomorphism, proving closure under retracts.

 For $V=\colim_iV_i$, naturality of \eqref{eq:hausdorff-source-continuity} gives the commutative square\begin{equation}\label{eq:hausdorff-comparison-limit-square}
 \xymatrix@C=3em@R=2.5em{
(\pi_0F_B(V))_{\mathrm{Hd}}\ar[r]\ar[d]_{\cong}
 &(\pi_0G_B(V))_{\mathrm{Hd}}\ar[d]^{\cong}\\
\lim_i(\pi_0F_B(V_i))_{\mathrm{Hd}}\ar[r]
 &\lim_i(\pi_0G_B(V_i))_{\mathrm{Hd}}.
}
\end{equation}
 If the comparisons at all $V_i$ are isomorphisms, then so is the lower map in \eqref{eq:hausdorff-comparison-limit-square}, and hence the upper map. 
\end{proof}

\begin{proof}[{Proof of Theorem~\ref{thm:overview-one-hausdorff}}]
 {Let $\mathcal V$ be the full subcategory in Corollary~\ref{cor:hausdorff-comparison-colimits}. Theorem~\ref{thm:nuclear-qd-hausdorff-comparison} shows that $j(D')$ belongs to $\mathcal V$ for every separable nuclear quasidiagonal algebra $D'$. Corollary~\ref{cor:hausdorff-comparison-colimits} supplies the closure properties in Definition~\ref{def:U-nq}, so}
\begin{equation}\label{eq:hausdorff-comparison-class-containment}
 \mathcal U_{\mathrm{nq}}\subseteq\mathcal V.
\end{equation}
 This proves the Hausdorff comparison in Theorem~\ref{thm:overview-one-hausdorff}. 

By Theorem~\ref{thm:intro-exactness}, ${\cR}_D$ is homological, so Proposition~\ref{prop:intro-splitting} applies. {Under this splitting, we have}
\begin{equation}\label{eq:hausdorff-remainder-kernel}
{\pi_0{\widehat{\cR}}_D(B)=\ker(\pi_0\alpha_{D,B}).}
\end{equation} Injectivity on Hausdorff quotients places this kernel in the closure of zero.
\end{proof}

\section{Quasi-representations and Willett's group}
\label{sec:willett}

\subsection{Controlled quasi-representations and {a shape system} for $C^{*}_{\max}(G)$}
\label{kohperthertgrtgertgretg}
Let $G$ be a countable discrete group, $S$ be a finite
symmetric subset of $G$ containing $1$, and  $\epsilon$ be in
$(0,\infty)$. We then form the
universal unital $C^*$-algebra  $C^{*}(S,\epsilon)$   generated by  a family of unitaries $(u_s)_{s\in S}$
subject to
\[
 u_1=1,\qquad u_{s^{-1}}=u_s^*,\qquad
 \|u_su_t-u_{st}\|\leq\epsilon\quad{(\forall s,t\in S\text{ with }st\in S)}.
\]

We introduce the poset $\mathcal I_G$ of pairs $(S,\epsilon)$,
where $S$ is a finite symmetric subset of $G$ containing $1$
and $\epsilon$ is in $(0,\infty)$. Its order is given by
\begin{equation}\label{eq:group-control-order}
 (S,\epsilon)\leq(T,\delta)
 \quad\Longleftrightarrow\quad
 S\subseteq T\ \text{and}\ \delta\leq\epsilon.
\end{equation}
In this case,  the universal property
gives a unital homomorphism
\[
 C^{*}(S,\epsilon)\longrightarrow C^{*}(T,\delta),
 \qquad u_s\longmapsto u_s\quad{(\forall s\in S)}.
\]
  Thus we obtain
a functor
\[
 C^{*}:\mathcal I_G\longrightarrow\Calg_{\sepa},
 \qquad (S,\epsilon)\longmapsto C^{*}(S,\epsilon).
\]
By comparing universal properties, we conclude that   \begin{equation}\label{eq:group-algebra-colimit}
 \colim_{(S,\epsilon)\in\mathcal I_G}C^{*}(S,\epsilon)
 \cong C^{*}_{\max}(G).
\end{equation}
 
Since $G$ is countable, we can choose an increasing sequence
$(S_n)_{n\in\nat_{\geq1}}$ of finite symmetric subsets containing
$1$ and exhausting $G$. Then the  sequence $(S_n,1/n)_n$ is cofinal
in $\mathcal I_G$. We put
\[
 D_n:=C^{*}(S_n,1/n).
\]
For {$S\subseteq T$ and} $\delta<\epsilon$, the transition maps
$C^{*}(S,\epsilon)\to C^{*}(T,\delta)$ are semiprojective by the finite-relations lifting argument of \cite[Proposition~2.35]{blackadar-shape}; see also \cite[Construction~3.3 and Proposition~3.5]{budu}. 
Thus by  \eqref{eq:group-algebra-colimit} and cofinality  we obtain a
shape system $(D_{n})$
\[
 D_1\longrightarrow D_2\longrightarrow D_3\longrightarrow\cdots,
 \qquad \colim_nD_n\cong C^{*}_{\max}(G)
\]
of $C^{*}_{\max}(G)$.

{Let $B$ be in $\nCalg_{\sepa}$. By \eqref{eq:R-shape-formula} and Definition~\ref{def:intro-T}, we have} 
\[
 {\cR}_{C^*_{\max}(G)}(B)
 \simeq\lim_n\cT_{D_n}(B)
 \simeq \lim_n\Fib_0\left(
 {\cM}_{D_n}((B\otimes\bbK)^u)^{\group}
 \xrightarrow{(\epsilon_{B\otimes\bbK})_*}
 {\cM}_{D_n}(\C)^{\group}\right).
\]
Let $B$ be in $\Calg_{\sepa}$.
\begin{kor}\label{cor:group-hausdorff-formula}
{  There is a canonical isomorphism of topological abelian groups
\begin{equation}\label{eq:group-hausdorff-formula}
 \bigl(\pi_0{\cR}_{C^*_{\max}(G)}(B)\bigr)_{\mathrm{Hd}}
 \cong\lim_{(S,\epsilon)\in\mathcal I_G}
 \pi_0\bigl({\cM}_{C^*(S,\epsilon)}(B)^{\group}\bigr)_{\mathrm{discrete}}.
\end{equation}
}
\end{kor}
\begin{proof}
{Apply Corollary~\ref{cor:unital-hausdorff-formula} to the shape system $(C^*(S_n,1/n))_n$ and use its cofinality in $\mathcal I_G$.}
\end{proof}

\subsection{Homotopy and stable approximate unitary equivalence}
\label{sec:willett-equivalence}

We recall the construction of Willett's quasi-representation
group \cite[Definitions~2.1, 2.5 and 2.7--2.11]{willett}.
We use the control poset $\mathcal I_G$ introduced in
\eqref{eq:group-control-order}; symmetric finite subsets are
cofinal among Willett's finite control sets.
Fix $(S,\epsilon)$ in $\mathcal I_G$.
An \emph{$(S,\epsilon)$-representation} is a map
$\varphi:S\to U(V_\varphi)$, where $V_\varphi$ is a
finite-dimensional complex Hilbert space, satisfying
\[
 \varphi(1)=1,\qquad
 \|\varphi(s)\varphi(t)-\varphi(st)\|<\epsilon
 \quad{(\forall s,t\in S\text{ with }st\in S)}.
\]
The zero-dimensional space gives the zero representation.
For two such representations, we define block sum by
$(\varphi\oplus\psi)(s):=\varphi(s)\oplus\psi(s)$.

We call two such maps $\varphi,\psi$ \emph{stably
$(S,\epsilon)$-equivalent} and write
$\varphi\sim_{S,\epsilon}\psi$ if there are an
$(S,\epsilon)$-representation $\theta$ and a unitary
\[
 v:V_\varphi\oplus V_\theta
       \longrightarrow V_\psi\oplus V_\theta
\]
such that
\begin{equation}\label{eq:willett-stable-equivalence}
 \max_{s\in S}
 \bigl\|v(\varphi(s)\oplus\theta(s))v^*
           -(\psi(s)\oplus\theta(s))\bigr\|<\epsilon.
\end{equation}
Stable $(S,\epsilon)$-equivalence is an equivalence relation
compatible with block sum \cite[Lemmas~2.6 and 2.8]{willett}.
We denote the resulting monoid and its group completion by
\[
 V_\epsilon(S):=
 \{(S,\epsilon)\text{-representations}\}/\sim_{S,\epsilon},
 \qquad K(S,\epsilon):=V_\epsilon(S)^{\group}.
\]
The monoid $V_\epsilon(S)$ is commutative, with the zero
representation as its unit.

For $(S,\epsilon)\leq(T,\delta)$ in $\mathcal I_G$,
restriction gives transition homomorphisms
$K(T,\delta)\to K(S,\epsilon)$. With the discrete topology on
each stage group, these form an inverse system of topological
abelian groups indexed by $\mathcal I_G$.

We follow Willett's definition of $QR(G)$ in
\cite[Definition~2.11, p.~11]{willett}, using the inverse system
introduced in \cite[Definition~2.10]{willett}.

\begin{ddd}\label{def:willett-QR}We define  the topological group of quasi-representations of $G$ by 
\[
 QR(G):=\lim_{(S,\epsilon)\in\mathcal I_G}K(S,\epsilon)_{\mathrm{discrete}}.
\]
 
\end{ddd}



\begin{proof}[Proof of Theorem~\ref{thm:overview-qr-comparison}]
 {Corollary~\ref{cor:group-hausdorff-formula}, applied with $B=\C$, gives an isomorphism of topological abelian groups
\begin{equation}\label{eq:theorem-C-mapping-groups}
 \bigl(\pi_0{\cR}_{C^*_{\max}(G)}(\C)\bigr)_{\mathrm{Hd}}
 \cong\lim_{(S,\epsilon)\in\mathcal I_G}
 \pi_0\bigl({\cM}_{C^*(S,\epsilon)}(\C)^{\group}\bigr)_{\mathrm{discrete}}.
\end{equation}
It remains to identify the right-hand side of \eqref{eq:theorem-C-mapping-groups} with $QR(G)$.}

For a finite symmetric subset $S$ of $G$ containing $1$ and
$\epsilon$ in $(0,\infty)$, we introduce
\begin{equation}\label{gerwfegwergwerf}H_\epsilon(S):=\pi_0\bigl({\cM}_{C^{*}(S,\epsilon)}(\C)^{\group}\bigr)
 \cong \bigl(\pi_0{\cM}_{C^{*}(S,\epsilon)}(\C)\bigr)^{\group}.
\end{equation}
We will construct mutually inverse morphisms between the pro-systems
$(H_\epsilon(S))_{S,\epsilon}$ and $(K(S,\epsilon))_{S,\epsilon}$.

We first describe the representatives of $H_\epsilon(S)$.
A homomorphism
$\rho:C^{*}(S,\epsilon)\to\bbK$ has support
\[
 p:=\rho(1),\qquad p^2=p=p^*,\qquad N:=\rk(p)<\infty.
\]
After identifying $p\ell^2\cong\C^N$, it determines a family
$a_s\in U(N)$ satisfying
\[
 a_1=1_N,\quad a_{s^{-1}}=a_s^*,\quad
 d_S(a):=\max_{s,t,st\in S}\|a_sa_t-a_{st}\|\leq\epsilon.
\]
Conversely such a family defines $\rho$ by $u_s\mapsto a_s$ on its
support and zero on the orthogonal complement. Furthermore,
\[
 \pi_0{\cM}_{C^{*}(S,\epsilon)}(\C)
 =\pi_0\Hom_{\nCalg_{\sepa,h}}(C^{*}(S,\epsilon),\bbK).
\]
Along a point-norm homotopy $\rho_u$, the projections $p_u:=\rho_u(1)$
have constant finite rank. Trivializing their ranges over $[0,1]$ turns
$\rho_u$ into a continuous path of such families in a fixed $U(N)$.
Its homotopy class is independent of the trivialization.
 
We next pass from homotopy to Willett's stable approximate
unitary equivalence. Fix $\delta$ in $(0,\epsilon)$ and a path
$a(u)$ of normalized families with $d_S(a(u))\leq\delta$.
Uniform continuity gives a subdivision
$0=u_0<\cdots<u_m=1$ {with $m\geq2$} such that, writing $a_j:=a(u_j)$,
\[
 \max_{s\in S}\|a_j(s)-a_{j-1}(s)\|<\epsilon
 \quad{(\forall j\in\{1,\ldots,m\})}.
\]
We use the Willett--Yu “speed for space” argument
\cite[proof of Lemma~4.15]{willett-yu-uct}: adding
$\theta:=\bigoplus_{j=1}^{m-1}a_j$ to both endpoints allows a
cyclic permutation unitary $v$ to pair consecutive, hence
$\epsilon$-close, samples of the path:
\[
 v(a_0\oplus\theta)v^*
 =a_{m-1}\oplus a_0\oplus\cdots\oplus a_{m-2},
 \qquad
 a_m\oplus\theta=a_m\oplus a_1\oplus\cdots\oplus a_{m-1}.
\]
Consequently,
\[
 \max_{s\in S}
 \|v(a_0(s)\oplus\theta(s))v^*
       -(a_m(s)\oplus\theta(s))\|
 \leq\max_{1\leq j\leq m}\max_{s\in S}
       \|a_{j-1}(s)-a_j(s)\|<\epsilon.
\]
Thus homotopic $(S,\delta)$-families are stably
$(S,\epsilon)$-equivalent. We write $[a]_W$ for the image in
$K(S,\epsilon)$ of the stable $(S,\epsilon)$-equivalence class
of $a$. We therefore obtain a homomorphism
\[
 a_{\delta,\epsilon}:H_\delta(S)\longrightarrow K(S,\epsilon),
 \qquad [a]-[b]\longmapsto[a]_W-[b]_W.
\]

For the reverse direction, we fix $\delta$ in
$(0,\min\{1/100,\epsilon/100\})$ and first normalize a Willett
$(S,\delta)$-representation $\varphi$ {as in \cite[Lemma~2.4]{willett}}.
{Normalization is necessary because $C^*(S,\epsilon)$ imposes the exact relations $u_{s^{-1}}=u_s^*$, whereas Willett's representations satisfy these relations only approximately.}
We set  $(N\varphi)(1):=1$.
 We can choose  one element  $s$ from
each inverse pair  in $G$ with $s\neq s^{-1}$. On the pairs contained in $S$ we define
\[
 (N\varphi)(s):=\varphi(s),\qquad
 (N\varphi)(s^{-1}):=\varphi(s)^*.
\]
For an involution $s$ in $S\setminus\{1\}$, we instead define
\[
 (N\varphi)(s):=\operatorname{sign}(\Ree\varphi(s)).
\]
The bound $\|\varphi(s)^2-1\|<\delta<1$ makes this functional
calculus well-defined. We claim that
\[
 \max_{s\in S}\|(N\varphi)(s)-\varphi(s)\|<\delta,
 \qquad d_S(N\varphi)<4\delta.
\]
For an inverse pair, the first estimate follows from
\[
 \|\varphi(s^{-1})-\varphi(s)^*\|
 =\|\varphi(s)\varphi(s^{-1})-1\|<\delta.
\]
For an involution, we use
$|z-1||z+1|<\delta$ on $\operatorname{spec}(\varphi(s))\subseteq S^1$.
For $s,t,st$ in $S$, the triangle inequality gives the defect
estimate:
\begin{align*}
 \|(N\varphi)(s)(N\varphi)(t)-(N\varphi)(st)\|
 &\leq\|(N\varphi)(s)-\varphi(s)\|
       +\|(N\varphi)(t)-\varphi(t)\|\\
 &\quad+\|\varphi(s)\varphi(t)-\varphi(st)\|
       +\|\varphi(st)-(N\varphi)(st)\|<4\delta.
\end{align*}
Normalization commutes with direct sums and unitary conjugation.

To check compatibility with stable equivalence, suppose that
$\varphi$ and $\psi$ are stably $(S,\delta)$-equivalent.
We choose a common $(S,\delta)$-representation $\theta$ and a
unitary $v$ realizing this relation:
\[
 \max_{s\in S}
 \|v(\varphi(s)\oplus\theta(s))v^*
       -(\psi(s)\oplus\theta(s))\|<\delta.
\]
The corresponding normalized families
\[
 A_s:=v\bigl((N\varphi)(s)\oplus(N\theta)(s)\bigr)v^*,
 \qquad B_s:=(N\psi)(s)\oplus(N\theta)(s)
\]
satisfy $d_S(A),d_S(B)<4\delta$ and
$\max_s\|A_s-B_s\|<3\delta$. We join them by the path
\[
 h_s(u):=A_s\exp\bigl(u\log(A_s^*B_s)\bigr),
\]
for $u$ in $[0,1]$, using the principal logarithm.
Since $3\delta<1/2$, we have
\[
 \max_{s,u}\|h_s(u)-A_s\|
 \leq2\max_s\|A_s-B_s\|<6\delta,
 \qquad d_S(h(u))<4\delta+18\delta<\epsilon.
\]
The path preserves the normalization. In fact $h_1(u)=1$ and
$h_{s^{-1}}(u)=h_s(u)^*$; the latter follows from the functional-calculus
identity
\[
 \log(A_sB_s^*)=-A_s\log(A_s^*B_s)A_s^*.
\]
In particular, the same path remains self-adjoint when $s=s^{-1}$.
It therefore defines a homotopy of homomorphisms from $C^{*}(S,\epsilon)$.
{We concatenate this homotopy with a conjugation path, using connectedness of the finite-dimensional unitary group.}
Cancelling the common summand $[N\theta]$ in $H_\epsilon(S)$
then gives $[N\varphi]=[N\psi]$. Thus normalization induces a
homomorphism
\[
 b_{\delta,\epsilon}:K(S,\delta)\longrightarrow H_\epsilon(S),
 \qquad [\varphi]_W-[\psi]_W
 \longmapsto[N\varphi]-[N\psi].
\]

It remains to check that these maps are inverse as morphisms
of pro-systems. We choose $\delta$ in
$(0,\min\{1/100,\epsilon/100\})$, $\gamma$ in $(0,\delta)$,
and $\eta$ in $(\epsilon,\infty)$. The composites are
\[
 b_{\delta,\epsilon}a_{\gamma,\delta}
   =\bigl(H_\gamma(S)\to H_\epsilon(S)\bigr),\qquad
 a_{\epsilon,\eta}b_{\delta,\epsilon}
   =\bigl(K(S,\delta)\to K(S,\eta)\bigr).
\]
The first equality uses $Na=a$ for normalized families.
For the second, the bound $\|N\varphi-\varphi\|_S<\delta<\eta$
gives stable $(S,\eta)$-equivalence by taking the
zero-dimensional representation and the identity unitary in
\eqref{eq:willett-stable-equivalence}, with tolerance $\eta$.
{All constructions commute with restriction to inverse-closed subsets and with the transition maps obtained by increasing the tolerance.} These subsets are cofinal among
finite control sets. We conclude that
\[
 (H_\epsilon(S))_{S,\epsilon}
 \cong(K(S,\epsilon))_{S,\epsilon}
 \quad\text{as pro-abelian groups}.
\]
{Taking limits gives
\[
 \lim_{(S,\epsilon)\in\mathcal I_G}H_\epsilon(S)_{\mathrm{discrete}}
 \cong\lim_{(S,\epsilon)\in\mathcal I_G}K(S,\epsilon)_{\mathrm{discrete}}
 =QR(G)
\]
as topological abelian groups. 
Together with \eqref{eq:theorem-C-mapping-groups}  and \eqref{gerwfegwergwerf}, this proves the theorem.}
\end{proof}

\section{Relations with other approaches and further questions}
\label{sec:research-context}

 The purpose of this section is to provide further relations of the present work with constructions from the $C^{*}$-literature.

\subsection{Deformation $K$-theory and $K$-theory of $BG$}

For a discrete group $G$, Carlsson's connective
\emph{deformation $K$-theory spectrum} $K^{\mathrm{def}}(G)$
is characterized by the equivalence in $\CGrp(\Ani)$
\[
 \Omega^\infty K^{\mathrm{def}}(G)\simeq
 \left(\coprod_{n\geq0}
 EU(n)\times_{U(n)}\Hom(G,U(n))\right)^{\group},
\]
where we equip $\Hom(G,U(n))$ with the topology of pointwise
convergence, and $U(n)$ acts by conjugation.
Here we group-complete the commutative monoid in anima
presented by the topological groupoid of finite-dimensional
unitary representations with block sum
\cite[Definitions~2.1--2.2 and Section~2]{ramras-deformation-excision}.
See also \cite{lawson-deformation,ramras-deformation-excision}
for multiplicative structures and excision results.

For the rest of this subsection, let $G$ be countable. We write
\[
 KU^{BG}:=\operatorname{map}_{\Sp}(\Sigma^\infty BG_+,KU)
\]
for the complex $K$-cohomology spectrum of $BG$ and
  give $\pi_0KU^{BG}$ the initial topology for the
restriction maps to $K^0(Y)$, where $Y$ ranges over the finite
subcomplexes of $BG$ and each $K^0(Y)$ carries the discrete
topology.
We discuss the commutative diagram
\begin{equation}\label{eq:deformation-quasirep-coassembly}
 \xymatrix@C=5em@R=3.5em{
 \Omega^\infty K^{\mathrm{def}}(G)
       \ar[r]^{\gamma_G}\ar[d]_{\mathrm{AS}_G}
   &{\cR}_{C^*_{\max}(G)}(\C)\ar@{..>}[dl]
       \ar[d]^{\alpha_{C^*_{\max}(G),\C}}\\
 \Omega^\infty KU^{BG}
   &\Omega^\infty E(C^*_{\max}(G),\C)
       \ar[l]_{\Omega^\infty\mu_G^*}
 }
\end{equation}
in $\CGrp(\Ani)$,
{The upper horizontal map $\gamma_G$ sends a finite-dimensional representation to the  compatible family of restrictions to $C^*(S,\epsilon)$. 
The right vertical map is the comparison
$\alpha_{C^*_{\max}(G),\C}$ constructed in
\eqref{eq:intro-alpha}. 
The lower horizontal map is induced by coassembly,
defined by pairing with the Mishchenko bundle over finite
subcomplexes of $BG$. The dotted arrow is its composite with
$\alpha_{C^*_{\max}(G),\C}$.
 {Finally, the left vertical map $\mathrm{AS}_G$ is the topological Atiyah--Segal map studied by Ramras \cite{ramras-atiyah-segal}. Commutativity of \eqref{eq:deformation-quasirep-coassembly} expresses that pairing a genuine representation with the Mishchenko bundle produces its associated flat vector bundle.}
Diagram~\eqref{eq:deformation-quasirep-coassembly} separates
the effects of passing from genuine representations to
quasi-representations, from quasi-representations to
$E$-theory classes, and then applying coassembly.

If $C^*_{\max}(G)\in \cU_{\mathrm{nq}}$,   then $\alpha_{C^{*}_{\max}(G),\C}$ is split surjective by
Proposition~\ref{prop:intro-splitting} and induces an
isomorphism on Hausdorff quotients of component groups by
Theorem~\ref{thm:overview-one-hausdorff}.

{If $G$ is torsion-free and has the Haagerup property, coassembly $\Omega^\infty\mu_G^*$ is an equivalence by the Higson--Kasparov theorem \cite{higson-kasparov}. Here $C^*_{\max}(G)$ is $KK$-equivalent to a nuclear algebra \cite[Section~3, p.~169]{dadarlat-group-quasirepresentations}, so the comparison from $KK$-theory to $E$-theory is an equivalence. For the Hausdorff version of coassembly, see \cite[Proposition~10.6]{willett}.}

\begin{kor}\label{cor:Haagerup-U-nq-coassembly}
{Assume that $G$ is torsion-free, has the Haagerup property, and satisfies $C^*_{\max}(G)\in\cU_{\mathrm{nq}}$. Then the dotted composite in \eqref{eq:deformation-quasirep-coassembly} is split surjective and induces an isomorphism of topological abelian groups}
\begin{equation}\label{eq:Haagerup-U-nq-Hausdorff-coassembly}
{QR(G)\cong
 \bigl(\pi_0{\cR}_{C^*_{\max}(G)}(\C)\bigr)_{\mathrm{Hd}}
 \xrightarrow{\cong}(\pi_0KU^{BG})_{\mathrm{Hd}}.}
\end{equation}
\end{kor}
\begin{proof}
{Coassembly is an equivalence by the preceding argument. We combine Proposition~\ref{prop:intro-splitting} and Theorems~\ref{thm:overview-one-hausdorff} and~\ref{thm:overview-qr-comparison}.}
\end{proof}

{Thus we obtain the comparison of \cite[Theorem~10.8]{willett} under the combined Haagerup and $\cU_{\mathrm{nq}}$ hypotheses. Willett assumes that $G$ is countable, torsion-free and has the Haagerup property, and that $C^*_{\max}(G)$ is quasidiagonal and has LLP; his conclusion is an isomorphism of topological rings.}

\begin{ex}\label{ex:Haagerup-U-nq-free-products}
{The group}
\begin{equation}\label{eq:Haagerup-U-nq-group-examples}
{G:=\prod_{i=1}^k\left(*_{j\in J_i}G_{ij}\right),}
\end{equation}
{with $k$ finite, each $J_i$ countable, and each $G_{ij}$ countable, torsion-free and amenable, satisfies the assumptions of Corollary~\ref{cor:Haagerup-U-nq-coassembly}. Indeed, such groups are torsion-free and have the Haagerup property, since these properties are preserved by countable free products and finite direct products. The algebras $C^*_{\max}(G_{ij})$ are nuclear and quasidiagonal. Lemma~\ref{lem:U-free-product-closure} and Proposition~\ref{prop:U-nq-max-tensor} therefore give $C^*_{\max}(G)\in\cU_{\mathrm{nq}}$.}
\hB
\end{ex}

\begin{ex}\label{ex:Haagerup-U-nq-non-LLP}
{The group $G:=F_r\times F_s$ 
satisfies the assumptions of Corollary~\ref{cor:Haagerup-U-nq-coassembly} for integers $r,s\geq2$. Indeed, it is torsion-free and has the Haagerup property, and $C^*_{\max}(G)\in\cU_{\mathrm{nq}}$ by Example~\ref{ex:U-nq-tensor-products}. Hence the Hausdorff comparison with $K$-cohomology of $BG$ applies. In contrast, Willett's Theorem~10.8 does not apply, since $C^*_{\max}(G)$ fails LLP by \cite[Theorem~A]{roydor-non-llp}.}
\hB
\end{ex}


Elements of $\pi_0K^{\mathrm{def}}(G)$ are differences of
deformation classes of finite-dimensional representations.
In general, these do not generate all components of
${\cR}_{C^*_{\max}(G)}(\C)$.

\begin{ex}\label{ex:deformation-Z2-counterexample}
{For $G=\mathbb Z^2$, we have}
\[
 \pi_0K^{\mathrm{def}}(\mathbb Z^2)\cong\mathbb Z,
 \qquad
 \bigl(\pi_0{\cR}_{C^*_{\max}(\mathbb Z^2)}(\C)\bigr)_{\mathrm{Hd}}
 \cong E_0(C(\mathbb T^2),\C)\cong\mathbb Z^2.
\]
The first  follows from connectedness of the character space $\mathbb T^2$, and the second from Theorem~\ref{thm:overview-qr-comparison} and
{Corollary~\ref{cor:overview-qr-E}}, together with the usual
$K$-homology calculation for $\mathbb T^2$.
Thus $\gamma_{\mathbb Z^2}$ is not surjective on components,
even after passing to the Hausdorff quotient of the target.
\hB
\end{ex}
Related failures of genuine representations to recover
$K$-cohomology through the topological Atiyah--Segal map
are discussed in
\cite[Introduction, Theorem~5.1 and Section~7]{ramras-surfaces}
and \cite[Proposition~8.8]{ramras-atiyah-segal}.
For characteristic classes of the bundles associated to almost
representations at prescribed finite controls, see
\cite{dadarlat-glebe-chern}.
\subsection{Controlled $KK$-theory and Milnor sequences}
\label{sec:controlled-KK}

The additional left-adjoint in \[
{\hat y:\EE_{\sepa}\leftrightarrows \Ind_{\aleph_{1}}(\EE_{\sepa}):\colim}
\]
{witnesses the fact that the category $\EE_{\sepa}$ is compactly assembled. Let $D$ and $B$ be in $\nCalg_{\sepa}$ and choose a presentation $\hat y(e_{\sepa}(D))\simeq\colim_n y(V_n)$ with $V_n$ in $\EE_{\sepa}$. We use the pro-abelian groups}
\begin{equation}\label{eq:shape-pro-mapping-groups}
{E_i(\hat y(e_{\sepa}(D)),B):=
\bigl(\pi_i\Hom_{\EE_{\sepa}}(V_n,e_{\sepa}(B))\bigr)_n
\qquad(i=0,1).}
\end{equation}
{Applying $\pi_0$ to the resulting limit of mapping anima gives a Milnor sequence natural in both variables}
\begin{equation}\label{gewrferfwf}
{0\to {\lim\nolimits^{1}}E_{1}(\hat y(e_{\sepa}(D)),B)\to E_0(D,B)\to \lim E_{0}(\hat y(e_{\sepa}(D)),B)\to 0,}
\end{equation}
{which is determined by the compactly assembled structure of $\EE_{\sepa}$.}

In the following  we explain that various Milnor sequences ocurring in the $C^{*}$-algebra literature
are equivalent to \eqref{gewrferfwf}. 
 
 {For a  shape system $(D_n)$ for $S^2D\otimes\bbK$, by \cite[Proposition~3.77(2)]{budu} and Bott periodicity, we have $\hat y(e_{\sepa}(D))\simeq\colim_n y(e_{\sepa}(D_n))$. Consequently,     \eqref{gewrferfwf} is isomorphic  to the Milnor sequence}
\begin{equation}\label{jtzhtrzhrtzj}0\to {\lim\nolimits_{n}^{1}} E_{1}(D_{n},B)\to E_{0}(D,B)\to \lim_{n} E_{0}(D_{n},B)\to 0.
\end{equation}
   
Willett--Yu define controlled $KK$-groups using projections
that approximately commute with a representation of the source
algebra, together with equivalence relations involving
prescribed controls. Their analysis of these cycles and
relations as the controls are refined gives a Milnor exact
sequence  for their controlled $KK$-groups \cite[Theorems~1.1--1.2]{willett-yu-controlled}.

If the source algebra belongs to the  UCT-class, Willett--Yu identify their sequence with Schochet's Milnor sequence \cite[Section~1.2]{willett-yu-controlled}. The latter has the following homotopical interpretation which we can safely formulate in terms of $E$-theory since the source algebras are all nuclear.
Let $D$ be a separable nuclear algebra satisfying the UCT, and let $B$ be separable. A $KK$-filtration {of $D$} is an increasing sequence $(A_{n})$ in $\nCalg_{\sepa}$ such that    $A_n^u\cong C(X_n)$ for finite CW complexes $X_n$, together with  an equivalence  $\colim_n e_{\sepa}(A_n)\simeq {e_{\sepa}(D)}$. {Such a filtration exists by Schochet\textquotesingle s theorem; see \cite[Section~1.2]{willett-yu-controlled}.} We get $E(D,B)\simeq\lim_n E(A_n,B)$, and  applying 
 $\pi_0$ gives Schochet's Milnor sequence
\begin{equation}\label{hrtegtrhbeth}0\to {\lim\nolimits_{n}^{1}} E_{1}(A_{n},B)\to E_{0}(D,B)\to \lim_{n} E_{0}(A_{n},B)\to 0\ .
\end{equation}{The filtration $(A_n)$ need not itself be a shape system for $D$. Since $X_n$ is a finite CW complex, $e_{\sepa}(C(X_n))$ belongs to the thick subcategory generated by the compact object $e_{\sepa}(\C)$. The split scalar projection $A_n^u\to\C$ exhibits $e_{\sepa}(A_n)$ as a retract of $e_{\sepa}(C(X_n))$, so $e_{\sepa}(A_n)$ is compact as well. We therefore obtain}
${\hat y(e_{\sepa}(D))\simeq \colim_n y(e_{\sepa}(A_n))}$. This implies that  Schochet's Milnor sequence \eqref{hrtegtrhbeth} is isomorphic to \eqref{gewrferfwf}.
%
%
  

For arbitrary separable $D$ and $B$, Carri\'on in \cite[Theorem~1.1]{carrion-milnor} gives another Milnor-sequence  description of  $E_0(D,B)$.   For a shape system   $(C_n)$  for $SD\otimes\bbK$ Carri\'on's Milnor sequence is, in the notation of the present paper,
\begin{equation}\label{eq:carrion-milnor} 
 0\longrightarrow{\lim\nolimits_n^1} \pi_{0}{\cM}_{C_n}(S^2B)
 \longrightarrow E_0(D,B)
 \longrightarrow\lim_n\pi_{0} {\cM}_{C_n}(SB)
 \longrightarrow0.
\end{equation}
For nuclear $D$, he compares it with the Willett--Yu sequence \cite[Proposition~4.5]{carrion-milnor}. 
Carri\'on constructs \eqref{eq:carrion-milnor} using asymptotic morphisms. However, it also has a standard homotopy-theoretic construction. The mapping-anima description of $E$-theory in \cite[Sections~3.5--3.6]{budu} gives an equivalence of commutative groups in anima
\begin{equation}\label{eq:E-suspended-mapping-limit}
 \Omega^\infty E(D,B)\simeq\lim_n {\cM}_{C_n}(SB).
\end{equation}
{Since ${\cM}_{C_n}(S^2B)\simeq\Omega {\cM}_{C_n}(SB)$, applying $\pi_0$ and the standard Milnor sequence to \eqref{eq:E-suspended-mapping-limit} gives \eqref{eq:carrion-milnor}, which by \cite[Sections~3.5--3.6]{budu} is isomorphic to \eqref{gewrferfwf}.} 


The isomorphisms between the Milnor sequences above show that they all induce the shape topology 
 on $E_{0}(D,B)$  which is, by definition, the one induced by  \eqref{gewrferfwf}, and which by Proposition~\ref{kophertgertr} coincides  with  the condensed topology.
 We get the formulas
\begin{equation}\label{gwergwerfwerg}
{EL(D,B):=\bigl(E_0(D,B)_{\mathrm{cond}}\bigl)_{\mathrm{Hd}} \cong \bigl(E_0(D,B)_{\mathrm{shape}}\bigl)_{\mathrm{Hd}}\cong \lim_{n}E_{0}(D_{n},B)_{{\mathrm{discrete}}}\cong 
 \lim_n\pi_0{\cM}_{C_n}(SB)}
\end{equation}
for the Hausdorff quotient.

The transformation $\alpha_{-,B}$ induces   a natural map  from the Milnor sequence  \eqref{gwegerfrewfw} for $\pi_{0}{\cR}_{D}(B)$ in Remark \ref{kopgrethertherthge} 
 to the sequence \eqref{gewrferfwf}.  
The difference between  \eqref{gwegerfrewfw}  and  \eqref{eq:carrion-milnor}  
 is the passage to an unsuspended representation construction. Formula~\eqref{eq:E-suspended-mapping-limit} uses a shape system for $SD\otimes\bbK$ and the suspended coefficient algebra $SB$; its mapping monoids are already grouplike. In contrast, ${\cR}_D(B)$ uses a shape system for $D$ and applies group completion and the reduced construction at its stages.

\subsection{Group completion and bivariant $K$-theory}

Group completion is a basic step in categorical constructions of $K$-theory \cite{segal-cohomology,gepner-groth-nikolaus}. In this paper we apply it to mapping monoids with block sum. We then compare its reduced, shape-resolved form ${\cR}_D(B)$ with the mapping groups in anima of the universal homological category $E_{\sepa}$.


In algebraic $K$-theory of stable {$\infty$-categories},
Efimov gives a related description of mapping spectra in the category $\operatorname{Mot}^{\mathrm{loc}}$ of localizing noncommutative motives \cite[Theorem~0.4]{efimov-rigidity}. Let $U_{\mathrm{loc}}$ denote the universal filtered-colimit-preserving localizing invariant of small idempotent-complete stable $\infty$-categories, and let $K$ denote nonconnective algebraic $K$-theory. For such a category $\mathcal C$, its countable Calkin category is
\begin{equation}\label{eq:efimov-calkin-category}
 \operatorname{Calk}_{\omega_1}(\mathcal C)
 :=\bigl(\operatorname{Ind}(\mathcal C)^{\omega_1}/\mathcal C\bigr)^{\mathrm{idem}},
\end{equation}
where the superscript $\omega_1$ selects the $\omega_1$-compact objects and $(-)^{\mathrm{idem}}$ denotes idempotent completion. Write $\operatorname{Calk}_{\omega_1}^2$ for its second iterate. Let $\mathcal B$ and $\mathcal C$ be small idempotent-complete stable $\infty$-categories, and assume that $\mathcal B$ is $\omega_1$-compact in $\operatorname{Cat}^{\mathrm{perf}}$. {We write $\Fun^{\mathrm{ex}}$ for the stable $\infty$-category of exact functors.} Efimov's equivalence is
\begin{equation}\label{eq:efimov-motivic-mapping-spectrum}
 \operatorname{map}_{\operatorname{Mot}^{\mathrm{loc}}}
 \bigl(U_{\mathrm{loc}}(\mathcal B),U_{\mathrm{loc}}(\mathcal C)\bigr)
 \simeq\Omega^2K\bigl(\Fun^{\mathrm{ex}}(\mathcal B,
                  \operatorname{Calk}_{\omega_1}^2(\mathcal C))\bigr).
\end{equation}
Applying $\Omega^\infty$ gives the corresponding description of mapping anima.

The common idea is to obtain a mapping object in a universal category by applying a $K$-theoretic construction to morphism data. In \eqref{eq:efimov-motivic-mapping-spectrum}, these data form a stable category of exact functors, and algebraic $K$-theory imposes the relations coming from cofiber sequences. Our input instead consists of stabilized $C^*$-algebra homomorphisms with block sum. We apply group completion, the reduced construction, and shape resolution. {For $\alpha$-admissible sources in the sense of Definition~\ref{def:alpha-admissible}, this construction recovers the $E$-theory mapping groups in anima; this includes $D\in\cU_{\mathrm{nqs}}$ by Corollary~\ref{cor:U-nqs-remainder}. For general $D\in\cU_{\mathrm{nq}}$, Theorems~\ref{thm:overview-one-splitting} and~\ref{thm:overview-one-hausdorff} give a splitting of the comparison and an isomorphism on Hausdorff component groups.} 


\section*{Acknowledgements}

The author acknowledges support by the Deutsche
Forschungsgemeinschaft (DFG, German Research Foundation)
through SFB~1785 \emph{Generalised Motivic Methods in Geometry},
project number~560334683.

\section*{Use of AI tools}

The author used ChatGPT and Codex (OpenAI) for mathematical
discussions, literature searches, and assistance with drafting
and revising the manuscript. The author takes responsibility
for the mathematical arguments, references, and final text.


\begin{thebibliography}{BD26}
\bibitem{bekhsy} T.~Bachmann, E.~Elmanto, M.~Hoyois, A.~A.~Khan,
V.~Sosnilo and M.~Yakerson,
\emph{On the infinite loop spaces of algebraic cobordism and the motivic sphere},
\'Epijournal de G\'eom\'etrie Alg\'ebrique \textbf{5} (2021), article~9.
\href{https://epiga.episciences.org/7494/pdf}{Full text}.

\bibitem{barlak-raum-dimension-drop}
{S.~Barlak and S.~Raum,
\emph{Cartan subalgebras in dimension drop algebras},
Journal of the Institute of Mathematics of Jussieu \textbf{20} (2021), no.~3, 725--755.
\href{https://arxiv.org/abs/1712.01957}{arXiv:1712.01957}.}

\bibitem{bartle-graves}
{R.~G.~Bartle and L.~M.~Graves,
\emph{Mappings between function spaces},
Transactions of the American Mathematical Society
\textbf{72} (1952), 400--414.
\href{https://doi.org/10.1090/S0002-9947-1952-0047910-X}{doi:10.1090/S0002-9947-1952-0047910-X}.}

\bibitem{beltita-beltita-ludwig-fourier}
{I.~Belti\c{t}\u{a}, D.~Belti\c{t}\u{a} and J.~Ludwig,
\emph{Fourier transforms of $C^*$-algebras of nilpotent Lie groups},
\href{https://arxiv.org/abs/1411.3254}{arXiv:1411.3254}, 2014, revised 2015.}

\bibitem{blackadar-shape} B.~Blackadar,
\emph{Shape theory for $C^*$-algebras},
Mathematica Scandinavica \textbf{56} (1985), no.~2, 249--275.
\href{https://doi.org/10.7146/math.scand.a-12100}{doi:10.7146/math.scand.a-12100}.

\bibitem{blackadar-kbook}
B.~Blackadar,
\emph{$K$-Theory for Operator Algebras},
second edition, Mathematical Sciences Research Institute Publications
\textbf{5}, Cambridge University Press, Cambridge, 1998.
\href{https://bruceblackadar.com/Mathematics/book6.pdf}{Author's online edition}.

\bibitem{blackadar-operator-algebras}
{B.~Blackadar,
\emph{Operator Algebras: Theory of $C^*$-Algebras and von Neumann Algebras},
Encyclopaedia of Mathematical Sciences \textbf{122},
Springer, Berlin, 2006.
\href{https://bruceblackadar.com/Mathematics/Cycr.pdf}{Author's online edition}.}

\bibitem{boca-free-products}
{F.~P.~Boca,
\emph{A note on full free product $C^*$-algebras, lifting and quasidiagonality},
in \emph{Proceedings of the 16th International Conference on Operator Theory (Timi\c{s}oara, 1996)},
The Theta Foundation, Bucharest, 1997, pp.~51--63.
\href{https://arxiv.org/abs/1508.03123}{arXiv:1508.03123}.}

\bibitem{brown-ozawa}
N.~P.~Brown and N.~Ozawa,
\emph{$C^*$-Algebras and Finite-Dimensional Approximations},
Graduate Studies in Mathematics \textbf{88},
American Mathematical Society, Providence, RI, 2008.
\href{https://doi.org/10.1090/gsm/088}{doi:10.1090/gsm/088}.

\bibitem{orbita} U.~Bunke,
\emph{$KK$- and $E$-theory via homotopy theory},
Orbita Mathematicae \textbf{1} (2024), no.~2, 103--210.
\href{https://doi.org/10.2140/om.2024.1.103}{doi:10.2140/om.2024.1.103};
\href{https://arxiv.org/abs/2304.12607}{arXiv:2304.12607}.

\bibitem{budu} U.~Bunke and B.~D\"unzinger,
\emph{$E$-theory is compactly assembled},
\href{https://arxiv.org/abs/2402.18228v4}{arXiv:2402.18228v4},
September 2026.

\bibitem{bunke-engel-land-kk}
U.~Bunke, A.~Engel and M.~Land,
\emph{A stable $\infty$-category for equivariant $KK$-theory},
\href{https://arxiv.org/abs/2102.13372}{arXiv:2102.13372}.

\bibitem{bunke-joachim-stolz-k}
{U.~Bunke, M.~Joachim and S.~Stolz,
\emph{Classifying spaces and spectra representing the $K$-theory of a graded $C^{*}$-algebra},
in \emph{High-dimensional manifold topology},
World Scientific, River Edge, NJ, 2003, 80--102.
\href{https://doi.org/10.1142/9789812704443_0003}{doi:10.1142/9789812704443\_0003}.}

\bibitem{carrion-milnor} {J.~R.~Carri\'on,
\emph{A Milnor exact sequence for $E$-theory},
\href{https://arxiv.org/abs/2608.06496}{arXiv:2608.06496}, 2026.}

\bibitem{carrion-schafhauser-topology} J.~R.~Carri\'on and C.~Schafhauser,
\emph{A topology on $E$-theory},
Journal of the London Mathematical Society {\textbf{109} (2024), no.~6,} article e12917.
\href{https://doi.org/10.1112/jlms.12917}{doi:10.1112/jlms.12917};
\href{https://arxiv.org/abs/2306.13757}{arXiv:2306.13757}.

\bibitem{choi-effros-lifting}
{M.-D.~Choi and E.~G.~Effros,
\emph{The completely positive lifting problem for $C^*$-algebras},
Annals of Mathematics (2) \textbf{104} (1976), no.~3, 585--609.
\href{https://doi.org/10.2307/1970968}{doi:10.2307/1970968}.}

\bibitem{clausen-scholze-analytic-stacks}
{D.~Clausen and P.~Scholze,
\emph{Analytic Stacks},
lecture course at IH\'ES and MPIM Bonn, winter semester 2023/24.
\href{https://people.mpim-bonn.mpg.de/scholze/AnalyticStacks.html}{Course webpage};
\href{https://www.youtube.com/watch?v=_4G582SIo28}{Lecture~2: Light condensed sets}.}

\bibitem{connes-higson-asymptotic}
{A.~Connes and N.~Higson, \emph{D\'eformations, morphismes
asymptotiques et $K$-th\'eorie bivariante},
C.~R.~Acad.~Sci.~Paris S\'er.~I Math.~\textbf{311} (1990), no.~2, 101--106.}

\bibitem{courtney-sherman-contraction}
{K.~Courtney and D.~Sherman,
{\emph{The universal $C^{*}$-algebra of a contraction}},
Journal of Operator Theory \textbf{84} (2020), no.~1, 153--184.
\href{https://arxiv.org/abs/1811.04043}{arXiv:1811.04043}.}

\bibitem{cuntz-k} J.~Cuntz,
\emph{$K$-theory for certain $C^*$-algebras},
Annals of Mathematics (2) \textbf{113} (1981), no.~1, 181--197.
\href{https://doi.org/10.2307/1971137}{doi:10.2307/1971137}.

\bibitem{dadarlat-shape} M.~Dadarlat,
\emph{Shape theory and asymptotic morphisms for $C^*$-algebras},
Duke Mathematical Journal \textbf{73} (1994), no.~3, 687--711.

\bibitem{dadarlat-topology} M.~Dadarlat,
\emph{On the topology of the Kasparov groups and its applications},
Journal of Functional Analysis \textbf{228} (2005), no.~2, 394--418.
\href{https://arxiv.org/abs/math/0410400}{arXiv:math/0410400}.

\bibitem{dadarlat-group-quasirepresentations}
{M.~Dadarlat,
\emph{Group quasi-representations and almost flat bundles},
Journal of Noncommutative Geometry \textbf{8} (2014), no.~1, 163--178.
\href{https://doi.org/10.4171/JNCG/152}{doi:10.4171/JNCG/152}.}

\bibitem{dadarlat-eilers-asymptotic} M.~Dadarlat and S.~Eilers,
\emph{Asymptotic unitary equivalence in $KK$-theory},
$K$-Theory \textbf{23} (2001), 305--322.
\href{https://www.math.purdue.edu/~mdd/Publications/kk-asympt-with-eilers.pdf}{Full text}.

\bibitem{dadarlat-glebe-chern} M.~Dadarlat and F.~Glebe,
\emph{On Chern classes of almost representations},
\href{https://arxiv.org/abs/2509.23950}{arXiv:2509.23950}, 2025.

\bibitem{dellambrogio-emerson-kandelaki-meyer}
{I.~Dell'Ambrogio, H.~Emerson, T.~Kandelaki and R.~Meyer,
\emph{A functorial equivariant $K$-theory spectrum and an equivariant Lefschetz formula},
preprint, 2011.
\href{https://arxiv.org/abs/1104.3441}{arXiv:1104.3441}.}

\bibitem{efimov-rigidity} {A.~I.~Efimov,
\emph{Rigidity of the category of localizing motives},
\href{https://arxiv.org/abs/2510.17010}{arXiv:2510.17010}, 2025.}

\bibitem{eilers-loring-pedersen-nccw}
{S.~Eilers, T.~A.~Loring and G.~K.~Pedersen,
\emph{Stability of anticommutation relations: An application of noncommutative CW complexes},
Journal f\"ur die reine und angewandte Mathematik \textbf{499} (1998), 101--143.
\href{https://doi.org/10.1515/crll.1998.055}{doi:10.1515/crll.1998.055}.}

\bibitem{elliott-evans-rotation}
{G.~A.~Elliott and D.~E.~Evans,
\emph{The structure of the irrational rotation $C^*$-algebra},
Annals of Mathematics (2) \textbf{138} (1993), no.~3, 477--501.
\href{https://doi.org/10.2307/2946553}{doi:10.2307/2946553}.}

\bibitem{enders} D.~Enders,
\emph{Semiprojectivity for Kirchberg algebras},
\href{https://arxiv.org/abs/1507.06091}{arXiv:1507.06091}, 2015.

\bibitem{gabe-classification} J.~Gabe,
\emph{Classification of $\mathcal O_\infty$-stable $C^*$-algebras},
Memoirs of the American Mathematical Society \textbf{293} (2024),
no.~1461, v+115~pp.
\href{https://doi.org/10.1090/memo/1461}{doi:10.1090/memo/1461};
\href{https://arxiv.org/abs/1910.06504}{arXiv:1910.06504}.

\bibitem{gabe-szabo} J.~Gabe and G.~Szab\'o,
\emph{The stable uniqueness theorem for equivariant Kasparov theory},
American Journal of Mathematics \textbf{147} (2025), no.~6, 1527--1576.
\href{https://doi.org/10.1353/ajm.2025.a975703}{doi:10.1353/ajm.2025.a975703};
\href{https://arxiv.org/abs/2202.09809}{arXiv:2202.09809}.

\bibitem{gepner-groth-nikolaus}
D.~Gepner, M.~Groth and T.~Nikolaus,
\emph{Universality of multiplicative infinite loop space machines},
Algebraic \& Geometric Topology \textbf{15} (2015), no.~6, 3107--3153.
\href{https://doi.org/10.2140/agt.2015.15.3107}{doi:10.2140/agt.2015.15.3107};
\href{https://arxiv.org/abs/1305.4550}{arXiv:1305.4550}.

\bibitem{gepner-haugseng-enriched}
{D.~Gepner and R.~Haugseng,
\emph{Enriched $\infty$-categories via non-symmetric
$\infty$-operads},
Advances in Mathematics \textbf{279} (2015), 575--716.
\href{https://arxiv.org/abs/1312.3178}{arXiv:1312.3178}.}

\bibitem{green-imprimitivity}
{P.~Green,
\emph{The local structure of twisted covariance algebras},
Acta Mathematica \textbf{140} (1978), 191--250.
\href{https://doi.org/10.1007/BF02392308}{doi:10.1007/BF02392308}.}

\bibitem{halmos-ten-problems}
{P.~R.~Halmos,
\emph{Ten problems in Hilbert space},
Bulletin of the American Mathematical Society \textbf{76} (1970), 887--933.
\href{https://doi.org/10.1090/S0002-9904-1970-12502-2}{doi:10.1090/S0002-9904-1970-12502-2}.}

\bibitem{harpaz-spans}
{Y.~Harpaz,
\emph{Ambidexterity and the universality of finite spans},
Proceedings of the London Mathematical Society (3)
\textbf{121} (2020), 1121--1170.
\href{https://doi.org/10.1112/plms.12367}{doi:10.1112/plms.12367}.}

\bibitem{higson-excision}
{N.~Higson, \emph{Categories of fractions and excision in $KK$-theory},
J.~Pure Appl.~Algebra \textbf{65} (1990), no.~2, 119--138.
\href{https://doi.org/10.1016/0022-4049(90)90114-W}{doi:10.1016/0022-4049(90)90114-W}.}

\bibitem{higson-kasparov} N.~Higson and G.~Kasparov,
\emph{$E$-theory and $KK$-theory for groups which act properly
and isometrically on Hilbert space},
Invent.\ Math.\ \textbf{144} (2001), 23--74.

\bibitem{javier-quasidiagonal}
{T.~F.~Javier,
\emph{Quasi-diagonal $C^*$-algebras},
Journal of Functional Analysis \textbf{25} (1977), 50--57.
\href{https://doi.org/10.1016/0022-1236(77)90031-3}{doi:10.1016/0022-1236(77)90031-3}.}

\bibitem{jiang-su}
{X.~Jiang and H.~Su,
\emph{On a simple unital projectionless $C^*$-algebra},
American Journal of Mathematics \textbf{121} (1999), no.~2, 359--413.
\href{https://doi.org/10.1353/ajm.1999.0012}{doi:10.1353/ajm.1999.0012}.}

\bibitem{joachim-ko}
{M.~Joachim,
\emph{A symmetric ring spectrum representing $KO$-theory},
Topology \textbf{40} (2001), no.~2, 299--308.
\href{https://doi.org/10.1016/S0040-9383(99)00063-4}{doi:10.1016/S0040-9383(99)00063-4}.}

\bibitem{kasparov-absorption} G.~G.~Kasparov,
\emph{Hilbert $C^*$-modules: theorems of Stinespring and Voiculescu},
Journal of Operator Theory \textbf{4} (1980), 133--150.
\href{https://www.theta.ro/jot/archive/1980-004-001/1980-004-001-007.pdf}{Full text}.

\bibitem{lance-nuclear}
{E.~C.~Lance,
\emph{On nuclear $C^*$-algebras},
J.~Funct.~Anal.~\textbf{12} (1973), 157--176.}

\bibitem{lawson-deformation} T.~Lawson,
\emph{The product formula in unitary deformation $K$-theory},
$K$-Theory \textbf{37} (2006), no.~4, 395--422.
\href{https://doi.org/10.1007/s10977-006-9003-9}{doi:10.1007/s10977-006-9003-9};
\href{https://arxiv.org/abs/math/0503468}{arXiv:math/0503468}.

\bibitem{loring-lifting}
T.~A.~Loring,
\emph{Lifting Solutions to Perturbing Problems in $C^*$-Algebras},
Fields Institute Monographs \textbf{8},
American Mathematical Society, Providence, RI, 1997.
\href{https://math.unm.edu/~loring/book/homepage.html}{Author's book page}.

\bibitem{may-einfty}
J.~P.~May, with contributions by F.~Quinn, N.~Ray and J.~Tornehave,
\emph{$E_\infty$ Ring Spaces and $E_\infty$ Ring Spectra},
Lecture Notes in Mathematics \textbf{577},
Springer-Verlag, Berlin--New York, 1977.
\href{https://www.math.uchicago.edu/~may/BOOKS/e_infty.pdf}{Author's online edition}.

\bibitem{meyer-nest-localisation}
{R.~Meyer and R.~Nest,
\emph{The Baum--Connes conjecture via localisation of categories},
Topology \textbf{45} (2006), no.~2, 209--259.
\href{https://doi.org/10.1016/j.top.2005.07.001}{doi:10.1016/j.top.2005.07.001};
\href{https://arxiv.org/abs/math/0312292}{arXiv:math/0312292}.}

\bibitem{phillips-sorensen-thiel}
N.~C.~Phillips, A.~P.~W.~S{\o}rensen and H.~Thiel,
\emph{Semiprojectivity with and without a group action},
Journal of Functional Analysis \textbf{268} (2015), 929--973.
\href{https://doi.org/10.1016/j.jfa.2014.11.005}{doi:10.1016/j.jfa.2014.11.005};
\href{https://arxiv.org/abs/1403.3440}{arXiv:1403.3440}.

\bibitem{pimsner-voiculescu-rotation}
{M.~Pimsner and D.~Voiculescu,
\emph{Imbedding the irrational rotation $C^*$-algebra into an AF-algebra},
J.~Operator Theory \textbf{4} (1980), 201--210.}

\bibitem{ramras-deformation-excision} D.~A.~Ramras,
\emph{Excision for deformation $K$-theory of free products},
{Algebraic \& Geometric Topology \textbf{7} (2007), 2239--2270.
\href{https://doi.org/10.2140/agt.2007.7.2239}{doi:10.2140/agt.2007.7.2239};}
\href{https://arxiv.org/abs/math/0703463}{arXiv:math/0703463}.

\bibitem{ramras-surfaces} D.~A.~Ramras,
\emph{Yang--Mills theory over surfaces and the Atiyah--Segal theorem},
Algebraic \& Geometric Topology \textbf{8} (2008), no.~4, 2209--2251.
\href{https://doi.org/10.2140/agt.2008.8.2209}{doi:10.2140/agt.2008.8.2209};
\href{https://arxiv.org/abs/0710.0681}{arXiv:0710.0681}.

\bibitem{ramras-atiyah-segal} D.~A.~Ramras,
\emph{The topological Atiyah--Segal map},
\href{https://arxiv.org/abs/1607.06430}{arXiv:1607.06430}.

\bibitem{roydor-non-llp}
{J.~Roydor,
\emph{$C^*(\mathbb F_2\times\mathbb F_2)$ fails the local lifting property},
\href{https://arxiv.org/abs/2611.31147}{arXiv:2611.31147}, 2026.}

\bibitem{schafhauser-tikuisis-white-problems}
{C.~Schafhauser, A.~Tikuisis and S.~White,
\emph{Nuclear $C^*$-algebras: 99 problems},
\href{https://arxiv.org/abs/2506.10902v2}{arXiv:2506.10902v2}, 2026.}

\bibitem{scholze-condensed}
{P.~Scholze,
\emph{Lectures on Condensed Mathematics},
lecture notes for the course at the University of Bonn, summer 2019, revised 2026.
\href{https://arxiv.org/abs/2605.03658}{arXiv:2605.03658}.}

\bibitem{segal-cohomology}
G.~Segal,
\emph{Categories and cohomology theories},
Topology \textbf{13} (1974), no.~3, 293--312.
\href{https://doi.org/10.1016/0040-9383(74)90022-6}{doi:10.1016/0040-9383(74)90022-6}.

\bibitem{shulman} T.~Shulman,
\emph{Homotopy lifting, asymptotic homomorphisms, and traces},
\href{https://arxiv.org/abs/2508.00125}{arXiv:2508.00125}, version~5, 2026.

\bibitem{sorensen-thiel} A.~P.~W.~S{\o}rensen and H.~Thiel,
\emph{A characterization of semiprojectivity for commutative $C^*$-algebras},
Proceedings of the London Mathematical Society (3) \textbf{105} (2012),
no.~5, 1021--1046.
\href{https://doi.org/10.1112/plms/pdr051}{doi:10.1112/plms/pdr051};
\href{https://arxiv.org/abs/1101.1856}{arXiv:1101.1856}.

\bibitem{skandalis-k-nuclearity}
{G.~Skandalis,
\emph{Une notion de nucl\'earit\'e en $K$-th\'eorie (d'apr\`es J.~Cuntz)},
$K$-Theory \textbf{1} (1988), no.~6, 549--573.
\href{https://doi.org/10.1007/BF00533786}{doi:10.1007/BF00533786}.}

\bibitem{szabo-uniqueness} G.~Szab\'o,
\emph{The uniqueness theorem for Kasparov theory},
\href{https://arxiv.org/abs/2601.23029}{arXiv:2601.23029}, 2026.

\bibitem{thiel-semiprojective-limits}
{H.~Thiel,
\emph{Inductive limits of semiprojective $C^*$-algebras},
Advances in Mathematics \textbf{347} (2019), 597--618.
\href{https://doi.org/10.1016/j.aim.2019.02.030}{doi:10.1016/j.aim.2019.02.030};
\href{https://arxiv.org/abs/1804.04908}{arXiv:1804.04908}.}

\bibitem{thomsen-absorbing} {K.~Thomsen,
\emph{On absorbing extensions},
Proceedings of the American Mathematical Society \textbf{129} (2001), no.~5, 1409--1417.
\href{https://doi.org/10.1090/S0002-9939-00-05936-0}{doi:10.1090/S0002-9939-00-05936-0}.}

\bibitem{weibel-kbook}
C.~A.~Weibel,
\emph{The $K$-book: An Introduction to Algebraic $K$-Theory},
Graduate Studies in Mathematics \textbf{145},
American Mathematical Society, Providence, RI, 2013.
\href{https://sites.math.rutgers.edu/~weibel/Kbook.html}{Author's online edition}.

\bibitem{willett} R.~Willett,
\emph{An Atiyah--Segal homomorphism for approximate representations of
infinite groups, and almost multiplicative ucp maps},
\href{https://arxiv.org/abs/2610.24252}{arXiv:2610.24252}, October 2026.

\bibitem{willett-yu-uct} R.~Willett and G.~Yu,
\emph{The Universal Coefficient Theorem for $C^*$-Algebras with
Finite Complexity},
Memoirs of the European Mathematical Society \textbf{8},
EMS Press, 2024.
\href{https://doi.org/10.4171/MEMS/8}{doi:10.4171/MEMS/8};
\href{https://arxiv.org/abs/2104.10766}{arXiv:2104.10766}.

\bibitem{willett-yu-controlled} R.~Willett and G.~Yu,
\emph{Controlled $KK$-theory and a Milnor exact sequence},
Documenta Mathematica \textbf{31} (2026), no.~3, 503--581.
\href{https://doi.org/10.4171/DM/1056}{doi:10.4171/DM/1056};
\href{https://arxiv.org/abs/2011.10906}{arXiv:2011.10906}.

\end{thebibliography}
\end{document}